\documentclass[11pt,a4paper]{article}
\usepackage{amsmath}
\usepackage{amsfonts}
\usepackage{amssymb}
\usepackage{amsthm}
\usepackage{mathrsfs}
\usepackage{dsfont}
\usepackage{paralist}
\usepackage{natbib}
\usepackage{bm}
\usepackage{color}
\usepackage[colorlinks=true,linkcolor=blue,citecolor=blue,pdfborder={0 0 0}]{hyperref}
\usepackage{graphicx}
\usepackage{tabularx}
\usepackage[left=2.5cm,right=2.5cm,top=2cm,bottom=2cm,includeheadfoot]{geometry}
\usepackage{comment}
\usepackage{setspace}

\newcount\Comments  
\newcommand{\kibitz}[2]{\ifnum\Comments=1\textcolor{#1}{#2}\fi}

\newtheoremstyle{normal}
{2ex}               
{3ex}               
{}                  
{}                  
{\bfseries} 
{}                  
{2pt}   
{\thmname{#1}\thmnumber{ #2.} \thmnote{(#3)}}

\newtheoremstyle{italic}
{2ex}
{3ex}
{\itshape}
{}
{\bfseries} 
{}
{2pt}
{\thmname{#1}\thmnumber{ #2.} \thmnote{(#3)}}

\theoremstyle{normal}
\newtheorem{definition}{Definition}[section]
\newtheorem{remark}[definition]{Remark}
\newtheorem{example}[definition]{Example}
\newtheorem{condition}[definition]{Assumption}

\theoremstyle{italic}
\newtheorem{theorem}[definition]{Theorem}
\newtheorem{lemma}[definition]{Lemma}
\newtheorem{proposition}[definition]{Proposition}
\newtheorem{corollary}[definition]{Corollary}

\newcommand\N{\mathbb{N}}

\newcommand\R{\mathbb{R}}

\newcommand\eps{\varepsilon}
\newcommand\Prob{\mathbb{P}}    

\newcommand\Dc{\mathcal{D}}
\newcommand\Fc{\mathcal{F}}
\newcommand\Gc{\mathcal{G}}

\newcommand\Bb{\mathbb{B}}
\newcommand\Cb{\mathbb{C}}
\newcommand\Db{\mathbb{D}}
\newcommand\Eb{\mathbb{E}}

\newcommand\Gb{\mathbb{G}}
\newcommand\Hb{\mathbb{H}}

\newcommand\linaeps{\ell^{\infty}(A_\eps)}
\newcommand\linbeps{\ell^{\infty}(B_\eps)}

\newcommand\izba{\mathcal I(\bar z)}
\newcommand\izei{\mathcal I(z_1)}
\newcommand\izzw{\mathcal I(z_2)}

\newcommand\linne{\ell^{\infty}([0,1])}

\newcommand\iz{\mathcal I(z)}

\newcommand\Ron{\mathbb R \setminus \lbrace 0 \rbrace}
\newcommand\Xn{X^{(n)}}
\newcommand\neps{(n, \varepsilon)}

\newcommand\kseqi{\zeta}
\newcommand\gseqi{\theta}
\newcommand\taili{z}
\newcommand\smooi{\varpi}
\newcommand\thrle{\varkappa}
\newcommand\preest{\hat \theta_{n}}
\newcommand\consest{\hat \gseqi^{(\eps)}_n}
\newcommand\chapo{\theta_0^{(\varepsilon)}}

\newcommand\hathrleei{\hat \thrle_{n,B}^{(\eps, \alpha)}(1)}
\newcommand\hathrle{\hat \thrle_{n,B}^{(\eps, \alpha)}(r)}

\newcommand\hacothrle{\hat \thrle_{\scriptscriptstyle n,B_n}^{(\eps, \alpha_n)}(r)}
\newcommand\hacothrleor{\hat \thrle_{\scriptscriptstyle n,B_n}^{(\eps, \alpha_n)}(r)}

\newcommand\bootsuhb{\hat \Hb^{(\eps)}}

\newcommand\emprbodif{K_{n,B}^{(\eps,r)}}
\newcommand\emprboindif{K_{n,B}^{(\eps,r)-}}

\newcommand\weak{\rightsquigarrow}
\newcommand\weakP{{\, {\weak_\xi}\ }}
\newcommand\pto{\stackrel{\scriptscriptstyle \mathbb P^*}{\rightarrow}}
\newcommand\probto{\stackrel{\scriptscriptstyle \mathbb P}{\rightarrow}}
\newcommand\defeq{:=}
\newcommand{\ip}[1]{\lfloor #1 \rfloor}

\allowdisplaybreaks[1]

\begin{document}
\title{Nonparametric inference of gradual changes  in the jump behaviour of time-continuous processes
}

\author{Michael Hoffmann\footnotemark[1]\  ,  Mathias Vetter\footnotemark[2] ~and Holger Dette\footnotemark[1]\ , \bigskip \\
{ Ruhr-Universit\"at Bochum \ \& \ Christian-Albrechts-Universit\"at zu Kiel}
}

\footnotetext[1]{Ruhr-Universit\"at Bochum,
Fakult\"at f\"ur Mathematik, 44780 Bochum, Germany.
{E-mail:} holger.dette@rub.de, michael.hoffmann@rub.de}

\footnotetext[2]{Christian-Albrechts-Universit\"at zu Kiel,
Mathematisches Seminar, Ludewig-Meyn-Str. 4, 24118 Kiel, Germany.
{E-mail:} vetter@math.uni-kiel.de}

\maketitle

\begin{abstract}

In applications the properties of a stochastic feature often change gradually rather than abruptly, that is: after a constant phase for some time they slowly start to vary. In this paper we discuss statistical inference for the detection and the localisation of gradual changes in the jump characteristic of a discretely observed Ito semimartingale. We propose a new measure of time variation for the jump behaviour of the process. The statistical uncertainty of a corresponding  estimate is analyzed by deriving new results on the weak convergence of a sequential empirical tail integral process and a corresponding multiplier bootstrap procedure.
\end{abstract}


\smallskip


\section{Introduction}
\def\theequation{1.\arabic{equation}}
\setcounter{equation}{0}

Stochastic processes in continuous time are widely used in the applied sciences nowadays, as they allow for a flexible modeling of the evolution of various real-life phenomena over time. Speaking of mathematical finance, of particular interest is the family of semimartingales, which is theoretically appealing as it satisfies a certain condition on the absence of arbitrage in financial markets and yet is rich enough to reproduce stylized facts from empirical finance such as volatility clustering, leverage effects or jumps. For this reason, the development of statistical tools modeled by discretely observed It\=o semimartingales has been a major topic over the last years, both regarding the estimation of crucial quantities used for model calibration purposes and with a view on tests to check whether a certain model fits the data well. 
For a detailed  overview of the state of the art we refer to  the recent monographs by \cite{JacPro12} and \cite{AitJac14}.

These statistical tools typically differ highly, depending on the quantities of interest. When the focus is on the volatility, most concepts are essentially concerned with discrete observations of the continuous martingale part. In this case one is naturally close to the Gaussian framework, and so a lot of classical concepts from standard parametric statistics turn out to be powerful methods. The situation is different with a view on the jump behaviour of the process, mainly for two reasons: On one hand there is much more flexibility in the choice of the jump measure than there is regarding the diffusive part. On the other hand even if one restricts the model to certain parametric families the standard situation is the one of $\beta$-stable processes, $0 < \beta < 2$, for which the mathematical analysis is quite difficult, at least in comparison to Brownian motion. To mention recent work besides the afore-mentioned monographs, see for example \cite{Nic16} and \cite{HofVet15} on the estimation of the jump distribution function of a L\'evy process or \cite{Tod15} on the estimation of the jump activity index from high-frequency observations.  

In the following, we are interested in the evolution of the jump behaviour over time
 in a completely non-parametric setting where we assume only stuctural conditions on the characteristic triplet of the underlying It\=o semimartingale. To be precise, let $X=(X_t)_{t \geq 0}$ be an It\=o semimartingale with a decomposition
\begin{multline} \label{ItoSem}
X_t = X_0 +  \int_0^t b_s \, ds + \int_0^t \sigma_s\, dW_s + \int_0^t \int_{\R} u 1_{\{|u| \leq 1\}} (\mu-\bar \mu)(ds,du) \\
+ \int_0^t \int_{\R} u 1_{\{|u|>1\}} \mu(du,dz),
\end{multline}
where $W$ is a standard Brownian motion, $\mu$ is a Poisson random measure on $\R^+ \times \R$, and the predictable compensator $\bar \mu$ satisfies $\bar \mu(ds,du) = ds \: \nu_s(du)$. The main quantity of interest is the kernel $\nu_s$ which controls the number and the size of the jumps around time $s$. 

In \cite{BueHofVetDet15} the authors are interested in the detection of abrupt changes in the jump measure of $X$. Based on high-frequency observations $X_{i \Delta_n}$, $i=0,\ldots,n$, with $\Delta_n \to 0$ they construct a test for a constant $\nu$ against the alternative 
\[
\nu_t^{(n)} = \mathtt 1_{ \{ t < \ip{n\theta_0} \Delta_n \} } \nu_1 +  \mathtt 1_{\{t \ge \ip{n\theta_0} \Delta_n \} } \nu_2.
\]
Here the authors face a similar problem as in the classical situation of changes in the mean of a time series, namely that the ``change point'' $\theta_0$ can only be defined relative to the length of the covered time horizon $n \Delta_n$ which needs to tend to infinity. In general, this problem cannot be avoided as there are only finitely  many large jumps on every compact interval, so consistent estimators for the jump measure have to be constructed over the entire positive half-line. 

There are other types of  changes in the jump behaviour of a process than just abrupt ones, though. In the sequel, we will deal with gradual (smooth, continuous) changes of $\nu_s$ and discuss how and how well they can be detected. A similar problem has recently been addressed in \cite{Tod16} who constructs a test for changes in the activity index. Since this index is determined by the infinitely many small jumps around zero, such a test can be constructed over a day. On the other hand, estimation of an index is obviously a simpler problem than estimation of an entire measure. 

While the  problem of detecting abrupt changes has been discussed intensively in a time series context 
(see   \cite{auehor2013}  and \cite{jandhyala:2013} for a  review of the literature), 
 detecting gradual changes is a much harder problem and  the  methodology 
  is not so well developed. Most  authors consider  nonparametric  location or parametric models  with independently distributed  observations, 
 and we refer to   \cite{bissel1984a}, \cite{gan1991},  \cite{siezha1994}, \cite{huskova1999}, \cite{husste2002}  and \cite{Mallik2013} among others.
 See  also \cite{aueste2002} for some results in  a time series model and \cite{VogDet15} who developed a
  nonparametric method for the analysis of smooth changes in locally stationary time series.

  The present paper is devoted to the development of nonparametric inference for gradual changes in the jump properties
  of a  discretely observed It\=o semimartingale. In Section \ref{sec2} we introduce the formal setup as well as a measure of time variation which is used to identify changes in the jump characteristic later on. 
  Section \ref{sec3} is concerned with weak convergence of an estimator for this measure, and as a consequence we also obtain weak convergence of related statistics which can be used for testing for a gradual change  and for localizing the first change point. As the limiting distribution depends in a complicated way on the unknown jump characteristic, a bootstrap procedure is discussed which can be used for a data driven choice of a regularization parameter of the change point estimator and for approximating quantiles of the test statistics. Section \ref{sec4} contains the formal derivation of an
  estimator of the    change point and  a test for gradual changes.  Finally, a  brief  simulation study can be found in Section \ref{sec5},  while
  all proofs are relegated to Section \ref{sec:AuxRes}. 

\section{Preliminaries and a measure of gradual changes }
\def\theequation{2.\arabic{equation}}
\setcounter{equation}{0}
\label{sec2}

In the sequel let $X^{\scriptscriptstyle (n)}=(X_t^{\scriptscriptstyle (n)})_{t \geq 0}$ be an It\=o semimartingale of the form  \eqref{ItoSem} with characteristic triplet $(b_s^{\scriptscriptstyle (n)}, \sigma^{\scriptscriptstyle (n)}_s, \nu_s^{\scriptscriptstyle (n)})$ for each $n \in \N$. We are interested in  investigating gradual changes in the evolution of the jump behaviour and we
assume throughout this paper that there is a driving law behind this evolution which is common for all $n \in \N$. Formally, we introduce a transition kernel $g(y,dz)$ from $([0,1], \Bb([0,1]))$ into $(\R, \Bb)$ such that
\begin{align}
\label{BasTrKAss}
\nu^{(n)}_{s}(dz) = g\Big(\frac{s}{n\Delta_n}, dz\Big)
\end{align}
for $s \in [0,n\Delta_n]$. This transition kernel shall be an element of the set $\Gc$ to be defined below. Throughout the paper $\Bb(A)$ denotes the trace $\sigma$-algebra of a set $A \subset \R$  with respect to the Borel $\sigma$-algebra. 
 

\begin{condition}
\label{mcGDef}Let $\mathcal G$ denote the set of all transition kernels $g(\cdot ,dz)$ from $([0,1], \Bb([0,1]))$ into $(\R, \Bb)$ such that 
\begin{enumerate}[(1)]
\item  For each $y \in [0,1]$ the measure $g(y,dz)$ 
         does not charge $\lbrace 0 \rbrace$. 
\item The function $y \mapsto \int (1 \wedge z^2) g(y,dz)$  is bounded on the interval $[0,1]$.
      \label{mcGDef2}
\item If
			\begin{align*}
			\mathcal I(z) \defeq \begin{cases}
                     [z, \infty), \quad &\text{ for } z > 0 \\
										 (- \infty, z],  \quad &\text{ for } z < 0 
										\end{cases} 
			\end{align*}
			denotes one-sided intervals and 
      $$g(y, z) :=g(y, \iz) = \int_{\iz } g(y,dx); \quad (y,z) \in [0,1] \times \R \setminus \lbrace 0 \rbrace,$$
			then for every $\taili \in \Ron$ there exists a finite set $M^{(\taili)} = \{t_1^{(\taili)}, \ldots, t_{n_z}^{(\taili)} \mid n_z \in \N\} \subset [0,1]$, such that the function $y \mapsto g(y,z)$ is continuous on $[0,1] \setminus M^{(\taili)}$.
			\label{mcGDef20}
\item   For each $y \in [0,1]$ the measure $g(y,dz)$ is absolutely continuous with respect to the Lebesgue measure with density $z \mapsto h(y,z)$,
where the measurable function $h \colon ([0,1] \times \R, \Bb([0,1]) \otimes \Bb) \rightarrow (\R, \Bb)$
  is continuously differentiable with respect to  $z \in \Ron$ for fixed $y \in [0,1]$. The function   $h(y,z) $ and its  derivative will be denoted by
$h_y(z)$ and   $h_y^{\prime}(z)$, respectively. Furthermore, we assume for each $\eps >0$ that
			$$\sup \limits_{y \in [0,1]} \sup \limits_{z\in M_\eps} \Big( h_y(z) + |h_y^{\prime}(z)| \Big) < \infty,$$
			\label{mcGDef3}
			where $M_\eps = (- \infty, -\eps] \cup [\eps, \infty)$. 
\end{enumerate}
\end{condition}

These assumptions are all rather mild. For Lebesgue almost every $y \in [0,1]$, the integral $\int (1 \wedge z^2) g(y,dz)$ needs to be finite by properties of the jump compensator, so part (\ref{mcGDef2}) just serves as a condition on uniform boundedness over time. Part (\ref{mcGDef20}) essentially says that for each $z$ only finitely many discontinuous changes of the jump measure $g(y, \cdot)$ are allowed. 
Finally, note that the existence of a density as in (\ref{mcGDef3}) is a standard condition when estimating a measure in a non-parametric framework. 

In order to investigate gradual changes  in  the jump behaviour of the underlying process we
 follow  \cite{VogDet15} and
 consider a measure of time variation for the jump behaviour which is defined by
 \begin{align}
\label{meaotivar}
D(\kseqi, \gseqi, \taili) \defeq \int \limits_0^\kseqi g(y, \taili) dy - \frac \kseqi \gseqi \int \limits_0^\gseqi g(y,z) dy,
\end{align}
where $ (\kseqi, \gseqi,z) \in  C \times \Ron$  and   
\begin{equation} \label{cset}
 C \defeq \lbrace (\kseqi, \gseqi) \in [0,1]^2 \mid \kseqi \leq \gseqi \rbrace.
 \end{equation}
 Here and throughout this paper  we use the convention $\frac 00 \defeq 1$.

The time varying measure \eqref{meaotivar} will be the main theoretical tool for our inference of gradual changes  in the jump behaviour
 of the process \eqref{ItoSem}. Our analysis will be based on the following observation: Due to $\bar\mu^{\scriptscriptstyle (n)}(ds,du) = ds \nu_s^{\scriptscriptstyle (n)}(du)$ the jump behaviour corresponding to the first $\ip{n \gseqi}$ observations for some $\theta \in (0,1)$ does not vary if and only if the kernel $g( \cdot, dz)$ is Lebesgue almost everywhere constant on the interval $[0, \gseqi]$. In this case we have $D(\kseqi, \gseqi, \taili) \equiv 0$ for all $0 \leq \kseqi \leq \gseqi$ and $z \in \Ron$, since $\kseqi^{-1} \int_0^\kseqi g(y,z) dy$ is constant on $[0, \gseqi]$ for each $z \in \Ron$. If on the other hand $D( \kseqi, \gseqi, \taili) =0$ for all $\kseqi \in [0,\theta]$ and $z \in \Ron$, then 
$$\int \limits_0^\kseqi g(y,z) dy = \kseqi \Big(\frac 1{\gseqi} \int \limits_0^\gseqi g(y,z) dy \Big) =: \kseqi A(z)$$
for each $\kseqi \in [0, \gseqi]$ and fixed $\taili \in \Ron$. Therefore by the fundamental theorem of calculus and Assumption \ref{mcGDef}\eqref{mcGDef20} for each fixed $\taili \in \Ron$ we have $g(y,z) = A(z)$ for every $y \in [0,\gseqi]\setminus M^{(z)}$. As a consequence 
\begin{align}
\label{Kernasgleich}
g(y,z)=A(z)
\end{align}
holds for every $\taili \in \mathbb Q\setminus\{0\}$ and each $y \in [0,\gseqi]$ outside the Lebesgue null set $\bigcup_{z \in \mathbb Q\setminus\{0\}} M^{(z)}$. Due to Assumption \ref{mcGDef}\eqref{mcGDef2} and dominated convergence $A(z)$ is left-continuous for positive $z \in \Ron$ and right-continuous for negative $z \in \Ron$. The same holds for $g(y,z)$ for each fixed $y \in [0, \gseqi]$. Consequently \eqref{Kernasgleich} holds for every $\taili \in \Ron$ and each $y \in [0,\gseqi]$ outside the Lebesgue null set $\bigcup_{z \in \mathbb Q\setminus\{0\}} M^{(z)}$. Thus by the uniqueness theorem for measures the kernel $g( \cdot, dz)$ is on $[0,\gseqi]$ Lebesgue almost everywhere equal to the L\'evy measure defined by $A(z)$.

In practice we restrict ourselves to $z$ which are bounded away from zero, as typically $g(y,z) \to \infty$ as $z \to 0$, at least if we deviate from the (simple) case of finite activity jumps. Below we discuss two standard applications of $D(\kseqi, \gseqi, \taili)$ we have in mind. 

\begin{itemize}
\item[(1)]  ({\it test for a gradual change}) If one defines
 \begin{equation} 
\label{tildDcepsDef}
 \tilde{\mathcal D}^{(\eps)}(\gseqi) \defeq \sup \limits_{|\taili| \geq \eps} \sup \limits_{0 \leq \kseqi \leq \gseqi} |D( \kseqi, \gseqi, \taili)| \\
\end{equation}
for some pre-specified constant $\eps >0$, one can characterize the existence of a change point as follows: There exists a gradual change in the behaviour of the jumps larger than $\eps$ of the process \eqref{ItoSem} if and only if
$$\sup \limits_{\gseqi \in [0,1]} \tilde{\mathcal D}^{(\eps)}(\gseqi)>0.
$$ 
Moreover for the analysis of gradual changes it is equivalent to consider
\begin{equation} \label{measjump}
\mathcal D^{(\eps)}(\gseqi) \defeq \sup \limits_{|\taili| \geq \eps} \sup \limits_{0 \leq \kseqi \leq \gseqi' \leq \gseqi} |D( \kseqi, \gseqi', \taili)|, \\
\end{equation}
because the first time points where $\mathcal D^{(\eps)}$ and $\tilde{\mathcal D}^{(\eps)}$ deviate from zero, if existent, coincide and this point is characteristic for a gradual change as we have seen previously. In this paper we consider $\Dc^{(\eps)}$ only, since due to its monotonicity it simplifies several steps in the proofs and our notation. 
In Section \ref{EstConstrSec}   we construct a consistent estimator, say $\Db_n^{(\eps)}$, of  $ \mathcal{D}^{(\eps)}$.
The  test for gradual changes in the behaviour of
the jumps larger than $\eps$ of the process \eqref{ItoSem} rejects the null hypothesis for large values of $\Db_n^{(\eps)}(1)$.
Quantiles for this test will be derived by a multiplier bootstrap (see Section \ref{sec:test} for details).
\item[(2)]  ({\it estimating the  gradual change point})  In Section \ref{EstConstrSec}   we construct an estimator for  the first  point where the jump behaviour changes (gradually). For this purpose we also use the time varying measure \eqref{meaotivar} and define
\begin{align}
\label{changepoint}
\gseqi^{(\eps)}_0  \defeq \inf \left\{ \gseqi \in [0,1] \mid \Dc^{(\eps)}(\gseqi) >0 \right\},
\end{align}
where we set $\inf \emptyset \defeq 1$. We call $\gseqi^{(\eps)}_0$
the  change point of the jumps larger than $\eps$ of the underlying process \eqref{ItoSem}.
 \end{itemize}
\begin{figure}[t!]
\centering
\includegraphics[width=0.48\textwidth]{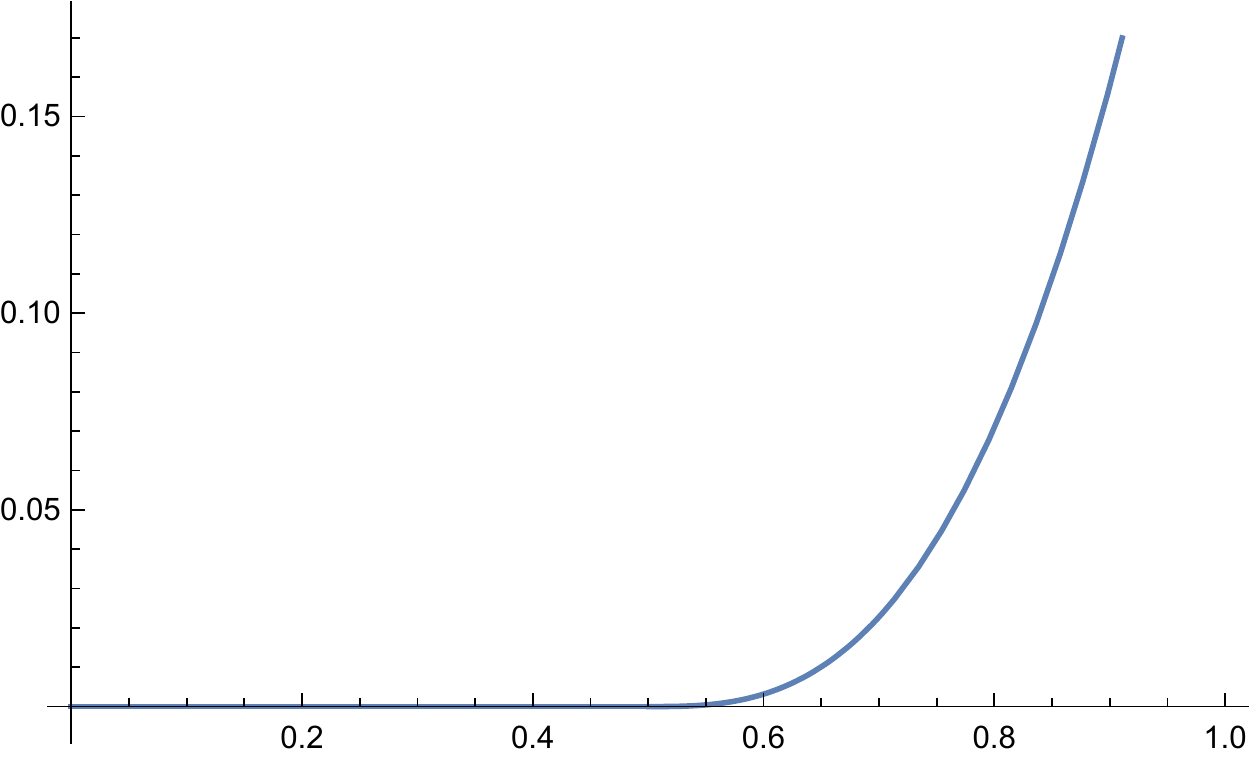}~~~ \includegraphics[width=0.48\textwidth]{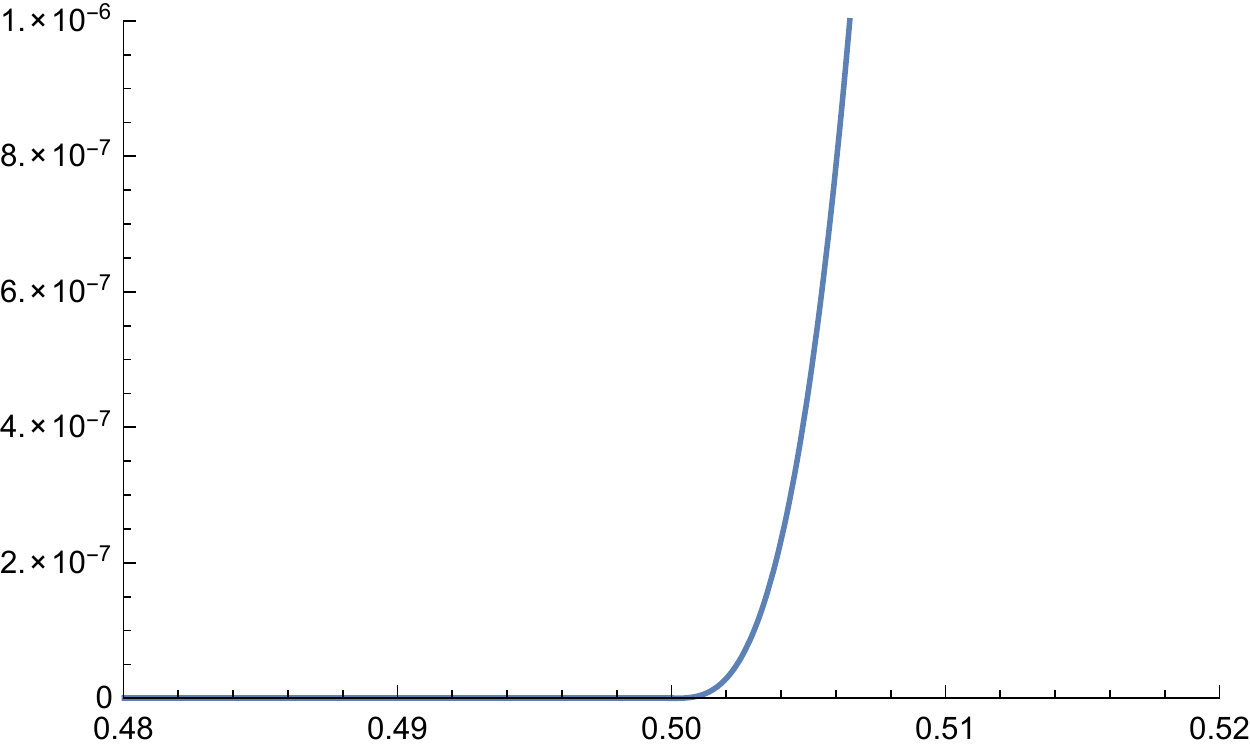}
\caption{\label{example}
\it The function ${\cal D}^{(\varepsilon)}$ for the transition kernel \eqref{gquad}, where $\varepsilon =1$. The ``true'' change point is located at $\theta_0 = 1/2$.
}
\end{figure}
A typical example is displayed in Figure \ref{example}. Here we show the function $\theta \mapsto {\cal D}^{(\varepsilon)} (\theta)$ defined in \eqref{measjump} for $\varepsilon =1$, where the transition kernel is given by
\begin{equation}
 \label{gquad}
g (y,z) = \left \{
\begin{array} {ccc}
10 e^{-|z| } & \mbox {if}  & y \in [0, \frac {1}{2}] \\
  10 \big( 1+3(y- \tfrac {1}{2})^2 \big ) e^{-|z| }& \mbox{if} & y \in [\frac {1}{2}, 1].
\end{array} \ 
\right.
\end{equation}
From the right panel it is clearly visible that the function $\mathcal{D}^{(\varepsilon)}$ is positive for all $\theta \in (\frac {1}{2}, 1]$, which identifies $\theta_0 = 1/2$ as the change point. Additionally we illustrate the previously introduced quantities in two further examples.

\begin{example}  \label{ex1} ({\it abrupt changes})
 {\rm
The classical change point problem, where the jump behaviour of the underlying process is  constant on two intervals,  is contained in our analysis.
 To be precise, assume that  $0 < \gseqi_0< 1$  and that $\nu_1$ and $\nu_2$ are L\'evy measures
  such that the transition kernel   $g$ satisfies
\begin{align}
\label{locojubekedef}
g(y,dz) = \begin{cases} \nu_1(dz), \quad &\text{ for } y \in [0,\gseqi_0] \\
           \nu_2(dz), \quad &\text{ for } y \in (\gseqi_0,1].
					\end{cases}
\end{align}
If each $\nu_j$  is absolutely continuous with respect to the Lebesgue measure and has a density $h_j$ which is continuously differentiable
at any point $z_0 \neq 0$, satisfying
\begin{align*}
\sup \limits_{|z|\ge \eps} \{h_j(z)+|h_j'(z)|\} < \infty
\end{align*}
for every $\eps >0$, then the kernel  $g$ satisfies  Assumption  \ref{mcGDef}. 

For a L\'evy measure $\nu$ on $(\R, \Bb)$ and $z \in \Ron$ let $\nu(z) \defeq \nu(\iz)$. If  $g \in \Gc$ is of  the form \eqref{locojubekedef} and $\eps >0$ is
 chosen sufficiently  small such that there exists a $\bar z \in \mathbb{R}$ with $|\bar z| \ge \eps $  and
\begin{align*}
 \nu_{1}(\bar z) \neq \nu_2(\bar z),
\end{align*}
then we have $D(\kseqi, \gseqi',\taili) =0$ for all $(\kseqi, \gseqi', \taili) \in B_\eps := C \times M_\eps$ with $\gseqi' \leq \gseqi_0$
and consequently $\Dc^{(\eps)}(\gseqi) =0$ for each $\gseqi \leq \gseqi_0$. On the other hand, if $\gseqi_0 < \gseqi ' <  1$ and $\kseqi \leq \gseqi_0$ we have
\begin{align*}
D(\kseqi,\gseqi',\taili) &= \kseqi \nu_1(z) - \frac{\kseqi}{\gseqi'} (\gseqi_0 \nu_1(z) + (\gseqi ' - \gseqi_0) \nu_2(z))
          = \kseqi (\nu_2(z) - \nu_1(z)) \Big( \frac {\gseqi_0}{\gseqi '} -1\Big)
\end{align*}
and obtain
\begin{align*}
\sup \limits_{\kseqi \le \gseqi_0} \sup \limits_{|z| \ge \eps} |D(\kseqi,\gseqi',\taili)| = V_\eps \gseqi_0 \Big( 1-\frac {\gseqi_0}{\gseqi'}\Big),
\end{align*}
where $V_\eps = \sup \limits_{|z| \ge \eps} |\nu_1(z) - \nu_2(z)| >0$. For $\gseqi_0 < \kseqi \leq \gseqi '$  a similar calculation yields
\begin{align*}
D(\kseqi,\gseqi ',\taili)
						&= \gseqi_0(\nu_2(z)-\nu_1(z))\Big(\frac\kseqi{\gseqi'} -1 \Big)
\end{align*}
which gives
\begin{align*}
\sup \limits_{\gseqi_0 < \kseqi \le \gseqi '} \sup \limits_{|z| \ge \eps} |D(\kseqi,\gseqi ',\taili)| = V_\eps \gseqi_0 \Big( 1-\frac{\gseqi_0}{\gseqi'}\Big).
\end{align*}
It follows that
the quantity defined \eqref{changepoint} is given by

 $\chapo = \gseqi_0$, because for $\gseqi > \gseqi_0$ we have
\begin{align}
\label{supcalcEq}
\Dc^{(\eps)}(\gseqi) = \sup \limits_{\gseqi_0 < \gseqi' \leq \gseqi} \max \Big\{ \sup \limits_{\kseqi \leq \gseqi_0} \sup \limits_{|\taili| \geq \eps} |D(\kseqi, \gseqi ', \taili)|, \sup \limits_{ \gseqi_0 < \kseqi \leq \gseqi'} \sup \limits_{|\taili| \geq \eps} |D(\kseqi, \gseqi ', \taili)| \Big\} = V_\eps \gseqi_0 \Big( 1-\frac{\gseqi_0}{\gseqi} \Big).
\end{align}
}
\end{example}

\begin{example}  \label{ex2} ({\it Locally symmetric $\beta$-stable jump behaviour}) {\rm 
A  L\'evy process is symmetric $\beta$-stable for some $0<\beta<2$ if and only if its Brownian part vanishes and its L\'evy measure has a Lebesgue density of the
form $h(z) = A/|z|^{1+\beta}$ with $A \in \R_+$ [see, for instance, Chapter 3 in \cite{Sat99}].
In this sense we say that an It\=o semimartingale with decomposition \eqref{ItoSem} satisfying
 \eqref{BasTrKAss} has locally symmetric
$\beta$-stable jump behaviour, if the corresponding transition kernel  $g$  is given by 
\begin{align}
\label{locbestagDef}
g(y,\iz) = g(y,z) \quad  \text{ and } \quad g(y,\{ 0 \}) = 0, \quad \text{ with } \quad g(y,z) = A(y)/|z|^{\beta(y)}
\end{align}
for $y \in [0,1]$ and $z \in \Ron$. Here  the functions $\beta : [0,1] \rightarrow (0,2)$ and $A:[0,1] \rightarrow (0,\infty)$
are continuous outside a finite set, $A$ is bounded and $\beta$ is bounded away from $2$. In Section \ref{sec:AuxRes} we show that a kernel of the form \eqref{locbestagDef} satisfies Assumption \ref{mcGDef}. 

Now, let $\gseqi_0 \in (0,1)$ with
\begin{align}
\label{erstkonstGl}
A(y) = A_0 \quad \text{ and } \quad \beta(y) = \beta_0 
\end{align} 
for all $y \in [0, \gseqi_0]$ with $A_0 \in (0, \infty)$, $\beta_0 \in (0,2)$. To model the continuous change after $\gseqi_0$ we assume that $\gseqi_0$ is contained in an open interval $U$ and that there exists a real analytic function $\bar A \colon U \rightarrow \R$
and an affine linear function $\bar \beta \colon U \rightarrow \R$ 
such that at least one of the functions $\bar A$, $\bar \beta$ is non-constant as well as 
\begin{align}
\label{analyticass}
A(y) = \bar A(y) \quad \text{ and } \quad \beta (y) = \bar \beta(y)
\end{align} 
for all $y \in [\gseqi_0,1) \cap U$. Then we also show in Section \ref{sec:AuxRes} that the quantity defined in \eqref{changepoint} is given by $\chapo = \gseqi_0$ for every $\eps >0$
(see Section \ref{proofex2}). 
}
\end{example}

We conclude this section with the main assumption for the characteristics of an It\=o semimartingale
which will be used throughout this paper. 

\begin{condition}
\label{Assump1}  For each $n \in \N$ let $X^{(n)}$ denote an It\=o semimartingale  of the form \eqref{ItoSem} with characteristics $(b^{(n)}_s, \sigma^{(n)}_s, \nu^{(n)}_s)$ defined on the probability space $(\Omega,\Fc, \Prob)$ that satisfies
\begin{enumerate}[(a)]
\item \label{RescaTime}There exists a $g \in \Gc$ such that 
       \begin{align*}
       \nu^{(n)}_{s}(dz) = g\Big(\frac s{n\Delta_n}, dz\Big)
        \end{align*}
				holds for all $s \in [0,n\Delta_n]$ and all $n \in \N$.
\item \label{DrVoMomBesch}The drift $b^{(n)}_s$ and the volatility $\sigma^{(n)}_s$ are predictable processes and
satisfy
$$\sup \limits_{n \in \N} \sup \limits_{s \in \R_+} \Big( \Eb|b_s^{(n)}|^{\alpha} \vee \Eb|\sigma_s^{(n)}|^p \Big) < \infty,$$
for some $p>2$, with $\alpha = 3p/(p+4)$.
\item The observation scheme $\{ X^{(n)}_{i \Delta_n} \mid i=0,\ldots,n \}$ satisfies 
\begin{align}
\label{RelSpeed}
\Delta_n \rightarrow 0, \quad n \Delta_n \rightarrow \infty, \quad \text{ and } \quad n \Delta_n^{1+\tau} \rightarrow 0,
\end{align}
for $\tau = (p-2)/(p+1) \in (0,1)$.
\end{enumerate}
\end{condition}

\begin{remark}
Assumption \ref{Assump1}\eqref{RescaTime} corresponds to the typical finding in classical change point analysis that a change point can only be defined relative to the length of the data set which in our case is the covered time horizon $n\Delta_n$. However, this assumption is no restriction in applications since a statistical method is always applied to a data set with a fixed length. Assumption \ref{Assump1}\eqref{DrVoMomBesch} is very mild as it requires only a bound on the moments of drift and volatility (note that $\alpha \in (1,3)$ as $p >2$). Typical for high frequency statistics is condition \eqref{RelSpeed}, which is used to control the bias arising from the presence of the continuous part of the underlying process. Moreover, if $b_s^{\scriptscriptstyle (n)}$ and $\sigma_s^{\scriptscriptstyle (n)}$ are deterministic and bounded functions in $s \in \R_+$ and $n \in \N$,  the process $X^{(n)}$ has independent increments and assumption \eqref{RelSpeed} can be weakened to 
\[ 
\Delta_n \rightarrow 0, \quad n \Delta_n \rightarrow \infty \quad \text{ and } \quad n \Delta_n^{3} \rightarrow 0.
\]
\end{remark}

\section{Weak convergence}
\def\theequation{3.\arabic{equation}}
\setcounter{equation}{0}
\label{sec3}

In order to estimate the measure of time variation  introduced in \eqref{meaotivar} we use
the sequential empirical tail integral process defined by
\begin{align*}
U_n(\gseqi, \taili) = \frac{1}{k_n} \sum \limits_{j=1}^{\ip{n \gseqi}} 1_{\lbrace \Delta_j^n X^{(n)} \in \iz \rbrace},
\end{align*}
where $ \Delta_j^n X^{(n)}= X^{(n)}_{ j \Delta_n } -X^{(n)}_{(j -1) \Delta_n} $,  $\gseqi \in [0,1]$, $\taili \in \Ron$  and  $k_n \defeq n \Delta_n$.
The process $U_n$ counts the number of increments that fall into $\iz$, as these are likely to be caused by a jump with the corresponding size, and will be the basic tool for estimating the measure of time variation defined  in \eqref{meaotivar}.
An estimate is given by 
\begin{align}\label{timevaryest}
\Db_n(\kseqi, \gseqi, \taili) \defeq U_n(\kseqi, \taili) - \frac \kseqi\gseqi U_n(\gseqi, \taili),~~ (\kseqi, \gseqi, \taili) \in C \times \Ron,
\end{align}
where the set $C$ is defined in \eqref{cset}. The statistic $U_n(1, \taili)$ has been considered by
 \cite{FigLop08} for observations of a  L\'evy process $Y$, so without a time-varying jump behaviour. In this case the author
 shows that $U_n(1,\taili)$  is in fact an $L^2$-consistent estimator for the tail integral $\nu(\iz) = U(z)$.
The  following theorem  provides a generalization of this statement. In particular, it provides the weak convergence  of the sequential empirical tail integral
\begin{align}
\label{gndef}
\Gb_n(\gseqi, \taili) \defeq \sqrt{k_n} \Big \{ U_n(\gseqi, \taili) - \int \limits_0^{\gseqi} g(y,z) dy \Big \}.
\end{align}
Throughout this paper we use the notation
$ 
A_\eps = [0,1] \times M_\eps,
$ 

and $R_1 \triangle  R_2$ denotes the symmetric difference of two sets $R_1,R_2$.

\begin{theorem}
\label{WeakConvThm}
If Assumption \ref{Assump1} holds, then the process $\Gb_n$ defined in \eqref{gndef} satisfies $\Gb_n \weak \Gb$ in $\linaeps$ for any $\eps >0$, where
$\Gb$ is a tight mean zero Gaussian process with covariance function
$$H(\gseqi_1, \taili_1; \gseqi_2, \taili_2) = \operatorname{Cov}(\Gb(\gseqi_1, \taili_1), \Gb(\gseqi_2, \taili_2)) = \int \limits_0^{\gseqi_1 \wedge \gseqi_2} g(y, \izei \cap \izzw) dy.$$
The sample paths of $\Gb$ are almost surely uniformly continuous with respect to the semimetric
$$\rho((\gseqi_1, \taili_1); (\gseqi_2, \taili_2)) \defeq \Big \{ \int \limits_0^{\gseqi_1 } g(y, \izei \triangle \izzw) dy + \int \limits_{\gseqi_1 }^{ \gseqi_2} g(y, \mathcal I(z_{2})) dy \Big \}^{\frac 12},$$
defined for $\gseqi_1 \le \gseqi_2$ without loss of generality. Moreover, the space $(A_\eps, \rho)$ is totally bounded.
\end{theorem}
\medskip

Recall the definition of the  measure of time variation for the jump behaviour  defined in  \eqref{meaotivar} and the definition of the set $C$ in \eqref{cset}.  For $B_\eps = C \times M_\eps$
consider the functional  $\Phi \colon \linaeps \rightarrow \linbeps$ defined by
\begin{equation} \label{phidef}
\Phi(f)(\kseqi, \gseqi, \taili) \defeq f(\kseqi, \taili) - \frac\kseqi\gseqi f(\gseqi, \taili).
\end{equation}
As $\| \Phi(f_1) - \Phi(f_2) \|_{B_\eps} \leq 2 \| f_1 - f_2 \|_{A_\eps}$   the mapping $\Phi$ is Lipschitz continuous.
Consequently, $\Hb \defeq \Phi(\Gb)$ is a tight mean zero
Gaussian process in $\linbeps$ with  covariance structure
\begin{align}
\label{HbGrProCov}
\operatorname{Cov}(\Hb(\kseqi_1, \gseqi_1, &\taili_1),\Hb(\kseqi_2, \gseqi_2, \taili_2)) = \nonumber \\
&= \int \limits_0^{\kseqi_1 \wedge \kseqi_2} g(y,\izei \cap \izzw) dy - \frac{\kseqi_1}{\gseqi_1} \int \limits_0^{\kseqi_2 \wedge \gseqi_1}
   g(y,\izei \cap \izzw) dy \nonumber \\
&\hspace{10mm}- \frac{\kseqi_2}{\gseqi_2} \int \limits_0^{\kseqi_1 \wedge \gseqi_2} g(y,\izei \cap \izzw) dy + \frac{\kseqi_1 \kseqi_2}{\gseqi_1 \gseqi_2} \int \limits_0^{\gseqi_1 \wedge \gseqi_2} g(y,\izei \cap \izzw) dy.
\end{align}
From the continuous  mapping theorem  we obtain weak convergence of the process
 \begin{align}
\label{bbhnDef}
\Hb_n(\kseqi, \gseqi, \taili) := \Phi(\Gb_n)(\kseqi, \gseqi, \taili) = \sqrt{k_n}(\Db_n(\kseqi, \gseqi, \taili) - D(\kseqi, \gseqi, \taili)).
\end{align}

\begin{theorem}
\label{SchwKtomotv}
If Assumption  \ref{Assump1} is satisfied, then  the process $\Hb_n$ defined in \eqref{bbhnDef} satisfies $\Hb_n \weak \Hb$ in $\linbeps$ for any $\eps >0$, where $\Hb$ is a tight mean zero Gaussian process with covariance function \eqref{HbGrProCov}.
\end{theorem}



For the statistical change-point inference proposed  in the following section we require the quantiles of
functionals of  the limiting distribution in Theorem \ref{SchwKtomotv}.   This  distribution depends in a complicated way on the unknown underlying kernel $g \in \Gc$
 and, as a consequence, corresponding quantiles are difficult to estimate.   
 
A typical approach to problems of this type are resampling methods.
One option  is to use  suitable estimates for drift, volatility  and  the unknown kernel $g$ to draw  independent samples of an It\=o semimartingale.
However,   such a method is computationally expensive  since one has to  generate independent It\=o
 semimartingales for each stage within the bootstrap algorithm. Therefore
 we  propose an alternative bootstrap method based on multipliers.
 For this resampling method one only needs to  generate $n$ i.i.d.\ random variables with mean zero and variance one. 

To be  precise let $X_1, \ldots, X_n$  and $\xi_1, \ldots, \xi_n$  denote random variables  defined on probability spaces  $(\Omega_X, \mathcal A_X, \mathbb P_X)$
and
 $(\Omega_{\xi}, \mathcal A_{\xi}, \mathbb P_{\xi})$, respectively,  and  consider  a random element $\hat Y_n = \hat Y_n(X_1, \ldots, X_n,$ $ \xi_1, \ldots, \xi_n)$
 on the product space
$(\Omega, \mathcal A,\mathbb P) \defeq (\Omega_X, \mathcal A_X, \mathbb P_X) \otimes (\Omega_{\xi}, \mathcal A_{\xi}, \mathbb P_{\xi})$
which maps into a  metric space, say  $\mathbb D$. Moreover, let $Y$ be a tight, Borel measurable  $\mathbb D$-valued random variable.
Following  \cite{Kos08}  we call   $\hat Y_n$ {\it weakly convergent to $Y$ conditional on  $X_1, X_2, \ldots$ in probability} if  the following
two conditions are satisfied
\begin{enumerate}
\item[(a)] $\sup \limits_{f \in \text{BL}_1(\mathbb D)} |\mathbb E_{\xi} f(\hat Y_n) - \mathbb E f(Y) | \pto 0,$
\item[(b)] $\mathbb E_{\xi} f( \hat Y_n)^{\ast} - \mathbb E_{\xi} f( \hat Y_n)_{\ast} \pto 0$ for all
      $f \in \text{BL}_1(\mathbb D).$
\end{enumerate}
Here, $\mathbb E_{\xi}$ denotes the conditional expectation  with respect to  $\xi_1,\ldots,\xi_n$ given  $X_1, \ldots, X_n$, 
whereas $\text{BL}_1(\mathbb D)$ is the space of all real-valued Lipschitz continuous functions $f$ on $\mathbb D$ with sup-norm $\| f \|_{\infty} \leq 1$ and Lipschitz constant $1$. Moreover, $f( \hat Y_n)^{\ast}$ and $f( \hat Y_n)_{\ast}$ denote a minimal measurable majorant and a maximal measurable minorant with respect to $\xi_1,\ldots,\xi_n$,  $X_1, \ldots, X_n$,   respectively.   Throughout this paper we denote this type of convergence by  $\hat Y_n \weakP Y$.


In the following we will work with  a multiplier bootstrap version of the process $\Gb_n $, that is
\begin{eqnarray}\label{hatGbndefeq}
\hat \Gb_n &=& \hat \Gb_n(\theta,z) =
\hat \Gb_n (X^{(n)}_{\Delta_n}, \ldots, X^{(n)}_{n \Delta_n}, \xi_1, \ldots, \xi_n;\theta,z) \nonumber \\ 
& \defeq &
\frac{1}{n \sqrt{k_n}} \sum \limits_{j=1}^{\lfloor n\theta \rfloor} \sum \limits_{i=1}^n \xi_j \lbrace \mathtt 1_{\lbrace \Delta_j^n X^{(n)} \in \iz \rbrace} - \mathtt 1_{\lbrace \Delta_i^n X^{(n)} \in \iz \rbrace} \rbrace  \nonumber \\ 
&=&
\frac{1}{\sqrt{k_n}} \sum \limits_{j=1}^{\lfloor n\theta \rfloor} \xi_j \lbrace \mathtt 1_{\lbrace \Delta_j^n X^{(n)} \in \iz \rbrace} -\eta_n(z) \rbrace,
\end{eqnarray}
 where  $\xi_1, \ldots , \xi_n$  are independent and identically distributed random variables with mean $0$ and variance $1$ and
$\eta_n(z)= n^{-1} \sum_{i=1}^n \mathtt 1_{\lbrace \Delta_i^n X^{(n)} \in \iz \rbrace}
$.
The following  theorem establishes conditional weak convergence of this bootstrap approximation for the sequential empirical tail integral process $\Gb_n$.

\begin{theorem}
\label{GnVorberKonvThm}
If Assumption \ref{Assump1} is satisfied  and   $(\xi_j)_{j \in \mathbb N}$ is a sequence of independent and identically distributed random variables with mean $0$ and variance $1$, defined on a distinct probability space as described above,
then the process $\hat \Gb_n$ defined in \eqref{hatGbndefeq} satisfies $\hat \Gb_n \weakP \mathbb \Gb$
in $\linaeps$ for any $\eps > 0$, where $\mathbb G$ denotes the limiting process of Theorem \ref{WeakConvThm}.
\end{theorem}
 \medskip

Theorem~\ref{GnVorberKonvThm} suggests to consider the following   counterparts of the process $\Hb_n$ defined in \eqref{bbhnDef} 
\begin{align}
\label{hatHbndefeq}
\hat \Hb_n (\kseqi, \gseqi, \taili)
&\defeq
\hat \Hb_n (X^{(n)}_{\Delta_n}, \ldots, X^{(n)}_{n \Delta_n}; \xi_1, \ldots, \xi_n; \kseqi, \gseqi, \taili )
\defeq
\hat \Gb_n(\kseqi,z) - \frac{\kseqi}{\gseqi} \hat\Gb_n(\gseqi,z) \nonumber \\
&=
\frac{1}{\sqrt {n \Delta_n}} \bigg[ \sum \limits_{j=1}^{\lfloor n\kseqi \rfloor} \xi_j \{ \mathtt 1_{\lbrace \Delta_j^n X^{(n)} \in \iz \rbrace} - \eta_n(z) \} - \nonumber \\
&\hspace{5cm}
- \frac{\kseqi}{\gseqi} \sum \limits_{j = 1}^{\ip{n \gseqi}} \xi_j \{ \mathtt 1_{\lbrace \Delta_j^n X^{(n)} \in \iz \rbrace} - \eta_n(z) \} \bigg ].
\end{align}
The following result establishes consistency of $\hat \Hb_n$. Its proof
 is a consequence of Proposition 10.7 in \cite{Kos08}, because we have $\hat \Hb_n = \Phi( \hat \Gb_n)$ and $\Hb = \Phi(\Gb)$
with the   Lipschitz continuous map $\Phi   $ defined in  \eqref{phidef}.

 \medskip

\begin{theorem}
\label{BootstrTeststThm}
If Assumption  \ref{Assump1} holds and   $(\xi_j)_{j \in \mathbb N}$ is a sequence of independent and identically distributed random variables with mean $0$ and variance $1$ defined on a distinct probability space, then the process $\hat \Hb_n$ defined in \eqref{hatHbndefeq} satisfies  $ \hat \Hb_n \weakP \mathbb H$
in $\linbeps$ for any $\eps > 0$, where    the  
process $\mathbb H$ is defined in Theorem \ref{SchwKtomotv}.
\end{theorem}


\section{Statistical inference for gradual changes}
\def\theequation{4.\arabic{equation}}
\setcounter{equation}{0}
\label{sec4}

As we have seen in Section \ref{sec2} the quantity $\Dc^{(\eps)}(\gseqi)$ defined in \eqref{measjump} is an indicator for a change
in the behaviour of the jumps larger than $\eps$ on $[0,\ip{n\gseqi} \Delta_n]$. Therefore we use the estimate $\Db_n(\kseqi, \gseqi, \taili)$  of the measure of time variation
$
D(\kseqi, \gseqi, \taili) $
    defined in \eqref{timevaryest}
to construct both a test for the existence and an estimator for the location of  a gradual change. We begin with the  problem of estimating the
point of such a gradual change in the jump behaviour. The discussion of corresponding tests
will be referred to Section \ref{sec:test}.

\subsection{Localizing change points}
\label{EstConstrSec}

 Recall the definition $$
 \Dc^{(\eps)}(\gseqi) = \sup \limits_{|\taili| \geq \eps} \sup \limits_{0 \leq \kseqi \leq \gseqi' \leq \gseqi} |D( \kseqi, \gseqi', \taili)| 
 $$
and the definition
of the change point $\gseqi^{(\eps)}_0$  in \eqref{changepoint}. By Theorem \ref{SchwKtomotv}  the process $ \Db_n(\kseqi,\gseqi,\taili) $  from \eqref{timevaryest}
is a consistent estimator of $D(\kseqi, \gseqi, \taili) $. Therefore we set
\begin{align*}
\Db_n^{(\eps)}(\gseqi) \defeq \sup \limits_{|\taili| \geq \eps} \sup \limits_{0 \leq \kseqi \leq \gseqi' \leq \gseqi} |\Db_n(\kseqi,\gseqi',\taili)|,
\end{align*}
and an application of the continuous mapping theorem and  Theorem \ref{SchwKtomotv}  yields  the following result.

\begin{corollary} \label{51}
If Assumption  \ref{Assump1} is satisfied, then    $k_n^{1/2} \Db_n^{(\eps)} \weak \Hb^{(\eps)}$ in $\ell^{\infty}\big([0,\gseqi^{(\eps)}_0]\big)$, where $\Hb^{(\eps)}$ is
the tight process in $\linne$ defined by
\begin{align}
\label{limhbepsDef}
\Hb^{(\eps)} (\gseqi) \defeq \sup \limits_{|z| \geq \eps} \sup \limits_{0 \leq \kseqi \leq \gseqi' \leq \gseqi} |\Hb(\kseqi, \gseqi', \taili)|,
\end{align}
with the centered Gaussian process $\Hb$ defined in Theorem \ref{SchwKtomotv}.
\end{corollary}


Intuitively, the estimation of $\gseqi^{(\eps)}_0$
becomes more difficult the flatter the curve $\theta \mapsto \Dc^{(\eps)}(\gseqi) $ is at  $\gseqi^{(\eps)}_0$.
Following \cite{VogDet15}, we describe the curvature  of $\theta \mapsto \Dc^{(\eps)}(\gseqi) $  by a local polynomial behaviour of the function $ \Dc^{(\eps)}(\gseqi) $
for values $\theta >\gseqi^{(\eps)}_0$. More precisely, we assume throughout this section  
that $\gseqi^{(\eps)}_0 <1$  and that there exist constants $\lambda, \eta, \smooi, c^{(\eps)}  >0$ such that $\Dc^{(\eps)}$ admits an  expansion of  the form
\begin{equation} \label{addass}
\Dc^{(\eps)}(\gseqi) = c^{(\eps)} \big( \gseqi - \gseqi_0^{(\eps)}\big)^{\smooi} + \aleph(\gseqi)
\end{equation}
for all  $\gseqi \in [\gseqi_0^{(\eps)}, \gseqi_0^{(\eps)} + \lambda]$, where the remainder term satisfies
 $|\aleph(\gseqi)| \leq K\big(\gseqi - \gseqi_0^{(\eps)} \big)^{\smooi + \eta}$ for some $K>0$.
 The construction  of an estimator for $\gseqi^{(\eps)}_0$ utilizes the fact 
  that,  by  Theorem \ref{SchwKtomotv},    $k_n^{1/2} \Db_n^{(\eps)}(\gseqi) \rightarrow \infty$ in probability  for any $\theta \in  (\gseqi^{(\eps)}_0,  1]$.
We now consider the statistic
$$ r^{(\eps)}_n(\gseqi) \defeq 1_{\lbrace k_n^{1/2} \Db_n^{(\eps)}(\gseqi) \leq \thrle_n \rbrace},$$
for a deterministic sequence  $\thrle_n \rightarrow \infty$. From the previous discussion we expect
$$r^{(\eps)}_n(\gseqi) \rightarrow \begin{cases}
                               1, \quad &\text{ if } \gseqi \leq \gseqi^{(\eps)}_0 \\
															 0, \quad &\text{ if } \gseqi > \gseqi^{(\eps)}_0
															 \end{cases}$$
in probability if the threshold level $\thrle_n$ is chosen appropriately.  Consequently, we define the estimator for the change point by
\begin{align}
\label{Estimator}
\hat \theta^{(\varepsilon)}_n = \hat \theta^{(\varepsilon)}_n(\thrle_n) \defeq \int \limits_0^1 r^{(\eps)}_n(\gseqi) d \gseqi.
\end{align}
Note that the estimate $\hat \theta^{(\varepsilon)}_n$ depends on the threshold $\thrle_n$ and we make this dependence visible in our notation whenever it is necessary. Our first result establishes consistency of the estimator $\hat \theta^{(\varepsilon)}_n$ under rather mild assumptions on the sequence $(\thrle_n)_{n \in \N}$.

\begin{theorem}
\label{SchaeKonvThm}
If   Assumption \ref{Assump1} is satisfied, $\gseqi^{(\eps)}_0 <1$, and  \eqref{addass} holds
for some $\smooi >0$, then
$$\consest - \gseqi^{(\eps)}_0 = O_\Prob \Big (  \big ( { \frac{\thrle_n}{\sqrt{k_n}}}\big)^{1/\smooi} \Big),
$$
for any  sequence $\thrle_n \rightarrow \infty$ with $\thrle_n / \sqrt{k_n} \rightarrow 0$.
\end{theorem}

Theorem \ref{SchaeKonvThm} makes the heuristic argument of the previous paragraph more precise. A lower
degree of smoothness in $\gseqi^{(\eps)}_0$  yields a  better  rate of convergence of the  estimator.
Moreover, the slower the threshold level $\thrle_n$ converges  to infinity the better the rate of convergence.
We will explain below how to choose this sequence to control the probability of over- and underestimation
by using bootstrap methods.
Before that we investigate the mean squared error 
\begin{align*}
\operatorname{MSE}^{(\eps)}(\thrle_n)= \Eb \Big[ \big( \hat \theta^{(\varepsilon)}_n(\thrle_n) - \gseqi^{(\eps)}_0 \big)^2 \Big]
\end{align*}
of the estimator $\hat \theta^{(\varepsilon)}_n$. Recall the definition of $\Hb_n$ in \eqref{bbhnDef} and define
\begin{align}
\label{Fnepsdef}
\Hb_n^{(\eps)}(\gseqi) \defeq \sup \limits_{|\taili| \ge \eps} \sup \limits_{0 \leq \kseqi \leq \gseqi ' \leq \gseqi} |\Hb_n(\kseqi, \gseqi ', \taili)|, \quad \gseqi \in [0,1], 
\end{align}
which measures the absolute distance between the estimator $\Db_n^{(\eps)}(\gseqi)$ and the true value $\Dc^{(\eps)}(\gseqi)$. For a sequence $\alpha_n \rightarrow \infty$ with $\alpha_n = o(\thrle_n)$ we decompose the MSE into 
\begin{eqnarray*}
\text{MSE}_1^{(\eps)}(\thrle_n,\alpha_n) & \defeq & \Eb \Big[ \big( \consest - \gseqi^{(\eps)}_0 \big)^2 1_{\left \{ \Hb_n^{(\eps)}(1) \leq \alpha_n \right\}} \Big], \\
\text{MSE}_2^{(\eps)}(\thrle_n,\alpha_n)&  \defeq  & \Eb \Big[ \big( \consest - \gseqi^{(\eps)}_0 \big)^2 1_{\left \{ \Hb_n^{(\eps)}(1) > \alpha_n \right\}} \Big] \leq \Prob \big( \Hb_n^{(\eps)}(1) > \alpha_n \big),
\end{eqnarray*}
which can be considered as the MSE due to small and large estimation error. With these notations the following theorem gives upper and lower bounds for the mean squared error.

\begin{theorem}
\label{MSEzerlThm}
Suppose that $\gseqi^{(\eps)}_0 <1$, Assumption \ref{Assump1} and \eqref{addass} are satisfied. Then for any
sequence $\alpha_n \rightarrow \infty$ with $\alpha_n = o(\thrle_n)$ we have
\begin{align}
\label{MSE1queeqn}
K_1 \big ( { \frac{\thrle_n}{\sqrt{k_n}} }\big)^{2/\smooi} \leq &\operatorname{MSE}_1^{(\eps)}(\thrle_n,\alpha_n) \leq K_2 \big ( { \frac{\thrle_n}{\sqrt{k_n}} }\big)^{2/\smooi}  \\
\nonumber
&\operatorname{MSE}_2^{(\eps)}(\thrle_n,\alpha_n) \leq  \Prob \big( \Hb_n^{(\eps)}(1) > \alpha_n \big),
\end{align}
for $n \in \N$ sufficiently large, where the constants $K_1$ and $K_2$ can be chosen as
\begin{align}
\label{KeinKzweDef}
K_1 = \bigg( \frac{1-\varphi}{c^{(\eps)}} \bigg)^{2/\smooi} \quad \text{ and } \quad K_2 = \bigg( \frac{1+ \varphi}{c^{(\eps)}} \bigg)^{2/\smooi}
\end{align}
for an arbitrary $\varphi \in (0,1)$.
\end{theorem}

In the remaining part of this section we discuss the choice of the regularizing sequence $\thrle_n$ for the estimator $\hat \theta^{(\varepsilon)}_n$. Our main goal here is to control the probability of over- and underestimation of the change point $\gseqi_0^{(\eps)} \in (0,1)$.

For this purpose let $\preest$ be a preliminary consistent estimator of  $ \theta^{(\varepsilon)}_0$. For example, if \eqref{addass} holds for some $\smooi >0$, one can take $\preest = \consest(\thrle_n)$ for a sequence $\thrle_n \rightarrow \infty$ satisfying the assumptions of Theorem \ref{SchaeKonvThm}. 
In the sequel, let $B\in\N$ be some large number and let $(\xi^{\scriptscriptstyle  (b)})_{b=1, \dots ,B}$ denote independent vectors of i.i.d.\ random variables, $\xi^{\scriptscriptstyle (b)} \defeq (\xi_j^{\scriptscriptstyle (b)})_{j=1, \dots, n}$,  with mean zero and variance one, which are defined on a probability space distinct to the one generating the data $\{X_{i \Delta_n}^{(n)} \mid i = 0, \ldots, n\}$.  We denote by $\hat \Gb_{\scriptscriptstyle  n, \xi^{(b)}}$ or $\hat \Hb^{(\eps)}_{\scriptscriptstyle  n, \xi^{(b)}}$ the particular bootstrap statistics calculated with respect to the data and the bootstrap multipliers $\xi^{\scriptscriptstyle (b)}_1, \ldots, \xi^{\scriptscriptstyle  (b)}_n$ from the $b$-th iteration, where
\begin{align}
\label{hatFnepsdef}
\hat \Hb_n^{(\eps)}(\gseqi) \defeq \sup \limits_{|\taili| \ge \eps} \sup \limits_{0 \leq \kseqi \leq \gseqi ' \leq \gseqi} |\hat \Hb_n(\kseqi, \gseqi ', \taili)|
\end{align}
for $\gseqi \in [0,1]$. With these notations and for $\eps >0$, $B,n \in \N$ and $r \in (0,1]$ we define the following empirical distribution function
\begin{align*}
\emprbodif (x) &= \frac 1B \sum \limits_{i=1}^B \mathtt 1_{\{  ( \bootsuhb_{\scriptscriptstyle n, \xi^{(i)}}(\preest) )^r \leq x  \}}, \\
\end{align*}
and denote by 
\begin{align*}
\emprboindif(y) := \inf \Big\{x \in \R \mid \emprbodif(x) \geq y \Big\}
\end{align*} 
its pseudoinverse.
Given a confidence level $0<\alpha<1$ we consider the threshold
\begin{align}
\label{thrledefeq}
\hathrle :=
\emprboindif(1-\alpha). 
\end{align}
This choice is optimal in the sense of the following two theorems. 

\begin{theorem}
\label{OptChoofThThm}
Let $\eps >0$, $0<\alpha<1$ and assume that  Assumption \ref{Assump1} is satisfied for some $g \in \Gc$ with $0<\chapo <1$. Suppose further that there exists some $\bar z \in M_\eps$ with
\begin{align}
\label{nichtdegenAss}
\int \limits_0^{\chapo} g(y, \bar z) dy >0.
\end{align}
Then the probability for underestimation of the change point $\chapo$ can be controlled by
\begin{align}
\label{underestabEq}
\limsup \limits_{B \rightarrow \infty} \limsup \limits_{n \rightarrow \infty} \Prob\Big(\consest(\hathrleei) < \chapo \Big) \leq \alpha.
\end{align}
\end{theorem}

\begin{theorem}
\label{KorThrChoiThm}
Let $\eps >0$, $0<r<1$. Assume that  Assumption \ref{Assump1} is satisfied for some  $g \in \Gc$ with $0<\chapo <1$ and  that \eqref{addass} holds for some $\smooi, c^{(\eps)} >0$. Furthermore suppose that there exist a constant $\rho > 0$ with $n \Delta_n^{1+ \rho} \rightarrow \infty$ and a  $\bar \taili \in M_\eps$ satisfying \eqref{nichtdegenAss}.
Additionally let the bootstrap multipliers be either bounded in absolute value or standard normal distributed. Then for each $K > \big(  1/c^{(\eps)} \big)^{1/\smooi}$ and all sequences $(\alpha_n)_{n \in \N} \subset (0,1)$ with $\alpha_n \rightarrow 0$ and $(B_n)_{n \in \N} \subset \N$ with $B_n \rightarrow \infty$ such that
\begin{enumerate}
\item
$\alpha_n^2 B_n \rightarrow \infty,$
\item $(n \Delta_n)^{\frac{1-r}{2r}} \alpha_n \rightarrow \infty,$
\end{enumerate}
we have
\begin{align}
\label{overestabEq}
\lim \limits_{n \rightarrow \infty} \Prob\Big(\consest(\hacothrleor) > \chapo + K \beta_n \Big) =0,
\end{align}
where $\beta_{n} = (\hacothrleor/\sqrt{k_n})^{1/\smooi} \probto 0$, while $\hacothrleor \probto \infty$.
\end{theorem}

Obviously, Theorem \ref{KorThrChoiThm} only gives a meaningful result since $\beta_{n} \probto 0$ can be guaranteed. Its proof shows that
a sufficient condition for this property is given by
\begin{align}
\label{imporEq}
 \Prob \Big( \bootsuhb_{n}(\preest)  > (\sqrt{k_n} x)^{1/r} \Big) \leq \Prob \Big( \bootsuhb_{n}(1)  > (\sqrt{k_n} x)^{1/r} \Big) = o(\alpha_n),
\end{align}
for arbitrary $x >0$. Moreover,  (\ref{imporEq}) follows  from $(n \Delta_n)^{\frac{1-r}{2r}} \alpha_n \rightarrow \infty$ without any further conditions. This explains why the threshold $0 < r < 1$ needs to be introduced, and it seems that the statement of \eqref{imporEq} can only be guaranteed under very restrictive assumptions in the case $r=1$.

Finally we illustrate that the polynomial behaviour introduced in \eqref{addass} is satisfied in the situations of Example \ref{ex1} and Example \ref{ex2}.

\begin{example}
\label{ex3} ~
\begin{enumerate}[(1)]
\item Recall the situation of an abrupt change in the jump characteristic
 considered in Example \ref{ex1}. In this case  it follows from \eqref{supcalcEq} that
 $$
 \Dc^{(\eps)}(\gseqi) = V_\eps \gseqi_0 \Big( 1-\frac{\gseqi_0}\gseqi\Big)
 = V_\eps(\gseqi - \gseqi_0) - \frac{V_\eps}\gseqi(\gseqi - \gseqi_0)^2
 >0,
 $$
 whenever $\gseqi_0 < \gseqi \leq 1 $. Therefore   assumption \eqref{addass} is satisfied with
$\aleph(\gseqi) = -\frac{V_\eps}\gseqi (\gseqi - \gseqi_0)^2$. Moreover, the transition kernel given by
 \eqref{locojubekedef} satisfies
 assumption \eqref{nichtdegenAss} if $\nu_1 \neq 0$ and $\eps >0$ is chosen small enough.
\item \label{ex3nr2} In the situation discussed in Example \ref{ex2} let $\bar g(y,z) = \bar A(y)/|z|^{\bar \beta(y)}$ for $y \in U$ and $z \in M_\eps$. Then we have for any $\eps >0$
            $$k_{0,\eps} := \min \Big\{ k \in \N \mid \exists z \in M_\eps : g_k(z) \neq 0 \Big\} < \infty,$$
						where for $k \in \N_0$ and $\taili \in \Ron$
						$$g_k(z) := \Big(\frac{\partial^k \bar g}{\partial y^k} \Big) \Big|_{(\gseqi_0, z)}$$
						denotes the $k$-th partial derivative of $\bar g$ with respect to $y$ at $(\gseqi_0, z)$, which is a bounded function on any $M_\eps$. Furthermore, for every $\eps >0$ there is a $\lambda >0$ such that 
				\begin{align}
				\label{smoexpanGl}
				\Dc^{(\eps)}(\gseqi) = \Big(\frac{1}{(k_{0,\eps} +1)!} \sup \limits_{|z| \geq \eps}|g_{k_{0,\eps}}(z)| \Big) (\gseqi - \gseqi_0)^{k_{0,\eps} +1} + \aleph(\gseqi)
				\end{align}
				on $[\gseqi_0, \gseqi_0 + \lambda]$ with $|\aleph(\gseqi)| \leq K\big(\gseqi - \gseqi_0 \big)^{k_{0,\eps} +2}$ for some $K>0$, so \eqref{addass} is satisfied. A proof of this result can be found in Section \ref{sec:AuxRes} as well. Again, \eqref{nichtdegenAss} holds also.
\end{enumerate}

\end{example}

\subsection{Testing for a gradual change} \label{sec:test}

In this section we want to derive test procedures for the existence of a gradual change in the data. In order to formulate suitable hypotheses for a gradual change point   recall
the definition of the measure of time variation for the jump behaviour in
\eqref{meaotivar} and define for  $\eps >0$, $z_0 \in \Ron$ and  $\gseqi \in [0,1]$ the quantities
\begin{align*}
\mathcal D^{(\eps)}(\gseqi) &\defeq \sup \limits_{|\taili| \geq \eps} \sup \limits_{0 \leq \kseqi \leq \gseqi' \leq \gseqi} |D( \kseqi, \gseqi', \taili)| \\
\mathcal D_{(z_0)}(\gseqi) &\defeq \sup \limits_{0 \leq \kseqi \leq  \gseqi' \leq \gseqi} |D(\kseqi, \gseqi', \taili_0)|.
\end{align*}
We also assume that  Assumption~\ref{Assump1} is satisfied and we are interested in the following
hypotheses
\begin{eqnarray} \label{h0global}
{\bf H}_0 (\eps):~  \mathcal D^{(\eps)}(1) =0 \quad ~\mbox{ versus } ~ \quad
{\bf H}_1 (\eps):~ \mathcal D^{(\eps)}(1) >0,
\end{eqnarray}
which refer to the global  behaviour of the tail integral. If one is interested
in the  gradual change in the tail integral for a fixed $z_0 \in \Ron$ one could consider
the hypotheses
\begin{eqnarray} \label{h0lokal}
{\bf H}_0^{(z_0)}:~
 \mathcal D_{(z_0)}(1) = 0 \quad ~\mbox{ versus } ~ \quad {\bf H}_1^{(z_0)}:~  \mathcal D_{(z_0)}(1) > 0.
\end{eqnarray}

\begin{remark}
Note that  the function $D$ in \eqref{meaotivar} is uniformly continuous in $(\kseqi,\gseqi) \in C$, uniformly in $\taili \in M_\eps$, that is for any $\eta >0$ there exists a $\delta >0$ such that
$$|D (\kseqi_1,\gseqi_1,\taili) - D (\kseqi_2, \gseqi_2, \taili)| < \eta$$
holds for each $\taili \in M_\eps$ and all pairs $(\kseqi_1,\gseqi_1), (\kseqi_2, \gseqi_2) \in C = \{(\kseqi, \gseqi) \in [0,1]^2 \mid \kseqi \leq \gseqi \}$ with maximum distance $\|(\kseqi_1,\gseqi_1)- (\kseqi_2, \gseqi_2)\|_\infty < \delta$.
Therefore the function $D^{(\eps)}(\kseqi, \gseqi) = \sup_{\taili \in M_\eps} |D(\kseqi, \gseqi, \taili)|$ is uniformly continuous on $ C$ as well, and as a consequence $\Dc^{(\eps)}$ is continuous on $[0,1]$. Thus the alternative ${\bf H}_1(\eps)$ holds if and only if the point $\theta_0^{(\varepsilon)} $ defined in  \eqref{changepoint} satisfies $\theta_0^{(\varepsilon)} < 1$.
\end{remark}

The null hypothesis in \eqref{h0global} and \eqref{h0lokal} will be rejected for large values  of the corresponding estimators
$$
\Db_n^{(\eps)}(1) \quad \mbox{and} \quad \sup_{(\kseqi, \gseqi) \in C} |\Db_n(\kseqi, \gseqi,z_0)|
$$
for $\mathcal D^{(\eps)}(1) $
and  $\mathcal D_{(z_0)}(1) $, respectively. The critical values
are obtained by the multiplier bootstrap introduced in  Section \ref{sec3}.
For this purpose we denote by  $\xi^{(b)}_1, \ldots, \xi^{(b)}_n$, $b=1,\ldots , B$, i.i.d.\ random variables
with mean zero and variance one. As before, we assume that these random variables are defined on a probability space distinct to the one generating the data $\{ X^{(n)}_{i \Delta_n} \mid i=0,\ldots,n\}$.
We denote by $\hat \Gb_{\scriptscriptstyle  n, \xi^{(b)}}$ and $\hat \Hb_{\scriptscriptstyle  n, \xi^{(b)}}$ the  statistics in \eqref{hatGbndefeq} and \eqref{hatHbndefeq} calculated from $\{ X^{(n)}_{i \Delta_n} \mid i=0,\ldots,n\}$  and the $b$-th bootstrap multipliers $\xi^{\scriptscriptstyle (b)}_1, \ldots, \xi^{\scriptscriptstyle  (b)}_n$.
For given $\eps >0$, $z_0 \in \Ron$ and a given level $\alpha \in (0,1)$, we
propose to reject  ${\bf H}_0(\eps)$ in favor of ${\bf H}_1(\eps)$, if
\begin{equation} \label{testglobal}
 k_n^{1/2} \Db_n^{(\eps)}(1) \geq \hat q^{(B)}_{1 - \alpha} \Big( {\Hb}^{(\eps)}_n(1) \Big),
\end{equation}
where $\hat q^{(B)}_{1 - \alpha} \Big({\Hb}^{(\eps)}_n(1)  \Big)$ denotes the $(1- \alpha)$-quantile of the sample $\hat{\Hb}_{\scriptscriptstyle n,\xi^{(1)}}^{(\eps)}(1), \ldots, \hat{\Hb}_{\scriptscriptstyle n, \xi^{(B)}}^{(\eps)}(1)$ with $\hat{\Hb}^{(\varepsilon)}_{n,\xi^{(b)}}$ defined in \eqref{hatFnepsdef}.
Note that under the null hypothesis it follows from the definition of the process $\Hb_n$ in \eqref{bbhnDef} that $k_n^{1/2}\Db^{(\varepsilon)}_n(1) = \Hb^{(\varepsilon)}_{n}(1)$, 
which by Theorem \ref{SchwKtomotv} and the continuous mapping theorem converges weakly to ${\Hb}^{(\varepsilon)}(1)$, defined in \eqref{limhbepsDef}. 
The bootstrap procedure mimics this behaviour.

 Similarly, the hypothesis
${\bf H}_0^{(z_0)}$ is rejected in favor of ${\bf H}_1^{(z_0)}$ if
\begin{equation} \label{testlokal}
W_n^{(z_0)} \defeq k_n^{1/2} \sup \limits_{(\kseqi, \gseqi) \in C} |\Db_n(\kseqi, \gseqi,z_0)| \geq \hat q^{(B)}_{1 - \alpha}(W^{(z_0)}_n),
\end{equation}
where $\hat q^{(B)}_{1 - \alpha}(W^{(z_0)}_n)$ denotes the $(1- \alpha)$-quantile of the sample $\hat W_{\scriptscriptstyle n, \xi^{(1)}}^{\scriptscriptstyle (z_0)}, \ldots, \hat W_{\scriptscriptstyle n, \xi^{(B)}}^{\scriptscriptstyle (z_0)}$, and
$$
\hat W_{\scriptscriptstyle n, \xi^{(b)}}^{\scriptscriptstyle (z_0)} \defeq \sup_{(\kseqi, \gseqi) \in C} | \hat \Hb_{n, \xi^{(b)}}(\kseqi, \gseqi,z_0) |.
$$ 

\begin{remark} 
Since $\eps > 0$ has to be chosen for an application of the test \eqref{testglobal}, one can only detect changes in the jumps larger than $\eps$.
From a practical point of view this is not a severe restriction as in most applications only the larger jumps are of interest. If one is interested in the entire jump measure, however, its estimation is rather difficult, at least in the presence of a diffusion component, as $\Delta_n^{1/2}$ provides a natural bound to disentangle jumps from volatility. We refer to \cite{Nic16} and \cite{HofVet15} for details in the case of L\'evy processes.
\end{remark}

The following two results  show that the  tests  \eqref{testglobal} and \eqref{testlokal} are consistent asymptotic level $\alpha$ tests.

\begin{proposition} \label{CorrGCP}
Under ${\bf H}_0(\eps)$ and ${\bf H}_0^{(z_0)}$, respectively, the tests \eqref{testglobal} and \eqref{testlokal}      have
 asymptotic level $\alpha$. More precisely,
 $$
\lim_{B \rightarrow \infty} \lim \limits_{n \rightarrow \infty} \mathbb P \Big( k_n^{1/2} \Db_n^{(\eps)}(1) \geq \hat q_{1- \alpha}^{(B)} (\Hb_n^{(\eps)}(1) ) \Big) = \alpha,
$$
if there exist $|\bar \taili|\geq \eps$, $\bar \kseqi \in (0,1)$ with $\int_0^{\bar \kseqi} g(y,\bar \taili) dy >0$, and
$$\lim \limits_{B \rightarrow \infty} \lim \limits_{n \rightarrow \infty} \mathbb P \Big( W_n^{(\taili_0)} \geq \hat q_{1 - \alpha}^{(B)}(W_n^{(\taili_0)}) \Big)  = \alpha,$$
if  there exists a $\bar \kseqi \in (0,1)$ with $\int_0^{\bar \kseqi} g(y,\taili_0) dy >0$.
\end{proposition}

\begin{proposition} \label{CorConsi}
The tests \eqref{testglobal} and \eqref{testlokal} are consistent in the following sense.
Under ${\bf H}_1(\eps)$ we have for all $B \in \mathbb N$ 
\begin{equation}
\label{neuA}
\lim \limits_{n \rightarrow \infty} \mathbb P \Big( k_n^{1/2}\Db_n^{(\eps)}(1) \geq \hat q_{1-\alpha}^{(B)} ( \Hb_n^{(\eps)}(1) ) \Big)  = 1.
\end{equation}
Under ${\bf H}_1^{\scriptscriptstyle (z_0)}$, we have for all $B \in \mathbb N$
$$
\lim \limits_{n \rightarrow \infty} \mathbb P \Big(W_n^{(z_0)} \geq \hat q_{1- \alpha}^{(B)}(W_n^{(z_0)}) \Big) =1.
$$
\end{proposition}


\section{Finite-sample properties} \label{sec5}
\def\theequation{5.\arabic{equation}}
\setcounter{equation}{0}

In this section we present the results of a simulation study, investigating the finite-sample properties of the new methodology for inference of gradual changes in the jump behaviour. The  design of this study is as follows.

\begin{itemize}
\item Each procedure is run $500$ times for any depicted combination of the involved constants. Furthermore in each run the number of bootstrap replications is $B=200$.
\item The estimators are calculated for $n=22,500$ data points with the combination of \textit{frequencies} $\Delta_n^{-1} = 450$, $225$, $90$ resulting in the choices $k_n =n \Delta_n= 50$, $100$, $250$ for the \textit{number of trading days}. For computational reasons, if not declared otherwise, we choose $n=10,000$ for each run of the tests \eqref{testglobal} and \eqref{testlokal} with \textit{frequencies} $\Delta_n^{-1} = 200$, $100$, $50$ corresponding to the \textit{number of trading days} $k_n = n \Delta_n = 50$, $100$, $200$.
\item We consider the following model for the transition kernel $g(y,d\taili)$ similar to Example \ref{ex2}:
\begin{align}
\label{SimModel}
g(y, \iz) = \begin{cases}
            \Big( \frac{\gamma(y)}{\pi \taili} \Big)^{1/2}, \quad &\text{ if } \taili > 0, \\
						0, \quad &\text{ if } \taili < 0,
						\end{cases}
\end{align}
with 
\begin{align}
\label{SimModel1}
\gamma(y) = \begin{cases}
            1, \quad &\text{ if } y \leq \gseqi_0, \\
						A(y- \gseqi_0)^w +1, \quad &\text{ if } y \geq \gseqi_0,
						\end{cases} \hspace{15mm} y \in [0,1]
\end{align}
for some $\gseqi_0 \in [0,1]$, $A > 0$ and $w >0$. In order to simulate pure jump It\=o semimartingale data according to such a gradual change we sample $15$ times more frequently and use a straight-forward modification of Algorithm 6.13 in \cite{ConTan04} to generate the increments $Y_j = \tilde X^{(j)}_{j \Delta_n /15} - \tilde X^{(j)}_{(j-1) \Delta_n /15}$ for $j = 1, \ldots, 15n$ of a $1/2$-stable pure jump L\'evy subordinator with characteristic exponent
\[
\Phi^{(j)}(u) = \int (e^{iu\taili} -1) \nu^{(j)}(d\taili),
\]
where $\nu^{(j)}(d\taili) = g(j/15n,d\taili)$. The resulting data vector $\{X_{\Delta_n}, \ldots, X_{n\Delta_n}\}$ is then given by
\[
X_{k\Delta_n} = \sum_{j =1}^{15k} Y_j
\]
for $k =1, \ldots, n$. 
\item In order to investigate the influence of a continuous component within the underlying It\=o semimartingale on the performance of our procedure we either use the plain data vector $\{X_{\Delta_n}, \ldots, X_{n\Delta_n}\}$ or $\{(X+S)_{\Delta_n}, \ldots, (X+S)_{n\Delta_n}\}$, where $S_t = W_t +t$.
\item For computational reasons the supremum of the tail parameter $\taili$ over $M_\eps = (-\infty,-\eps] \cup [\eps, \infty)$ in each statistic is approximated by taking the maximum over a finite grid. In the pure jump case we use $M=\{ 0.1$, $0.15$, $0.25$, $1$, $2\}$, resulting in $\eps = 0.1$. In the case including a continuous component we consider $M=\{j\sqrt{\Delta_n} \mid j = 2$, $3.5$, $5$, $6.5$, $8\}$, resulting in $\eps = 2 \sqrt{\Delta_n}$. In the latter case, we choose $\eps$ depending on $\sqrt{\Delta_n}$ since jumps of smaller size may be dominated by the Brownian component leading to a loss of efficiency of the statistical procedures.
\end{itemize}

\subsection{Finite-sample properties of the estimator \( \hat \theta^{(\varepsilon)}_n \)}

We implement our estimation method as follows:
\begin{enumerate}
\item[\textit{Step 1.}] Choose a preliminary estimate \(\hat \gseqi^{(pr)} \in (0,1)\), a probability level \(\alpha \in (0,1) \) and a parameter \(r \in (0,1]\).
\item[\textit{Step 2.}] Initial choice of the tuning parameter \(\thrle_n\): \\
       Evaluate \eqref{thrledefeq} for $\hat \gseqi^{(pr)}, \alpha$ and $r$ (with \(\eps\) and \( B\) as described above) and obtain \( \hat\thrle^{(in)}\).
\item[\textit{Step 3.}] Intermediate estimate of the change point. \\
   Evaluate \eqref{Estimator} for \(\hat\thrle^{(in)}\) and obtain \(\hat \gseqi^{(in)}\).
\item[\textit{Step 4.}] Final choice of the tuning parameter \(\thrle_n\): \\
   Evaluate \eqref{thrledefeq} for $\hat \gseqi^{(in)}, \alpha, r$ and obtain \( \hat\thrle^{(fi)}\).
\item[\textit{Step 5.}] Estimate $\gseqi_0$. \\
    Evaluate \eqref{Estimator} for \(\hat\thrle^{(fi)}\) and obtain the final estimate \(\hat \gseqi\) of the change point.
\end{enumerate}

Furthermore in order to measure the quality of the resulting estimates $\Theta = \{\hat \gseqi_1, \ldots, \hat \gseqi_{500}\}$ we subsequently use the mean absolute deviation to the true value $\gseqi_0$, that is
\[
\ell^1(\Theta, \gseqi_0) := \frac 1{500} \sum \limits_{j=1}^{500} |\hat \gseqi_j - \gseqi_0|.
\]

\begin{figure}[t!]
\centering
\includegraphics[width=0.48\textwidth]{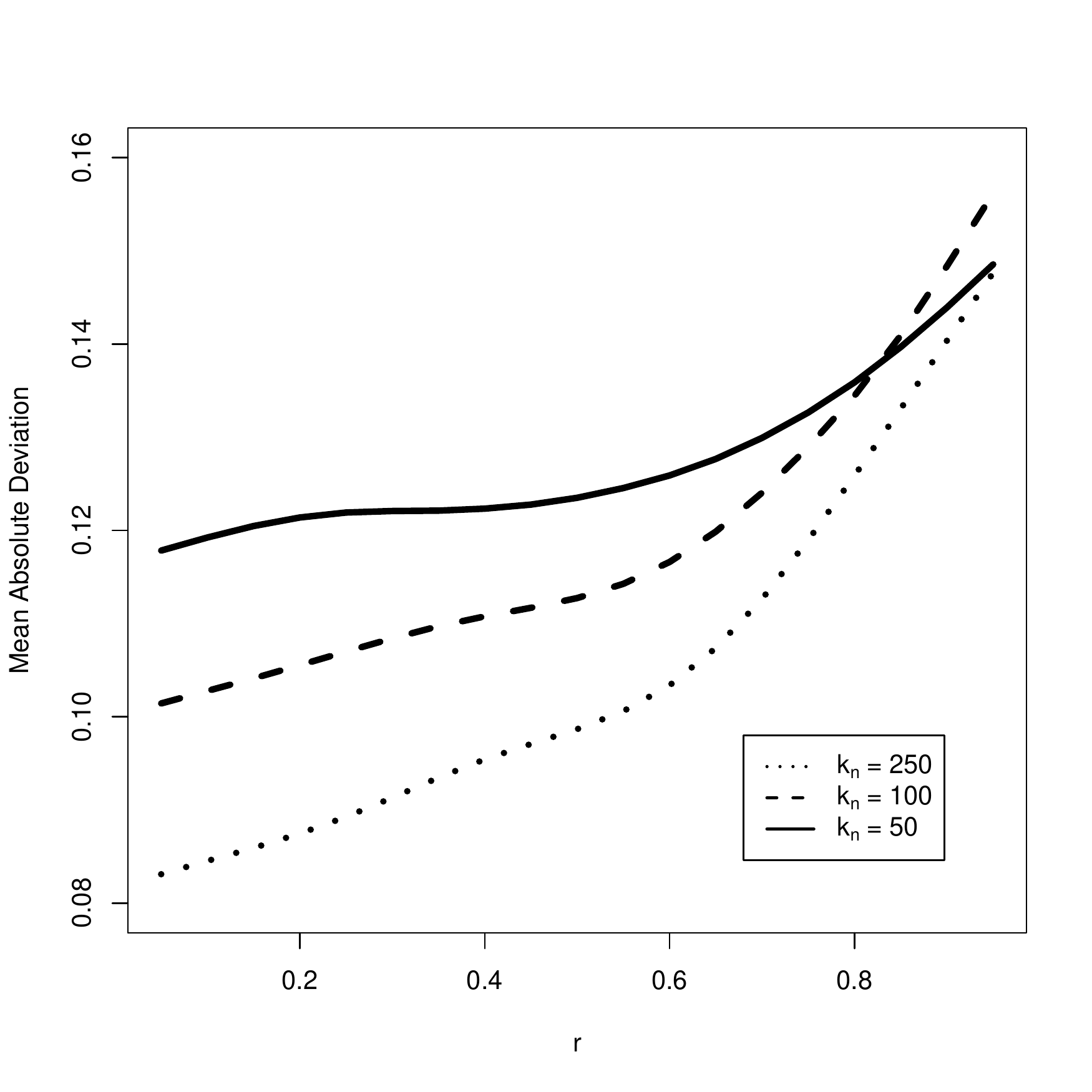}~~~ \includegraphics[width=0.48\textwidth]{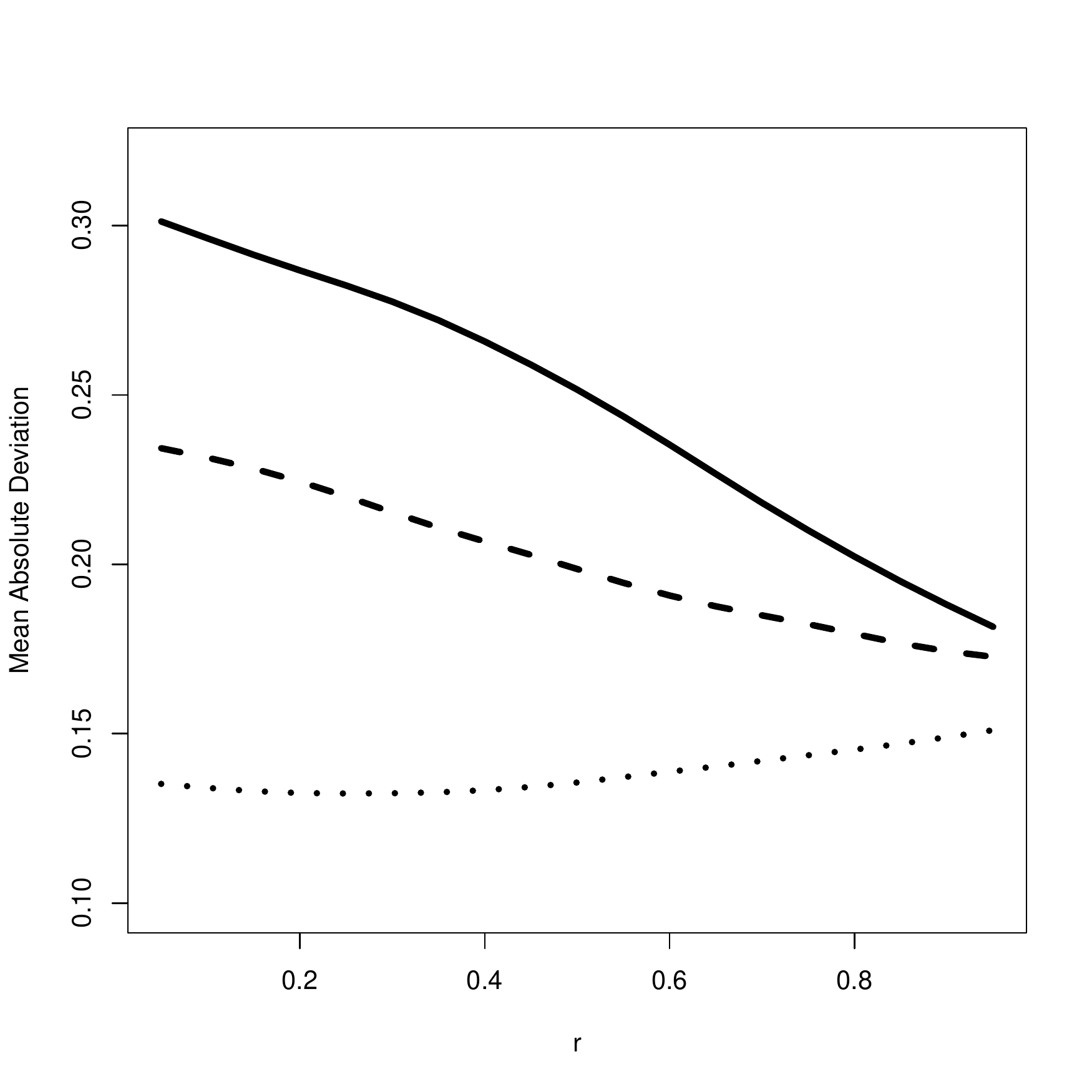}
\caption{\label{rDurchl}
\it Mean absolute deviation for different choices of the parameter $r$ for pure jump data (left-hand side) and with an additional continuous component (right-hand side).
}
\end{figure}

The results presented in Figure \ref{rDurchl} show the mean absolute deviation for varying parameter $r \in (0,1]$. 
Here and below we use the probability level $\alpha = 0.1$. As a preliminary estimate we choose \(\hat \gseqi^{(pr)} = 0.1\), whereas the true change point is located at \(\gseqi_0 = 0.4\). Moreover we simulate a linear change, that is we have $w=1$ in model \eqref{SimModel} while the constant $A$ in 
\eqref{SimModel1} is chosen such that the characteristic quantity $\Dc^{\scriptscriptstyle (\eps)}(1)$ for a gradual change satisfies $\Dc^{\scriptscriptstyle (\eps)}(1)=5$. From the left panel it becomes apparent that in the pure jump case the mean absolute deviation is increasing in $r$. Consequently we choose $r=0.01$ for pure jump It\=o semimartingale data in all following considerations. The behaviour for data including a continuous component (right panel) is different. For the choices $k_n =50$ and $k_n =100$ it is decreasing while $k_n =250$ results in a curve which is nearly constant. Thus we choose $r=1$ in the following investigations whenever a continuous component is involved. 

\begin{figure}[t!]
\centering
\includegraphics[width=0.48\textwidth]{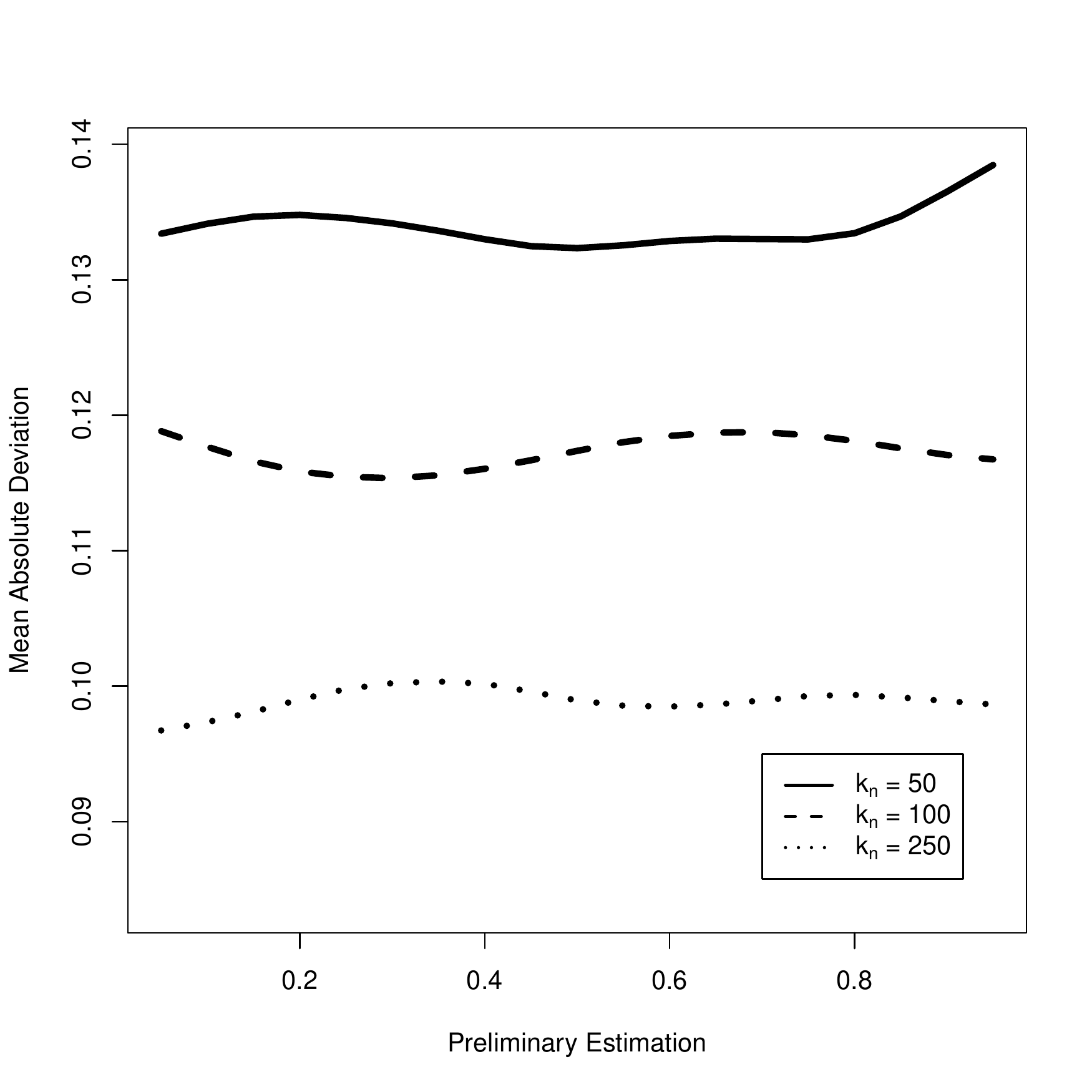}~~~ \includegraphics[width=0.48\textwidth]{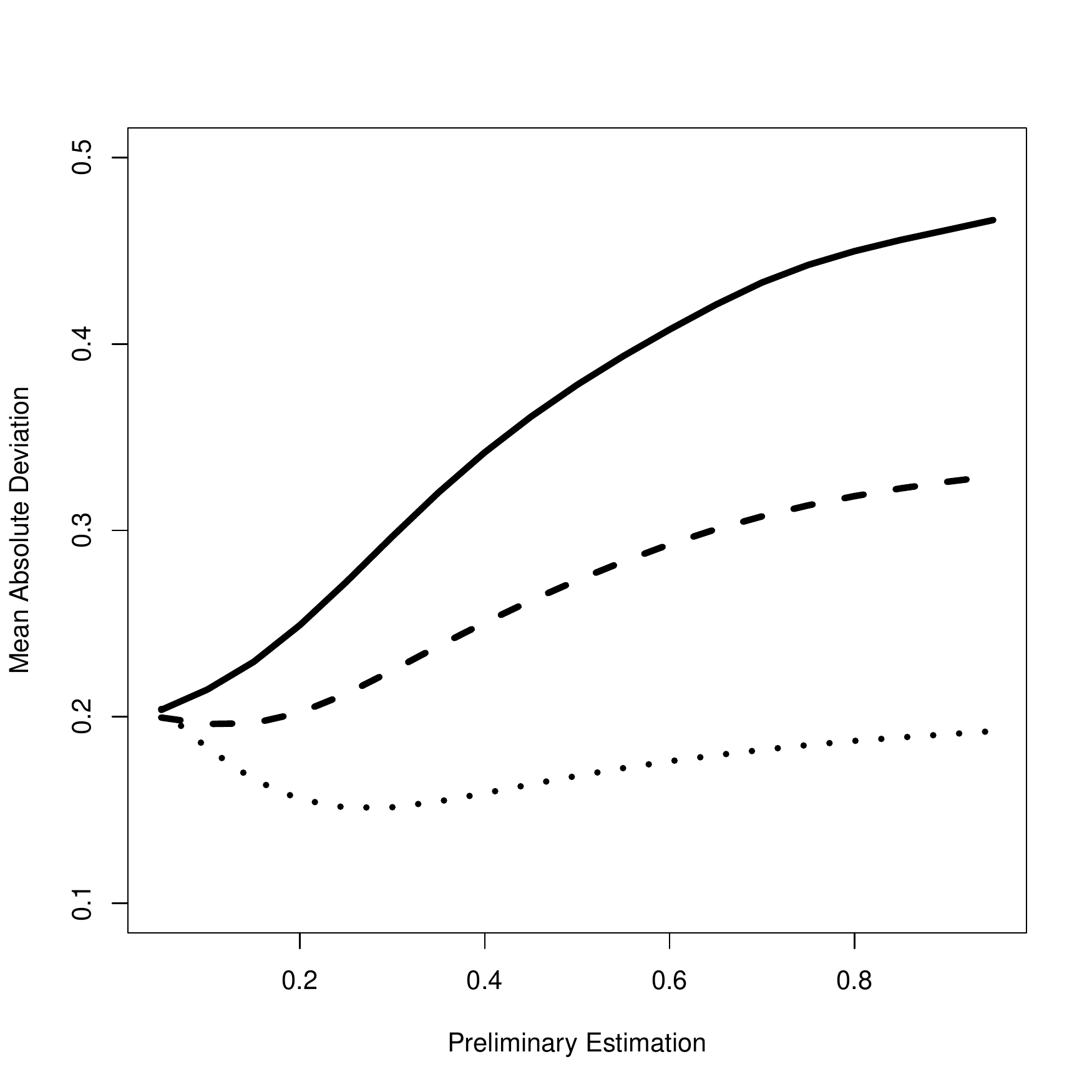}
\caption{\label{PEDuLDu}
\it Mean absolute deviation of the estimator $\hat \theta^{(\varepsilon)}_n$ for different choices of the preliminary estimate \(\hat \gseqi^{(pr)}\) for pure jump It\=o semimartingales (left panel) and with a continuous component (right panel).
}
\end{figure}

Figure \ref{PEDuLDu} shows the mean absolute deviation of the estimator $\hat \theta^{(\varepsilon)}_n$ for different choices of the preliminary estimate \(\hat \gseqi^{(pr)}\). Here the change is again linear $(w=1)$ and it is located  at $\gseqi_0 = 0.4$, while the choice of the constant $A$ in \eqref{SimModel1} corresponds to $\Dc^{\scriptscriptstyle (\eps)}(1) =3$. In the pure jump case (left-hand side) the performance of the estimator is nearly independent of the preliminary estimate. However, it can be seen from the right panel that the mean absolute deviation becomes minimal  for small choices of $\hat \gseqi^{(pr)}$. These findings are  confirmed by a further simulation study, which is not depicted here for the sake of brevity. Our results show that large choices of the preliminary estimate yield an over-estimation, which can be explained by the fact that larger values of \(\hat \gseqi^{(pr)}\) induce  larger tuning parameters $\thrle_n$ in the estimates in Steps 2--5 as well. 

\begin{figure}[t!]
\centering
\includegraphics[width=0.48\textwidth]{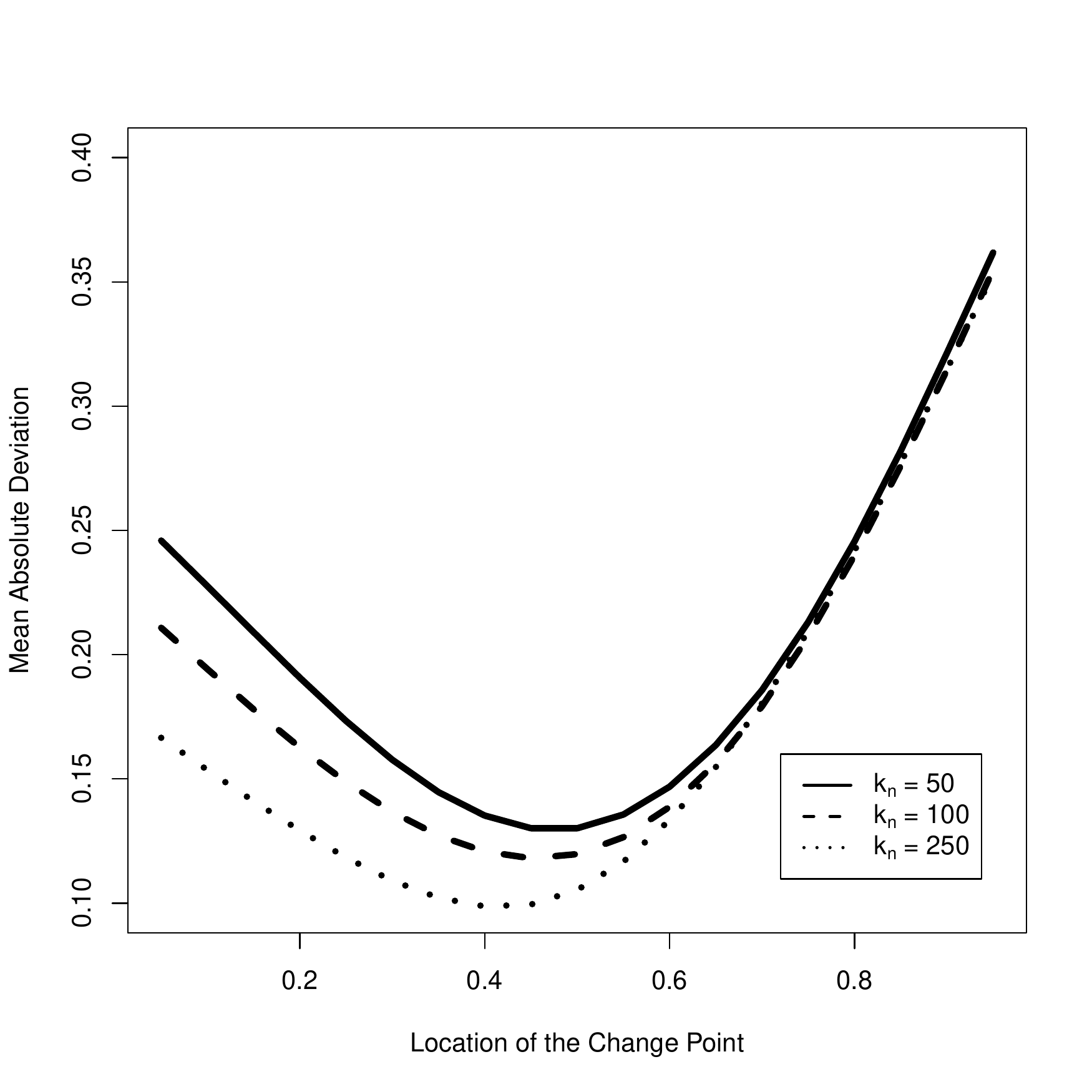}~~~ \includegraphics[width=0.48\textwidth]{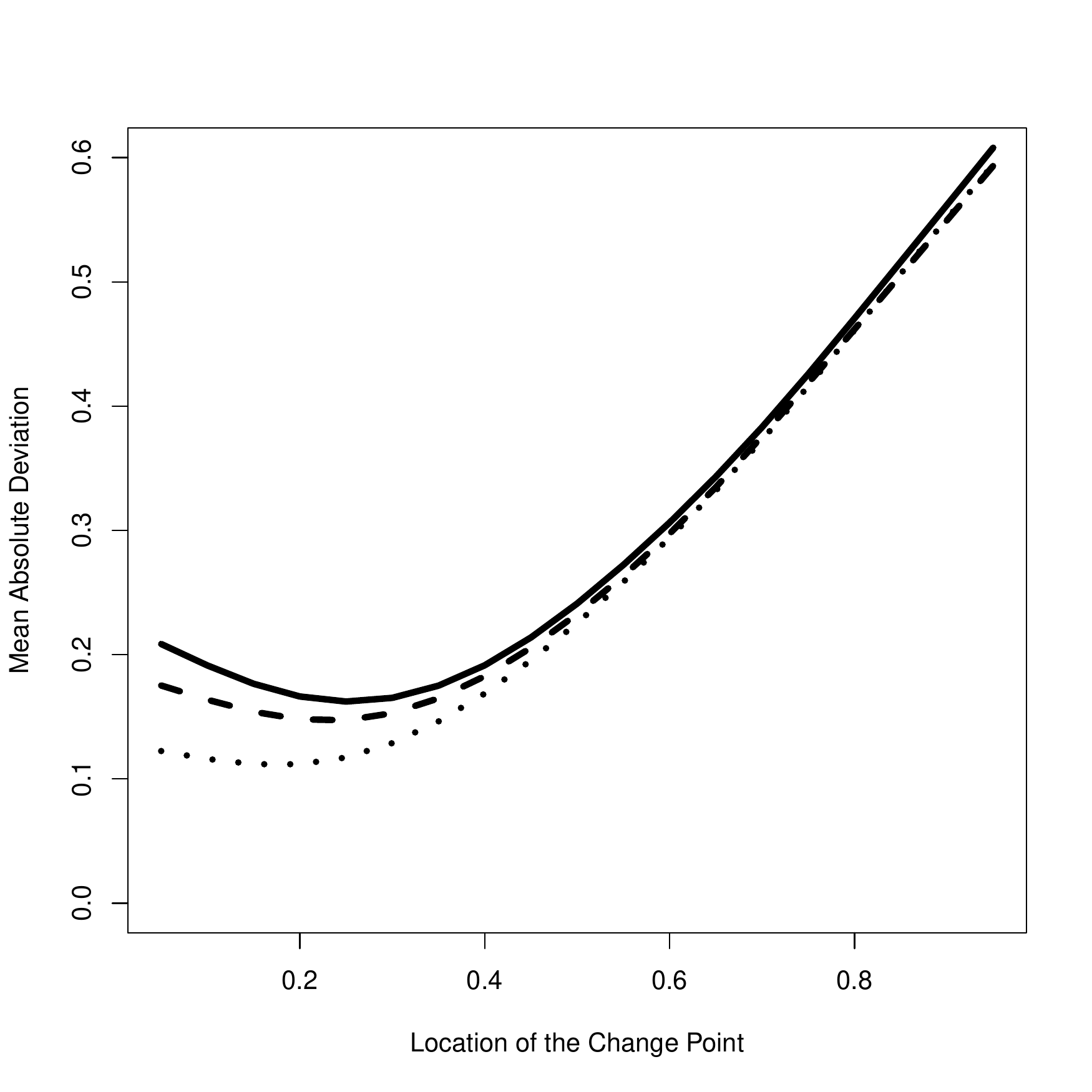}
\caption{\label{Fig:LDurchl}
\it Mean absolute deviation for different locations of the change point $\gseqi_0$ for pure jump data (left-hand side) and with an additional continuous component (right-hand side).
}
\end{figure}

In Figure \ref{Fig:LDurchl} we display the performance of the estimation procedure for different locations of the change point. The change is linear $(w=1)$
and the constant $A$  in \eqref{SimModel1} is chosen appropriately such that $\Dc^{\scriptscriptstyle (\eps)}(1) =3$ holds in each scenario. Furthermore the preliminary estimate is chosen as \(\hat \gseqi^{(pr)} = 0.1\). The left panel suggests that a change point can be estimated best, if it is located around $\gseqi_0 =0.5$. This result corresponds to the findings in Figure 2 in \cite{BueHofVetDet15}, who demonstrated that  under the alternative of an abrupt change the power of the classical CUSUM test is maximal for a change point around $\gseqi_0 =0.5$. However, from the right panel in Figure \ref{Fig:LDurchl} it is clearly visible that, if an additional continuous component is present, large values of the change point lead to a large estimation error. This is a consequence of the shape of the model \eqref{SimModel}: For $\gseqi_0$ close to $1$ the behaviour of the underlying process is similar to the null hypothesis.

\begin{figure}[t!]
\centering
\includegraphics[width=0.48\textwidth]{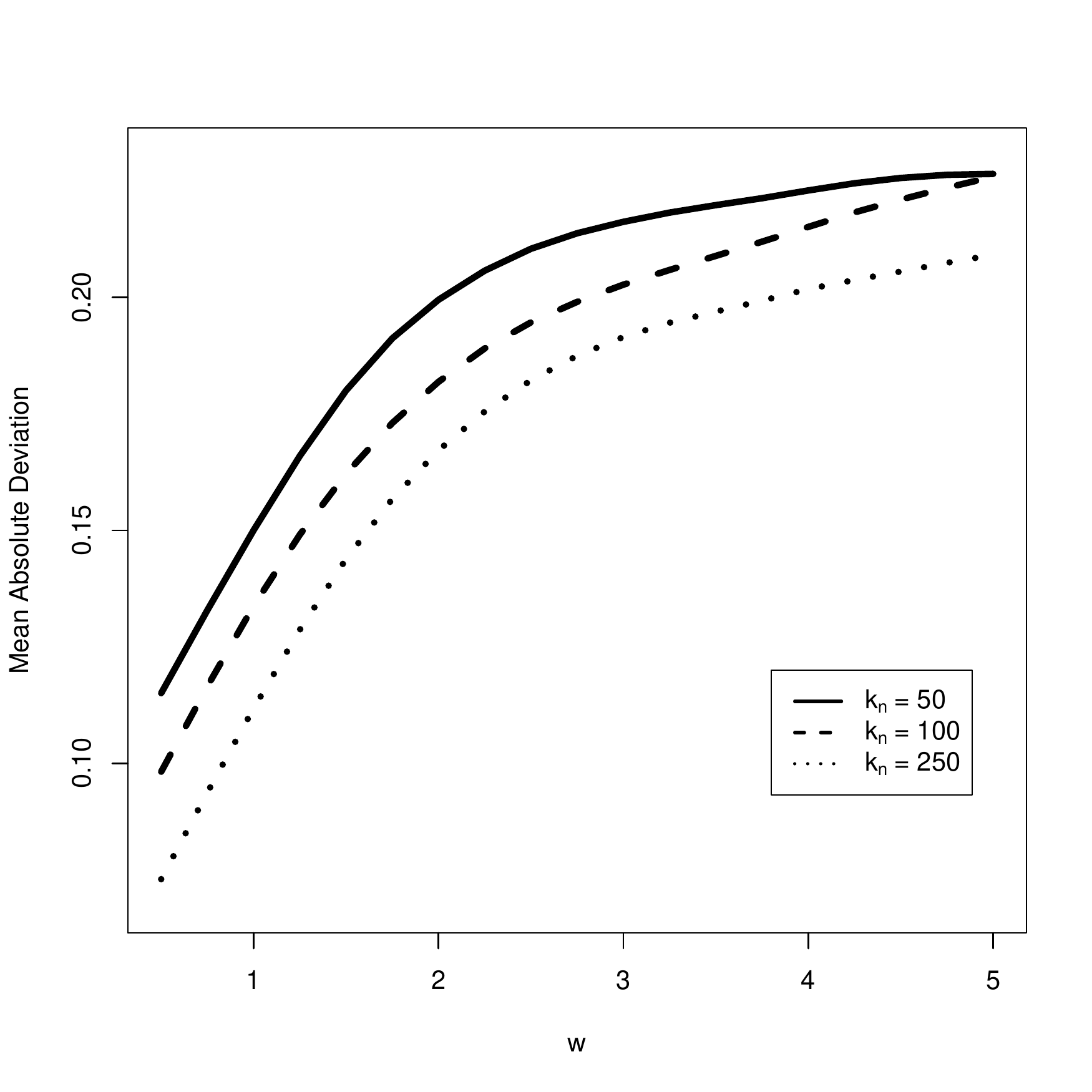}~~~ \includegraphics[width=0.48\textwidth]{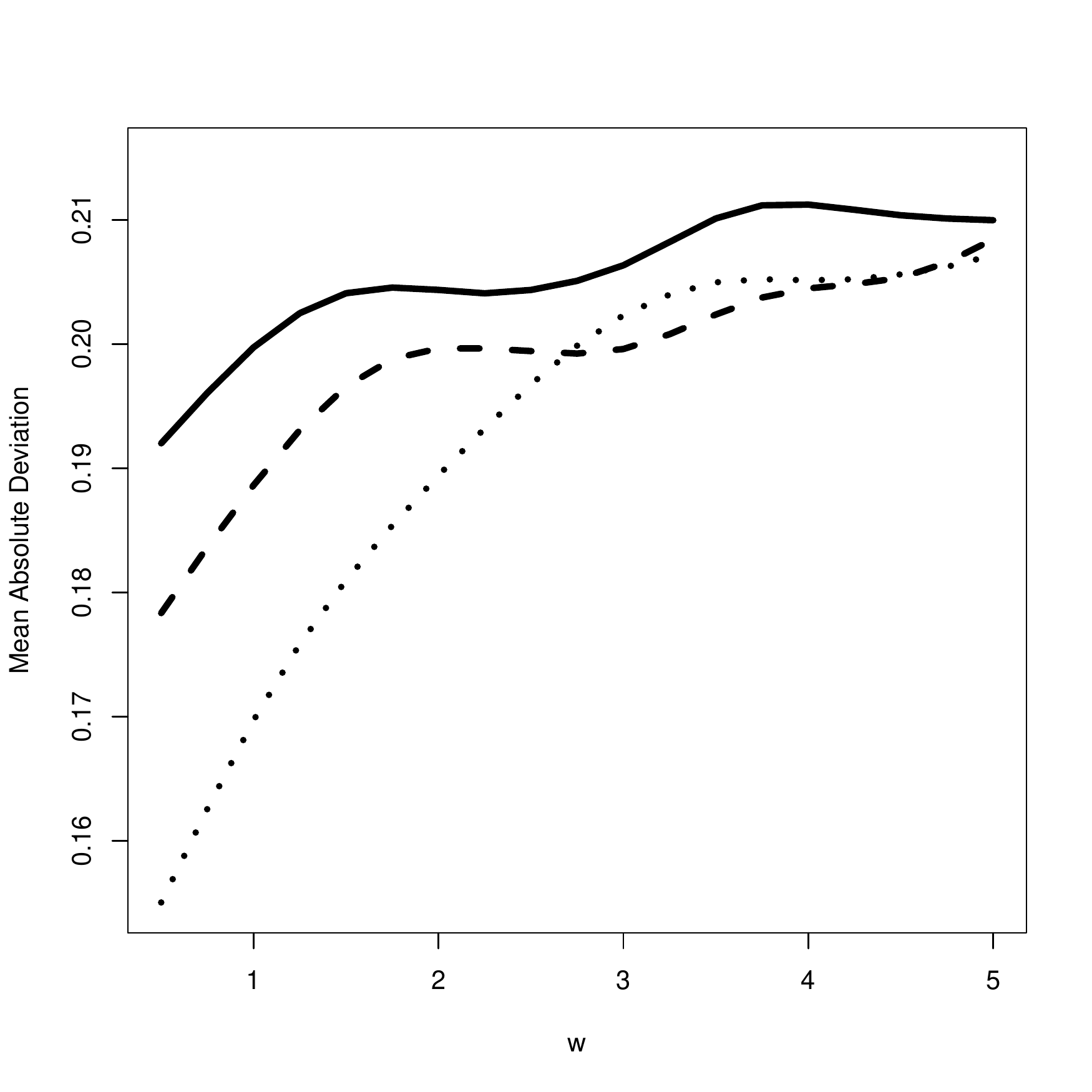}
\caption{\label{sDurchl}
\it Mean absolute deviation of the estimator $\hat \theta^{(\varepsilon)}_n$ for different degrees of smoothness $w$ of the change for pure jump It\=o semimartingales (left panel) and with an additional continuous component (right panel).
}
\end{figure}

In Figure \ref{sDurchl} we depict the results for different degrees of smoothness $w$ in \eqref{SimModel1}. The true location of the change point is $\gseqi_0 = 0.4$, while we choose \(\hat \gseqi^{(pr)} = 0.1\) for the preliminary estimate and in order to keep the results comparable the constant $A$ in \eqref{SimModel1} is chosen such that $\Dc^{\scriptscriptstyle (\eps)}(1) =2$.
Notice that the graphic on the right-hand side has a different scale of the $y$-axis, such that we obtain a slightly higher estimation error as well if a Brownian component and an effective drift are present. The results in Figure \ref{sDurchl} are as expected: The higher the degree of smoothness, the more difficult a break point can be detected resulting in a larger estimation error. 

\subsection{Finite-sample properties of the test procedures}

In order to investigate the finite-sample properties of the test procedures \eqref{testglobal} and \eqref{testlokal} we choose a level of significance $\alpha = 5\%$ in each of the following simulations. Table \ref{tab:h01} and \ref{tab:h02} contain the relative frequencies of rejections of both tests where the sample size is $n=22,500$ under \({\bf H}_0 (\eps)\) and \({\bf H}_0^{\scriptscriptstyle (z_0)}\), respectively, that is for $\gseqi_0 =1$ in \eqref{SimModel1}. For the test \eqref{testglobal} the supremum of the tail index $\taili$ is again approximated by the maximum over the finite grid $M=\{0.1$, $0.15$, $0.25$, $1$, $2\}$ in the pure jump case and over the finite grid $M=\{j \sqrt{\Delta_n} \mid j\in \{ 2$, $3.5$, $5$, $6.5$, $8\}\}$, if a continuous component is involved. We observe a reasonable approximation of the nominal level in all cases under consideration.

\begin{table}[t!]
\begin{center}
\begin{tabular}{ c ||c||c|c|c|c|c| }
\hline
\multicolumn{1}{|c||}{ } & \multicolumn{1}{c||}{Test \eqref{testglobal}} & \multicolumn{5}{c|}{Test \eqref{testlokal}} \\
\hline
\multicolumn{1}{|c||}{$k_n$} & \multicolumn{1}{c||}{M} & 
\multicolumn{1}{c|}{$z_0 = 0.1$} & \multicolumn{1}{c|}{$z_0 = 0.15$} & \multicolumn{1}{c|}{$z_0 = 0.25$} & 
\multicolumn{1}{c|}{$z_0 = 1$} & \multicolumn{1}{c|}{$z_0 = 2$} \\
\hline
\multicolumn{1}{|c||}{$50$} & \multicolumn{1}{c||}{0.062} &  
\multicolumn{1}{c|}{$0.066$} & \multicolumn{1}{c|}{$0.054$} & \multicolumn{1}{c|}{$0.068$} & 
\multicolumn{1}{c|}{$0.052$} & \multicolumn{1}{c|}{$0.054$} \\
\hline
\multicolumn{1}{|c||}{75} & \multicolumn{1}{c||}{0.064} &  
\multicolumn{1}{c|}{$0.066$} & \multicolumn{1}{c|}{$0.046$} & \multicolumn{1}{c|}{$0.066$} & 
\multicolumn{1}{c|}{$0.066$} & \multicolumn{1}{c|}{$0.066$} \\
\hline
\multicolumn{1}{|c||}{100} & \multicolumn{1}{c||}{0.062} & 
\multicolumn{1}{c|}{$0.066$} & \multicolumn{1}{c|}{$0.058$} & \multicolumn{1}{c|}{$0.064$} & 
\multicolumn{1}{c|}{$0.056$} & \multicolumn{1}{c|}{$0.058$} \\
\hline
\multicolumn{1}{|c||}{150} & \multicolumn{1}{c||}{0.056} & 
\multicolumn{1}{c|}{$0.054$} & \multicolumn{1}{c|}{$0.058$} & \multicolumn{1}{c|}{$0.078$} & 
\multicolumn{1}{c|}{$0.074$} & \multicolumn{1}{c|}{$0.064$} \\
\hline
\multicolumn{1}{|c||}{250} & \multicolumn{1}{c||}{0.060} &
\multicolumn{1}{c|}{$0.052$} & \multicolumn{1}{c|}{$0.066$} & \multicolumn{1}{c|}{$0.070$} & 
\multicolumn{1}{c|}{$0.052$} & \multicolumn{1}{c|}{$0.042$} \\
\hline
\end{tabular}
\caption{\it Simulated rejection probabilities of the test \eqref{testglobal} and the test \eqref{testlokal}, using $500$ pure jump It\=o semimartingale data vectors under the null hypotheses \({\bf H}_0 (\eps)\) and \({\bf H}_0^{\scriptscriptstyle (z_0)}\), respectively.} 
\label{tab:h01}
\end{center}
\vspace{-.1cm}
\end{table}

\begin{table}[t!]
\begin{center}
\begin{tabular}{ c ||c||c|c|c|c|c| }
\hline
\multicolumn{1}{|c||}{ } & \multicolumn{1}{c||}{Test \eqref{testglobal}} & \multicolumn{5}{c|}{Test \eqref{testlokal}} \\
\hline
\multicolumn{1}{|c||}{$k_n$} & \multicolumn{1}{c||}{M} & 
\multicolumn{1}{c|}{$z_0 = 2\sqrt{\Delta_n}$} & \multicolumn{1}{c|}{$z_0 = 3.5\sqrt{\Delta_n}$} & \multicolumn{1}{c|}{$z_0 = 5\sqrt{\Delta_n}$} & 
\multicolumn{1}{c|}{$z_0 = 6.5\sqrt{\Delta_n}$} & \multicolumn{1}{c|}{$z_0 = 8\sqrt{\Delta_n}$} \\
\hline
\multicolumn{1}{|c||}{$50$} & \multicolumn{1}{c||}{0.056} &  
\multicolumn{1}{c|}{$0.056$} & \multicolumn{1}{c|}{$0.066$} & \multicolumn{1}{c|}{$0.062$} & 
\multicolumn{1}{c|}{$0.064$} & \multicolumn{1}{c|}{$0.058$} \\
\hline
\multicolumn{1}{|c||}{75} & \multicolumn{1}{c||}{0.042} &  
\multicolumn{1}{c|}{$0.042$} & \multicolumn{1}{c|}{$0.038$} & \multicolumn{1}{c|}{$0.044$} & 
\multicolumn{1}{c|}{$0.038$} & \multicolumn{1}{c|}{$0.044$} \\
\hline
\multicolumn{1}{|c||}{100} & \multicolumn{1}{c||}{0.042} & 
\multicolumn{1}{c|}{$0.042$} & \multicolumn{1}{c|}{$0.054$} & \multicolumn{1}{c|}{$0.042$} & 
\multicolumn{1}{c|}{$0.070$} & \multicolumn{1}{c|}{$0.046$} \\
\hline
\multicolumn{1}{|c||}{150} & \multicolumn{1}{c||}{0.064} & 
\multicolumn{1}{c|}{$0.064$} & \multicolumn{1}{c|}{$0.066$} & \multicolumn{1}{c|}{$0.062$} & 
\multicolumn{1}{c|}{$0.058$} & \multicolumn{1}{c|}{$0.060$} \\
\hline
\multicolumn{1}{|c||}{250} & \multicolumn{1}{c||}{0.066} &
\multicolumn{1}{c|}{$0.066$} & \multicolumn{1}{c|}{$0.062$} & \multicolumn{1}{c|}{$0.046$} & 
\multicolumn{1}{c|}{$0.048$} & \multicolumn{1}{c|}{$0.060$} \\
\hline
\end{tabular}
\caption{\it Simulated rejection probabilities of the test \eqref{testglobal} and the test \eqref{testlokal}, using $500$ pure jump It\=o semimartingale data vectors plus an effective drift and a Brownian motion under the null hypotheses \({\bf H}_0 (\eps)\) and \({\bf H}_0^{\scriptscriptstyle (z_0)}\), respectively.} 
\label{tab:h02}
\end{center}
\vspace{-.1cm}
\end{table}

\begin{figure}[t!]
\centering
\includegraphics[width=0.48\textwidth]{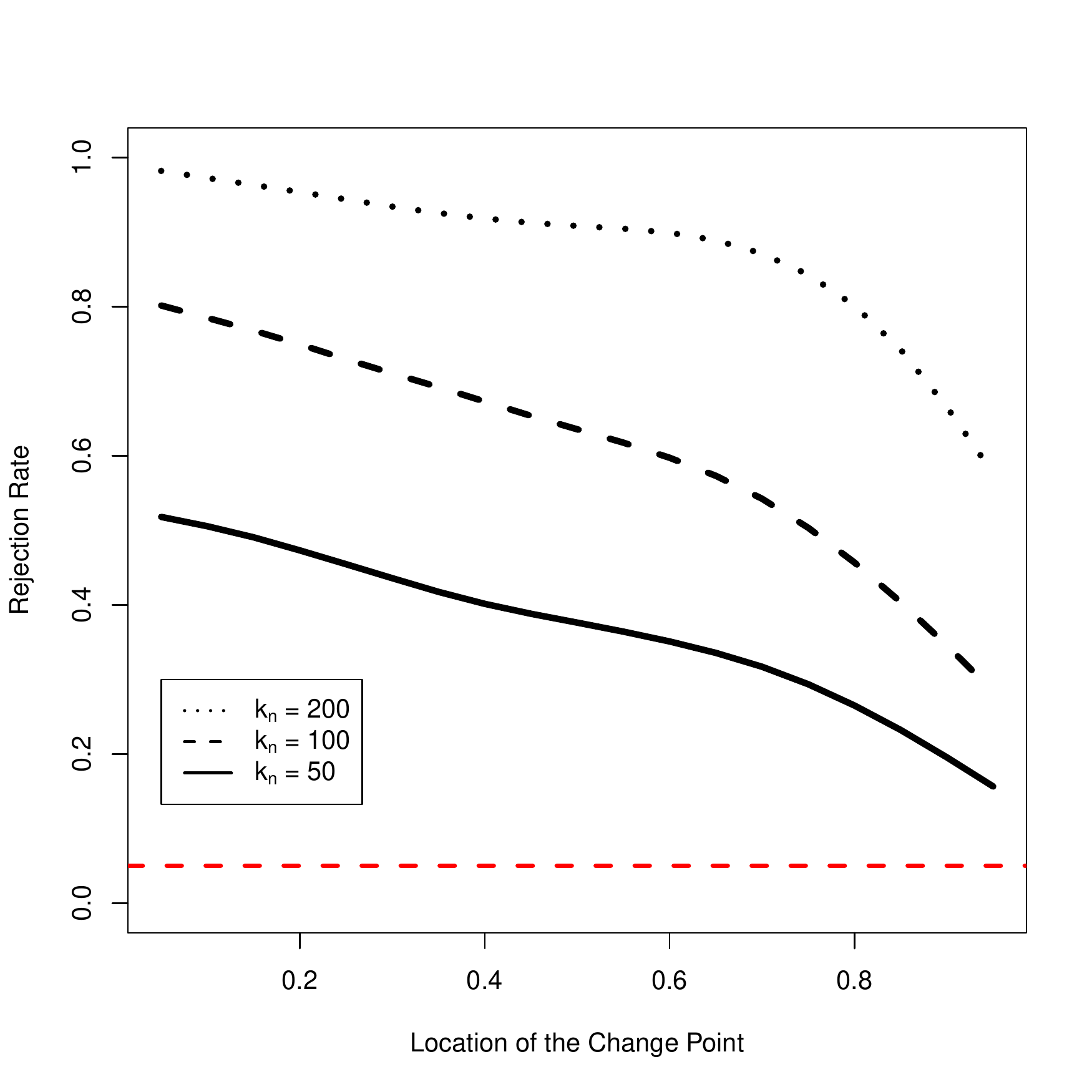}~~~ \includegraphics[width=0.48\textwidth]{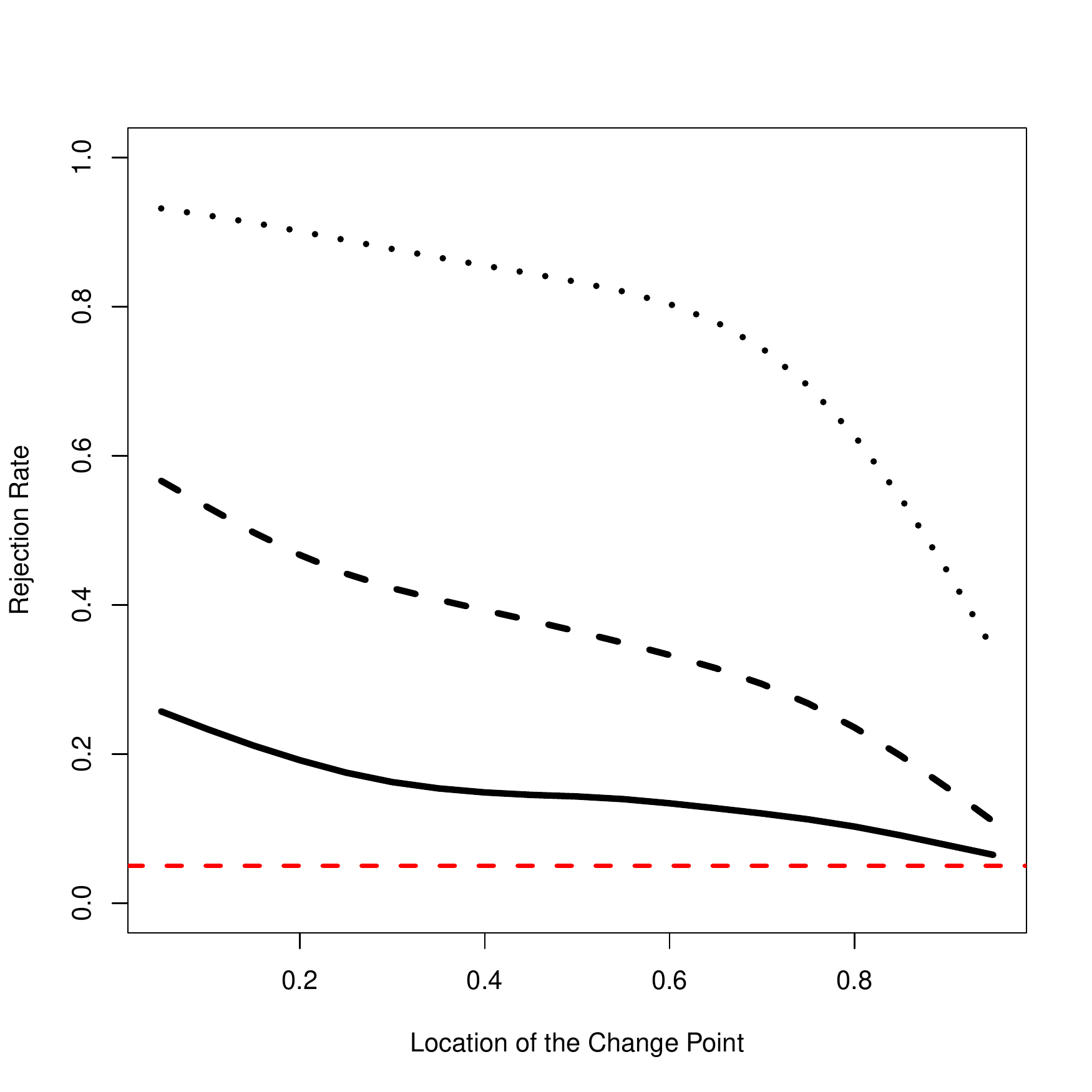}
\caption{\label{H1Lvz0v}
\it Simulated rejection probabilities of the test \eqref{testglobal} for different locations of the change point $\gseqi_0$ for pure jump It\=o semimartingale data (left-hand side) and with an additional continuous component (right-hand side). The dashed red line indicates the nominal level $\alpha = 5\%$.
}
\end{figure}

Figure \ref{H1Lvz0v} shows the simulated rejection probabilities of the test \eqref{testglobal} for a linear change $(w=1)$ at different locations of the change point $\gseqi_0$. The constant $A$ in \eqref{SimModel1} is chosen such that $\Dc^{\scriptscriptstyle (\eps)}(1) =0.8$. In this special case the presence of a continuous component leads to a relatively strong loss of power of the test. This can be explained by the small value of $\Dc^{\scriptscriptstyle (\eps)}(1)$, for which the rather small change in the jump behaviour is predominated stronger by the Brownian component. Furthermore the power of the test is decreasing in $\gseqi_0$, which is a consequence of the fact that for large $\gseqi_0$ the data is ''closer'' to the null hypothesis (recall the right-hand side of Figure \ref{Fig:LDurchl}). 

\begin{figure}[t!]
\centering
\includegraphics[width=0.48\textwidth]{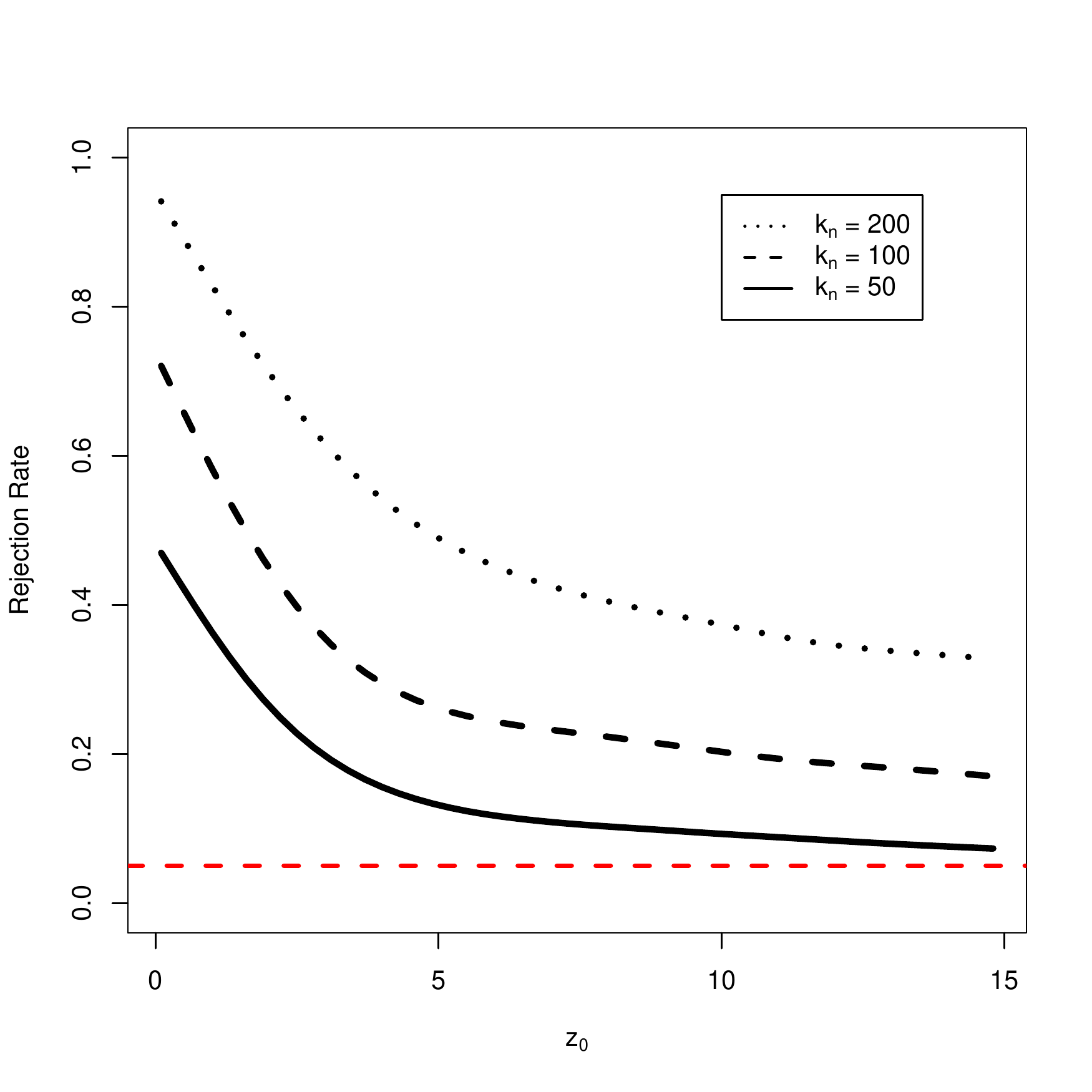}~~~ \includegraphics[width=0.48\textwidth]{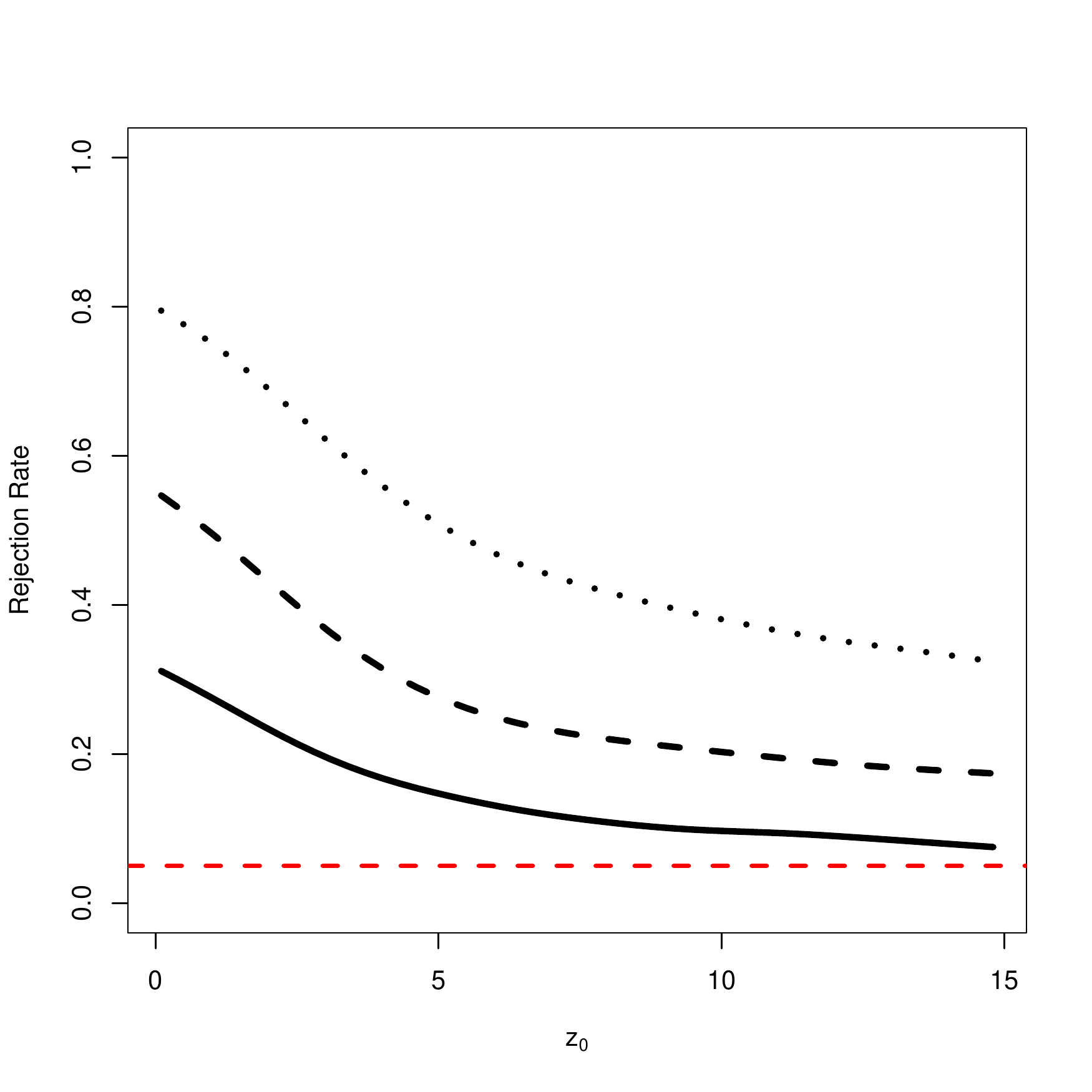}
\caption{\label{Fig:H1z0var}
\it Simulated rejection probabilities of the test \eqref{testlokal} for different choices of $\taili_0$ for pure jump data (left panel) and plus an effective drift and a Brownian motion (right panel). The dashed red line indicates the nominal level $\alpha = 5\%$.
}
\end{figure}

Figure \ref{Fig:H1z0var} reveals the behaviour of the test \eqref{testlokal} for different choices of the tail parameter $\taili_0$. The data exhibits a linear change at $\gseqi_0 = 0.4$ with $A=5$. The curves are obviously decreasing in $\taili_0$, which is a consequence of very few large jumps even for this rather large choice of $A=5$.

\begin{figure}[t!]
\centering
\includegraphics[width=0.48\textwidth]{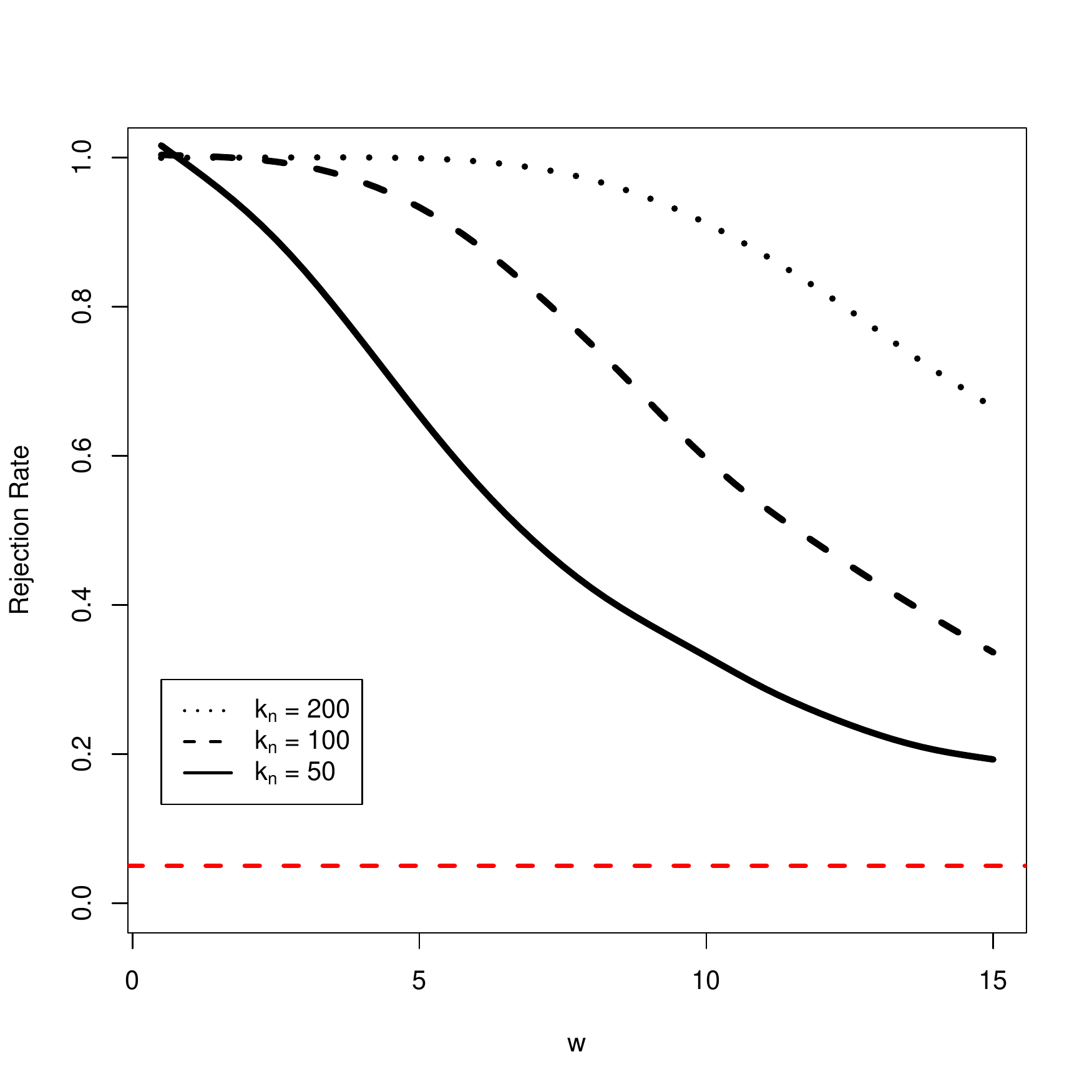}~~~ \includegraphics[width=0.48\textwidth]{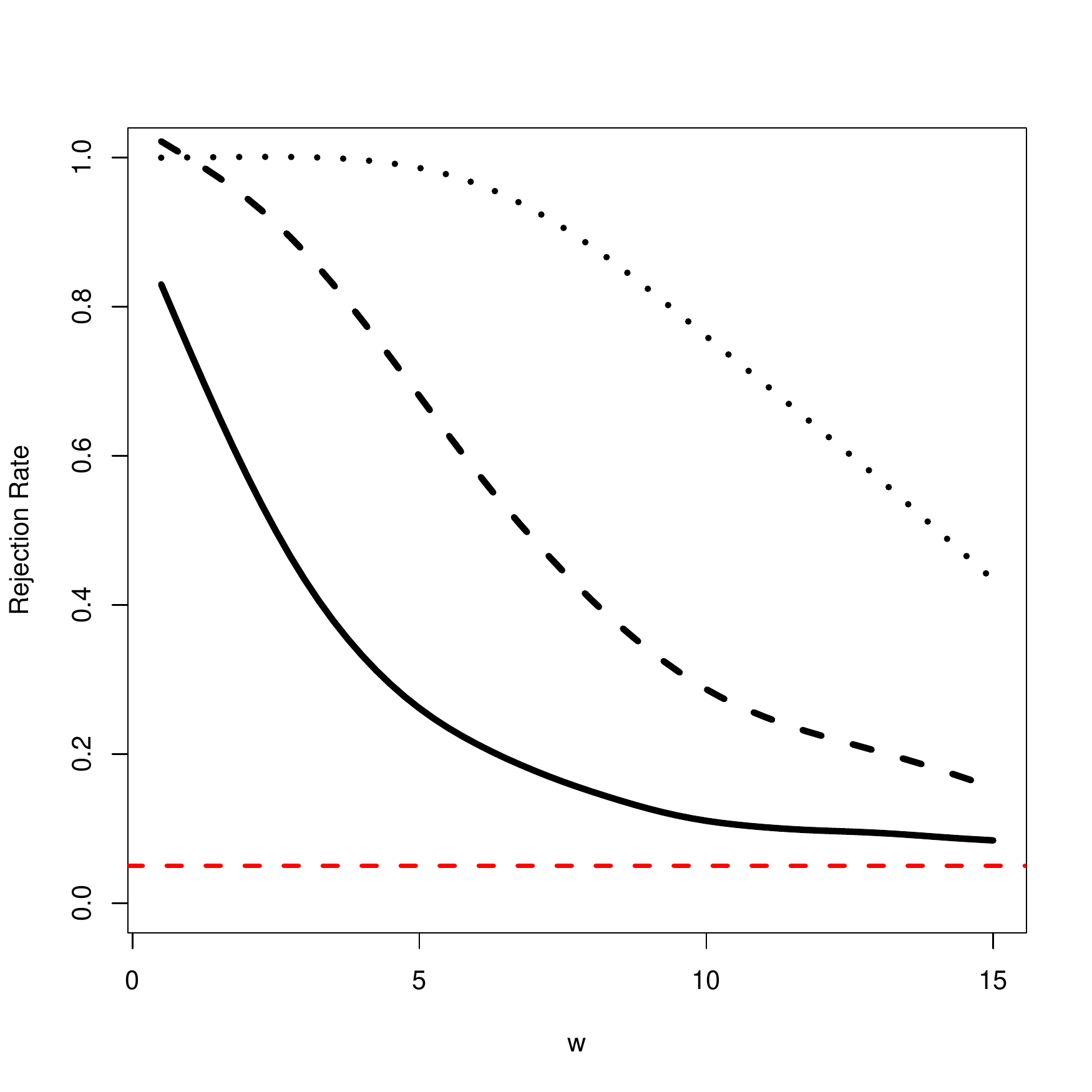}
\caption{\label{H1sDurchl}
\it Simulated rejection probabilities of the test \eqref{testglobal} for different degrees of smoothness $w$ for pure jump data (left-hand side) and with an additional continuous part (right-hand side). The dashed red line indicates the nominal level $\alpha = 5\%$.
}
\end{figure}

Finally, Figure \ref{H1sDurchl} shows the dependence of the test \eqref{testglobal} on the degree of smoothness $w$. The change is located at $\gseqi_0=0.4$, and in order to keep the results comparable the constant $A$  in \eqref{SimModel1} is chosen  such that $\Dc^{\scriptscriptstyle (\eps)}(1) =3$ in each scenario. Overall, the power of the test in this experiment is relatively high, which is a consequence of the intermediate choice $\gseqi_0 = 0.4$ (recall Figure \ref{Fig:LDurchl} and Figure \ref{H1Lvz0v}). From the graphic on the right-hand side we also observe that the presence of a continuous component leads to a loss of power in this case as well. Moreover, the fact that the curves in Figure \ref{H1sDurchl} are decreasing in $w$ coincides with the heuristic idea that smoother changes are more difficult to detect.

\medskip
\medskip

{\bf Acknowledgements}
This work has been supported in part by the
Collaborative Research Center "`Statistical modeling of nonlinear
dynamic processes"' (SFB 823, Projects A1 and C1) of the German Research Foundation (DFG).
We  would like to thank Martina  Stein who typed parts of this manuscript with considerable
technical expertise and Michael Vogt for some helpful discussions.

\bibliographystyle{apalike}
\setstretch{1.25}
\setlength{\bibsep}{1pt}
\begin{small} 
\bibliography{biblio}
\end{small}

\section{Proofs and technical details}
\def\theequation{6.\arabic{equation}}
\setcounter{equation}{0}
\label{sec:AuxRes}

The following assumptions will be used frequently in the sequel.

\begin{condition}
\label{Cond1}
For each $n \in \N$ let $X^{(n)}$ be an It\=o semimartingale  of the form \eqref{ItoSem} with characteristics $(b^{(n)}_s, \sigma^{(n)}_s, \nu^{(n)}_s)$ and the following properties:
\begin{enumerate}[(a)]
\item There exists a $g \in \Gc$ such that 
       \begin{align*}
       \nu^{(n)}_{s}(dz) = g\Big(\frac s{n\Delta_n}, dz\Big)
        \end{align*}
				holds for all $s \in [0,n\Delta_n]$ and all $n \in \N$ as measures on $(\R,\Bb)$.
\item The drift $b^{(n)}_s$ and the volatility $\sigma^{(n)}_s$ are deterministic and Borel measurable functions on $\R_+$. Moreover, these
      functions are uniformly bounded in $s \in \R_+$ and $n \in \N$.
\item The observation scheme $\{ X^{(n)}_{i \Delta_n} \mid i=0,\ldots,n \}$ satisfies
      $\Delta_n \rightarrow 0 $ and $  n \Delta_n \rightarrow \infty.$
\end{enumerate}
\end{condition}

We begin with an auxiliary result which is a generalization of Lemma 2 in \cite{RueWoe02}. Throughout this section $K$ denotes a generic constant which typically changes from line to line and may depend on certain bounds and parameters, but not on $n$. 
 
\begin{lemma}
\label{TransProbAbs}
Let $T>0$ and let $Y$ be an It\=o semimartingale with $Y_0=0$ having a representation as in \eqref{ItoSem} with characteristics $(b_s,\sigma_s,\nu_s)$, where $b_s$ and $\sigma_s$ are uniformly bounded in $\omega \in \Omega$ and $s \leq T$ and $\nu_s$ is deterministic. Suppose that there are constants $0 < A, t_0 \leq 1$ such that the support of the measure $\int \limits_0^{t_0} \nu_s(dz) ds$ is contained in the set $\lbrace z \mid |z| \leq A \rbrace$. Furthermore assume that there is a $g \in \mathcal G$ with $g(y,dz) = \nu_{yT}(dz)$ for all $y \in [0,1]$. Then for each $z \in \Ron$ and $\zeta >0$ there are $K>0$ and $0< t_1 \leq t_0 \wedge T$, which depend only on $A,z, \zeta$, the bound on $g$   in Assumption \ref{mcGDef}\eqref{mcGDef2} and the bounds on $b_s$ and $\sigma_s$, such that the transition probability is bounded by
$$\Prob(Y_t \in \mathcal I(z)) \leq K t^{\frac{|z|}{2A} - \zeta}$$
for all $0 \leq t \leq t_1$.
\end{lemma}

\noindent \textbf{Proof.}
We will only show the inequality for $z>0$ fixed, because otherwise we can consider the process $-Y$ which has the same properties. 

The H\"older inequality and the upper Burkholder-Davis-Gundy inequality yield for any $0<t \leq t_0 \wedge T$ and any $m \in \N$:
\begin{align*}
\Eb \big | \int \limits_0^t b_s ds \big |^m \leq
      t^m \Eb \bigg ( \frac{1}{t} \int \limits_0^t |b_s|^m ds \bigg ) \leq K t^m
\end{align*}
and
\begin{align*}
\Eb \big | \int \limits_0^t \sigma_s dW_s \big |^m \leq K t^{m/2} \Eb \bigg ( \frac{1}{t} \int \limits_0^t | \sigma_s|^2 ds \bigg )^{m/2}
\leq K t^{m/2}.
\end{align*}
Therefore Markov inequality ensures that the claim follows if we can show the lemma for each It\=o semimartingale with $b_s = \sigma_s \equiv 0$. 

Let $Y$ be such an It\=o semimartingale and let $0< t \leq t_0 \wedge T$. Then by Theorem II.4.15 in \cite{JacShi02} $Y$ is a process with independent increments with characteristic function
\begin{align}
\label{YtCharFkt}
  \Eb \big[ \exp \left\{ iu Y_t \right \}\big] = \exp \Big\{ \int \limits_0^t \int (e^{iuz} -1- iuz) \nu_s(dz) ds \Big \} = \exp \{ \Psi_t(iu) \}, \quad (u \in \R)
\end{align}
since $t \leq t_0$ and $A \leq 1$ with
\begin{align*}
\Psi_t(u) \defeq \int \limits_0^t \int (e^{uz} -1- uz) \nu_s(dz) ds.
\end{align*}
$\Psi_t(u)$ exists for all such $t$ and all $u \in \mathbb R$ by a Taylor expansion of the integrand and the assumption on the support of $\int_0^t \nu_s(dz) ds$ as well as item \eqref{mcGDef2} in Assumption \ref{mcGDef}. Furthermore the first two derivatives of $\Psi_t$ are given by
\begin{align*}
\Psi_t^{\prime}(u) = \int \limits_0^t \int (e^{uz} -1)z \nu_s(dz) ds; \qquad
\Psi_t^{\prime \prime}(u) = \int \limits_0^t \int z^2 e^{uz} \nu_s(dz) ds
\end{align*}
where we have exchanged differentiation and integration by the differentiation lemma of measure theory and the assumption on the support of $\int_0^t \nu_s(dz) ds$.
Without loss of generality we may assume that the measure $\int_0^t \nu_s(dz) ds$ is not zero, because otherwise
$Y_t = 0$ a.s. and the assertion of the lemma is obvious. Therefore $\Psi_t^{\prime \prime}(u) >0$ for all $u \in \R$ and $\Psi_t^{\prime}$ is a strictly increasing function with $\Psi_t^{\prime}(0)=0$ and $\lim_{u \rightarrow \infty} \Psi_t^{\prime}(u) = B \in (0,\infty]$. Thus it has a strictly increasing, differentiable inverse function $\tau_t \colon [0,B) \rightarrow \R_+$ with $\tau_t(0)=0$. Moreover, it is sufficient to show the claim for all $0< z \neq B$, because for $B$ and $\zeta >0$ we can find some $\tilde z < B$ and $0 < \tilde \zeta < \zeta$ with
$$\frac{B}{2A} - \zeta = \frac{\tilde z}{2A} - \tilde \zeta$$
and $\Prob(Y_t \geq B) \leq \Prob(Y_t \geq \tilde z)$. By Corollary 1.50 in \cite{Hof16} and by the bounded support of the measure $\int_{\scriptscriptstyle 0}^{\scriptscriptstyle t} \nu_s(dz) ds$ the right-hand side of \eqref{YtCharFkt} can be extended to an entire function on $\Cb$. Consequently, by Lemma 25.7 in \cite{Sat99} $\Eb (\exp\{c|Y_t|\}) < \infty$ for every $c>0$ and with another application of Corollary 1.50 in \cite{Hof16} the mapping $w \mapsto \Eb(\exp\{w Y_t\})$ is an entire function on $\Cb$. Therefore according to the identity theorem of complex analysis \eqref{YtCharFkt} holds for all $u \in \Cb$. Thus by Markov inequality we have for arbitrary $s \geq 0$:
\begin{align}
\label{ChebAbschEq}
\Prob( Y_t \geq z) \leq \Eb \left\{ \exp \{ s Y_t - sz \} \right\} = \exp \{ \Psi_t(s) - sz \}.
\end{align}
First suppose that $z > B$. Then we obtain
\begin{align*}
\Prob( Y_t \geq z) \leq     \limsup \limits_{s \rightarrow \infty} \exp \Big\{ \int \limits_0^s \Big( \Psi_t^{\prime}(y) - z \Big) dy \Big\} \leq \lim \limits_{s \rightarrow \infty} \exp \{ (B-z) s \} = 0
\end{align*}
and the claim obviously follows. Therefore for the rest of the proof we may assume $z < B$. In this case \eqref{ChebAbschEq} yields (recall that $\tau_t$ is the inverse function of $\Psi_t'$)
\begin{align}
\label{FiPabseqn}
\Prob( Y_t \geq z) &\leq
                     \exp \Big\{ \int \limits_0^{\tau_t(z)} \Psi_t^{\prime}(y) dy - z \tau_t(z) \Big\}
									= \exp \Big\{ \int \limits_0^z w \tau_t^{\prime}(w) dw - z \tau_t(z) \Big\} \nonumber \\
									&= \exp \Big\{ - \int \limits_0^z \tau_t(w) dw \Big\}.
\end{align}
By a Taylor expansion we have
\begin{align}
\label{PsitTaylAbsch}
(e^{Dy}-1)y \leq e^{D A} D y^2
\end{align}
for $D>0$ and $|y| \leq A$. Therefore if we set $D = \tau_t(w)$ in \eqref{PsitTaylAbsch} we obtain
\begin{align}
\label{tauAbscheqn}
w &= \Psi^{\prime}_t(\tau_t(w)) = \int \limits_0^t \int (e^{\tau_t(w)y} -1)y \nu_s(dy) ds \leq e^{A \tau_t(w)} \tau_t(w)
         \int \limits_0^t \int y^2 \nu_s(dy) ds \nonumber \\
 &\leq e^{A \tau_t(w)} \tau_t(w) Kt,
\end{align}
for arbitrary $0 \leq w < B$ and $0< t \leq t_0 \wedge T$, where the constant $K>0$ depends only on the bound on $g$    in Assumption \ref{mcGDef}\eqref{mcGDef2}. By a series expansion of the exponential function we have
\begin{align}
\label{logtauAbscheq}
\log(\tau_t(w)) \leq A \tau_t(w)
\end{align}
if $\tau_t(w) \geq \frac{2}{A^2}$ and this is the case if
\begin{align}
\label{BedantAbscheq}
\Psi_t^{\prime}\Big(\frac{2}{A^2} \Big) &= \int \limits_0^t \int \left(\exp \left\{ \frac{2}{A^2} y \right\} -1\right)y \nu_s(dy) ds \nonumber \\
       &\leq \frac{2}{A^2} e^{\frac{2}{A}} \int \limits_0^t \int y^2 \nu_s(dy)ds \leq \frac{2}{A^2} e^{\frac{2}{A}} Kt =: K_0(t) \leq w,
\end{align}
where we used \eqref{PsitTaylAbsch} again. Combining \eqref{tauAbscheqn}, \eqref{logtauAbscheq} and \eqref{BedantAbscheq} gives
\begin{align}
\label{integrabseqn}
\log \bigg (\frac{w}{t} \bigg) \leq \log(K) + 2A \tau_t(w) \Longleftrightarrow \frac{1}{2A} \log \bigg(\frac{w}{Kt} \bigg) \leq \tau_t(w)
\end{align}
for $K_0(t) \leq w < B$. Let $0 <\bar t_1 \leq t_0 \wedge T$ be small enough such that $K_0(t) \leq z <B$ for each $0 \leq t \leq \bar t_1$. Then \eqref{integrabseqn} together with \eqref{FiPabseqn} yield the estimate
\begin{align*}
\Prob(Y_t \geq z) &\leq \exp \Big \{ - \int \limits_0^z \tau_t(w) dw \Big\} \leq \exp \Big\{ - \frac{1}{2A} \int \limits_{K_0(t)}^z \log \bigg( \frac{w}{Kt} \bigg) dw \Big\} \\
 &= \exp \Big\{ - \frac{Kt}{2A} \int \limits_{K_0(t)/Kt}^{z/Kt} \log(u) du \Big\} = \exp \Big\{ -\frac{Kt}{2A}
     \Big[ -u+u\log(u) \Big]_{K_0(t)/Kt}^{z/Kt} \Big\} \\
 &= \exp \Big\{ - \frac{1}{2A} \Big[-z+z \log \Big ( \frac{z}{Kt} \Big) +K_0(t)-K_0(t) \log \Big(\frac{K_0(t)}{Kt} \Big)
    \Big] \Big\} \\
 &= \exp \Big\{ - \frac{1}{2A} \Big[-z+z \log \Big ( \frac{z}{K} \Big) +K_0(t)-K_0(t) \log \Big(\frac{K_0(t)}{K} \Big)
    \Big] \Big\} \times \\
 &\hspace{7cm}\times \exp \Big\{ \frac{z}{2A} \log(t) - \frac{K_0(t)}{2A} \log(t) \Big\} \\
 &\leq \exp \Big\{  \frac{1}{2A} \Big[z-z \log \Big ( \frac{z}{K} \Big) + 1 \Big] \Big\} t^{\frac{z}{2A} - \zeta} ~= ~K t^{\frac{z}{2A} - \zeta}
\end{align*}
for each $0 < t \leq t_1 \leq \overline t_1 \leq t_0 \wedge T$ with a $0< t_1 \leq \overline t_1$ small enough such that
$$\Big | K_0(t)-K_0(t) \log \Big(\frac{K_0(t)}{K} \Big) \Big | \leq 1$$
and $\frac{K_0(t)}{2A} < \zeta$
for every $0 < t \leq t_1$.
\hfill $\Box$

\medskip

The preceding result is helpful to deduce the following claim which is the main tool to establish consistency of $\Db_n(\kseqi, \gseqi, \taili)$ as an estimator for $D(\kseqi, \gseqi, \taili)$, when it is applied to the It\=o semimartingale $Y^{(j,n)}_s = X^{(n)}_{s+(j-1) \Delta_n} - X^{(n)}_{(j-1) \Delta_n}$. 

\begin{lemma}
\label{WahrschAbLem}
Suppose that Assumption \ref{Cond1} is  satisfied  and let  $\delta >0$. If $\Xn_0 = 0$ for all $n \in \N$, then there exist
constants $K=K(\delta)>0$ and $0 < t_0 = t_0(\delta) \leq 1$ such that
$$\Big| \Prob( \Xn_t \in \iz) - \int \limits_0^t \nu^{(n)}_s(\iz) ds \Big| \leq K t^2$$
holds for all $|z| \geq \delta$, $0 \leq t < t_0$ and $n \in \N$ with $n \Delta_n \geq 1$.
\end{lemma}

\noindent \textbf{Proof of Lemma \ref{WahrschAbLem}.}
Let $\eps < (\delta/6 \wedge 1)$ and pick a smooth cut-off function $c_\eps: \mathbb R \to \mathbb R$ satisfying
\begin{align*}
1_{[-\eps/2,\eps/2]}(u) \leq c_\eps(u) \leq 1_{[-\eps,\eps]}(u).
\end{align*}
We also define the function $\bar c_\eps$ via $\bar c_\eps(u)=1-{c_\eps}(u)$. For $n \in \N$ and $t \in \R_+$ let $M^{\neps}_t$ be the measure defined by
$$M^{\neps}_t(A)= \int \limits_0^t \int \limits_A \bar c_\eps(u) \nu^{(n)}_s(du) ds,$$
for $A \in \Bb$ which has total mass
\begin{align}
\label{totmassAbseqn}
\lambda^{\neps}_t \defeq M^{\neps}_t(\R) = \int \limits_0^t \int_{\mathbb{R}} \bar c_\eps(u) \nu^{(n)}_s(du) ds = \int \limits_0^t \int_{\mathbb{R}} \bar c_\eps(u) g\Big(\frac{s}{n\Delta_n},du\Big) ds \leq Kt,
\end{align}
where $K$ depends only on the bound on $g$ in Assumption \ref{mcGDef}\eqref{mcGDef2} and on $\eps$ and therefore on $\delta$. Furthermore, let
$$d^{(n,\eps)}_s \defeq \int u \mathtt 1_{\lbrace |u| \leq 1 \rbrace} \bar c_\eps(u) \nu^{(n)}_s(du).$$
By Theorem II.4.15 in \cite{JacShi02} for each $n \in \N$ and $t \in \R_+$ with $t \leq n \Delta_n$ we can decompose $X^{(n)}_t$ in law by
\begin{equation} \label{d1}
X^{(n)}_t =_d X_t^{(n,\eps)} + \widetilde{X}_t^{(n,\eps)},
 \end{equation}
 where $X^{(n,\eps)}$ and $\widetilde X^{(n, \eps)}$ are independent It\=o semimartingales starting in zero, with  characteristics $(b^{(n, \eps)}_s, \sigma^{(n)}_s, c_\eps(u) \nu^{(n)}_s(du))\ $,  $\ b^{(n, \eps)}_s \defeq b^{(n)}_s - d^{(n,\eps)}_s $, and $ (d^{(n,\eps)}_s,0, \bar c_\eps(u) \nu^{(n)}_s(du))$, respectively. 
 
$\widetilde X^{(n, \eps)}$ can be seen as a generalized compound Poisson process. To be precise, let $\hat \mu^{\neps}$ be a Poisson random measure independent of $X^{\neps}$ with predictable compensator $\hat \nu^{\neps}(ds,du) = \bar c_\eps(u)\nu_s^{(n)}(du) ds$ and consider the process
$$N_t^{\neps} \defeq \hat \mu^{\neps}([0,t] \times \R).
$$
By  Theorem II.4.8 in \cite{JacShi02} $N_t^{\neps} $  is a process with independent increments and distribution
$$
N_t^{\neps} -N_s^{\neps}  =_d \ \text{Poiss}\big(\lambda_t^{\neps}-\lambda_s^{\neps}\big), \quad (0 \leq s \leq t)
$$
(here we use  the convention that $\text{Poiss}(0)$ is the Dirac measure with mass in zero). Moreover, for each $n \in \N$ let $((Z^{\neps}_j(t))_{t \in [0,n \Delta_n]})_{j \in \N}$ be a sequence of independent processes, which is also independent of the Poisson random measure $\hat \mu^{\neps}$ and of the process $X^{\neps}$, such that for each $j \in \N$ and $t \in [0, n \Delta_n]$ its distribution is given by
$$ Z^{\neps}_j(t)  =_d \begin{cases}
                       M_t^{\neps}/\lambda_t^{\neps}, \quad &\text{ if } \lambda_t^{\neps} >0 \\
											 \mbox{Dirac(0)},                         \quad &\text{ if } \lambda_t^{\neps} =0.
											\end{cases}$$
Then we have for any  $n \in \N$ and $0 \leq t \leq n \Delta_n$
\begin{align}
\label{distreqn}
 \widetilde X^{(n, \eps)}_t   =_d    \sum \limits_{j=1}^{\infty} Z^{\neps}_j(t) \mathtt 1_{\left\{j \leq N_t^{\neps} \right\}},
\end{align}
because by using independence of the involved quantities we calculate the characteristic function for $w \in \R$ and $\lambda_t^{\neps} >0$ as follows:
\begin{align}
\label{CharFktBer}
\Eb &\exp \Big\{ i w \sum \limits_{j=1}^{\infty} Z^{\neps}_j(t) \mathtt 1_{\left\{j \leq N_t^{\neps} \right \} } \Big\} = \sum \limits_{j=0}^{\infty} \Eb \Big( \exp \Big\{ iw \sum \limits_{k=1}^j Z^{\neps}_k(t) \Big\} \mathtt 1_{\left\{ N_t^{\neps} = j \right\}} \Big) \nonumber \\
  &= \exp \left\{ -\lambda^{\neps}_t \right\} \sum \limits_{j=0}^{\infty} \frac{1}{j!} \big ( \Phi(M_t^{\neps})(w) \big)^j =
	   \exp \left\{ \Phi(M_t^{\neps})(w) -\lambda^{\neps}_t \right\} \nonumber \\
	&= \exp \Big\{ iw \int \limits_0^t d^{\neps}_s ds + \int \limits_0^t \int \Big( e^{iwu} -1 - iwu \mathtt 1_{\lbrace |u| \leq 1 \rbrace} \Big) \bar c_\eps(u) \nu^{(n)}_s(du)ds \Big\} \nonumber \\
	&=  \mathbb{E} [ \exp (i w \tilde X_t^{(n,\varepsilon)})].
\end{align}
In the above display $\Phi(M)$ denotes the characteristic function of a finite Borel measure $M$. The last equality in \eqref{CharFktBer} follows  from Theorem II.4.15 in \cite{JacShi02}. Furthermore note that in the case
 $\lambda_t^{\neps}=0$ the distributions in \eqref{distreqn} are obviously equal.
 
Let $z \in \Ron$ with $|z| \geq \delta$, define  $f(x)=1_{\{x \in \iz \}}$ and recall the decomposition in \eqref{d1} and the representation \eqref{distreqn}
for $t \leq n \Delta_n$. As the processes  $X^{(n,\eps)}$ and $\widetilde{X}^{(n,\eps)}$ are independent, we can calculate
\begin{align}
\label{eq:fdec}
\Eb \Big [f \Big (X^{(n)}_t \Big ) \Big] &= \sum_{j=0}^{\infty} \exp \Big \{-\lambda^{\neps}_t \Big \} \Big(\lambda_t^{\neps} \Big )^j \frac{1}{j!} \Eb \Big[f \Big(X^{(n)}_t \Big) \Big|N^{\neps}_t = j \Big] \nonumber \\
      &= \exp \Big\{ - \lambda^{\neps}_t \Big\} \Eb\Big[f\Big(X^{\neps}_t\Big)\Big]  \nonumber \\
			&+ \exp \Big\{ -\lambda^{\neps}_t \Big\} \lambda_t^{\neps} \Eb\Big[f\Big(X^{\neps}_t + Z_1^{\neps}(t)\Big)\Big] \nonumber \\
			&+ \sum_{j=2}^{\infty} \exp \Big\{- \lambda_t^{\neps} \Big\} \Big(\lambda^{\neps}_t\Big)^j \frac{1}{j!} \Eb\Big[f\Big(X^{\neps}_t+\sum_{\ell=1}^j
			   Z^{\neps}_\ell(t) \Big) \Big].
\end{align}
For the first summand on the right-hand side of the last display we use Lemma \ref{TransProbAbs} with $t_0=1, A= \eps, T= n \Delta_n$ and $\zeta = 1$ and obtain
\begin{align}
\label{FitermAbsch}
\exp \Big\{ - \lambda^{\neps}_t \Big\} \Eb\Big[f\Big(X^{\neps}_t\Big)\Big] \leq \Prob \Big( \Big| X^{\neps}_t \Big| \geq \delta \Big) \leq 2K t^{\delta/2\eps - \zeta} \leq K t^2
\end{align}
for $0 \leq t \leq \hat t_1$, where $K$ and $\hat t_1$ depend only on $\delta$, the bound
for the transition kernel $g$ in Assumption \ref{mcGDef}\eqref{mcGDef2} and the bounds on $b_s$ and $\sigma_s$. Note that $d^{\neps}_s$ is bounded for $s \leq n \Delta_n$ by a bound which depends on $\eps$, thus on $\delta$, and on the previously mentioned bound on $g$. Also, for the third term on the right-hand side of \eqref{eq:fdec}, we have
\begin{align}
\label{LasttermAbsch}
\sum_{j=2}^{\infty} \exp \Big\{- \lambda_t^{\neps} \Big\} \Big(\lambda^{\neps}_t\Big)^j \frac{1}{j!}
      \Eb\Big[f\Big(X^{\neps}_t+\sum_{\ell=1}^j Z^{\neps}_\ell(t) \Big) \Big] \leq \Big( \lambda^{\neps}_t \Big)^2 \leq Kt^2
\end{align}
by \eqref{totmassAbseqn} since $f$ is bounded by $1$. Now, if $\lambda_t^{\neps}=0$ the second term in \eqref{eq:fdec} and $\int \limits_0^t \nu_s^{(n)}(\iz) ds$ vanish. Hence the lemma follows from \eqref{FitermAbsch} and \eqref{LasttermAbsch}. Thus in the following we  assume $\lambda_t^{\neps} >0$ and  consider the term $\Eb\Big[f\Big(X^{\neps}_t + Z^{\neps}_1(t)\Big)\Big]$. For $t \leq n \Delta_n$ the distribution of $Z_1^{\neps}(t)$ has the Lebesgue density
$$u \mapsto \bar h_t^{\neps}(u) \defeq \int \limits_0^t \bar c_\eps(u) h\Big(\frac{s}{n \Delta_n}, u\Big) ds /\lambda_t^{\neps}.$$
As a consequence (for $t \leq n \Delta_n$), the function
\begin{align*}
\rho_t^{\neps}(x) := \Eb\Big[f\Big(x+ Z_1^{\neps}(t)\Big)\Big] = \Prob\Big(x + Z_1^{\neps}(t) \in \iz\Big)
\end{align*}
is twice continuously differentiable and it follows
\begin{align}
\label{Abl12Absch}
\sup \limits_{x \in \R} \left\{ \left| \Big( \rho_t^{\neps} \Big)^{\prime} (x) \right| + \left| \Big( \rho_t^{\neps} \Big)^{\prime \prime}(x) \right| \right \} \leq \frac{Kt}{\lambda_t^{\neps}},
\end{align}
where the constant $K>0$ depends only on the bound in Assumption \ref{mcGDef}\eqref{mcGDef3} for some $\eps' >0$ with $\eps' \leq \eps/2$ but not on $n$ or $t$. Using the independence of $X^{\neps}$ and $Z_1^{\neps}$ it is sufficient to discuss $\Eb[\rho_t^{\neps}(X_t^{\neps} ) ]$.   It\^o formula (Theorem I.4.57 in \cite{JacShi02}) gives
\begin{align} \label{decomp} 
	\rho_t^{\neps}(X_r^{\neps}) &= \rho_t^{\neps}(  X_0^{\neps} ) + \int_{0}^r \big(\rho_t^{\neps} \big)'(X_{s-}^{\neps}) dX_s^{\neps} \nonumber \\
	&\hspace{5mm}+ \frac 12 \int_{0}^r \big(\rho_t^{\neps} \big)''(X_{s-}^{\neps} ) d \langle X^{\neps,c}, X^{\neps,c} \rangle_s \nonumber \\
	&\hspace{5mm}+ \sum_{0<s \leq r} \Big(\big(\rho_t^{\neps} \big)(X_s^{\neps}) - \big(\rho_t^{\neps} \big)(X_{s-}^{\neps}) - \big(\rho_t^{\neps} \big)'(   X_{s-}^{\neps})
	 \Delta X_s^{\neps}) \Big),
 \end{align}
where $t \leq n \Delta_n, \ r \geq 0, \ \langle X^{\neps,c}, X^{\neps,c} \rangle_s $ denotes the predictable quadratic variation of the continuous local martingale part of $X^{\neps}$, and
$\Delta X_s^{\neps}$ is the jump size at time $s$. We now discuss each of the four summands in \eqref{decomp} separately for $r=t$: 
first, $u \in \iz$ implies $\bar c_{\eps}(u) = 1$ by definition of $\eps$. Thus, with $X_0^{\neps} = 0$
\begin{align*}
\rho_t^{\neps}\Big(X_0^{\neps}\Big) &= \Prob\Big(Z_1^{\neps}(t) \in \iz \Big)
=
\frac{1}{\lambda_t^{\neps}} \int 1_{\{ u \in \iz\}} \bar c_\eps(u) \int \limits_0^t h\Big(\frac{s}{n \Delta_n}, u\Big) ds du \\
&=
\frac{1}{\lambda_t^{\neps}} \int 1_{\{ u \in \iz\}} \int \limits_0^t g\Big(\frac{s}{n \Delta_n}, du\Big) ds
=
\frac{1}{\lambda_t^{\neps}} \int \limits_0^t \nu_s^{(n)}(\iz) ds.
\end{align*}
By the canonical representation of semimartingales (Theorem II.2.34 in \cite{JacShi02}) we get the decomposition
$$X_t^{\neps} = \int \limits_0^t b_s^{\neps} ds + Y_t^{\neps},$$
where $Y^{\neps}$ is a local martingale with characteristics $(0,\sigma_s^{(n)}, c_\eps(u) \nu^{(n)}_s(du))$ which starts at zero and has bounded jumps. Consequently $Y^{\neps}$ is a locally square integrable martingale and by Proposition I.4.50 b), Theorem I.4.52 and Theorem II.1.8 in \cite{JacShi02} its predictable quadratic variation is given by
$$\langle Y^{\neps}, Y^{\neps} \rangle_t = \int \limits_0^t (\sigma^{(n)}_s)^2 ds + \int \limits_0^t \int u^2 c_\eps(u) \nu^{(n)}_s(du) ds.$$
Thus for $t \leq n \Delta_n$ and because of the boundedness of $(\rho_t^{\neps})^{\prime}$ and the construction of the stochastic integral the integral process $((\rho_t^{\neps})'(X^{\neps}_{s-}) \cdot Y^{\neps})^t$ stopped at time $t$ is in fact a square integrable martingale because
$$\Eb \int \limits_0^t ((\rho_t^{\neps})'(X^{\neps}_{s-}))^2 d \langle Y^{\neps}, Y^{\neps} \rangle_s < \infty.$$
Therefore we obtain
$$\Eb \int \limits_0^t (\rho_t^{\neps})'(X^{\neps}_{s-}) dY^{\neps}_s =0$$
and according to \eqref{Abl12Absch} we get a bound for the second term in \eqref{decomp}:
\begin{align}
\label{intbousect}
\Big \vert \Eb\Big[\int \limits_{0}^t  (\rho_t^{\neps}  )'(X^{\neps}_{s-}) dX^{\neps}_s\Big] \Big \vert \leq \int \limits_{0}^t \Big \vert \Eb\Big[ (\rho_t^{\neps}  )'(X^{\neps}_{s-})\Big] b^{\neps}_s \Big \vert ds \leq \frac{Kt^2}{\lambda_t^{\neps}},
\end{align}
where $K>0$ depends only on $\delta$, the bounds on the characteristics and the bounds of Assumption \ref{mcGDef}\eqref{mcGDef2} and \eqref{mcGDef3} for an appropriate $\eps' >0$. For the third term in \eqref{decomp} it is immediate to get an estimate as in \eqref{intbousect}. Finally, let $\mu^{\neps}(ds,du)$ denote the random measure associated with the jumps of $X^{\neps}$ which has the predictable compensator $\nu^{\neps}(ds,du) = ds \nu^{(n)}_s(du) c_\eps(u)$. Therefore Theorem II.1.8 in \cite{JacShi02} yields for the expectation of the last term in \eqref{decomp}
\begin{align}
\label{stochintabsch}
\Eb &\Big\{ \sum_{0<s \leq t} \big(\rho_t^{\neps} (X^{\neps}_s) -  \rho_t^{\neps}  (X^{\neps}_{s-}) - (\rho_t^{\neps}  )'(X^{\neps}_{s-}) \Delta X^{\neps}_s \big) \Big\} \nonumber \\
  &=\Eb \bigg\{ \int \limits_0^t \int \Big(\rho_t^{\neps} (X^{\neps}_{s-} +u) - \rho_t^{\neps} (X^{\neps}_{s-})  \nonumber -  (\rho_t^{\neps})'(X^{\neps}_{s-}) u\Big) \mu^{\neps}(ds,du) \bigg\} \nonumber \\
	&=\Eb \bigg\{ \int \limits_0^t \int \Big( \rho_t^{\neps} (X^{\neps}_{s-} +u) - \rho_t^{\neps} (X^{\neps}_{s-}) \nonumber
	-  (\rho_t^{\neps})'(X^{\neps}_{s-}) u\Big) c_\eps(u) \nu^{(n)}_s(du) ds \bigg\} \nonumber \\
	&\leq \frac{Kt}{\lambda_t^{\neps}} \int \limits_0^t \int u^2 c_\eps(u) g\Big(\frac{s}{n \Delta_n},du\Big) ds \leq \frac{Kt^2}{\lambda_t^{\neps}}.
\end{align}
Note that the integrand in the second line in \eqref{stochintabsch} is a concatenation of a Borel measurable function on $\R^2$ and the obviously predictable function $(\omega,r,u) \mapsto (X^{\neps}_{r-}(\omega),u)$ from $\Omega \times \R_+ \times \R$ into $\R^2$. Consequently, this integrand is in fact a predictable function and Theorem II.1.8 in \cite{JacShi02} can in fact be applied. The first inequality in \eqref{stochintabsch} follows with \eqref{Abl12Absch} and a Taylor expansion of the integrand. Accordingly the constant $K$ after the last inequality in \eqref{stochintabsch} depends only on the quantities as claimed in the assertion of this lemma. Thus we have
\begin{align*}
\bigg \vert \Eb\Big[\rho_t^{\neps}\Big(X_t^{\neps}\Big)\Big]-\frac{1}{\lambda_t^{\neps}} \int \limits_0^t \nu^{(n)}_s(\iz) ds \bigg \vert \leq \frac{Kt^2}{\lambda_t^{\neps}},
\end{align*}
which together with $|1-\exp(-\lambda_t^{\neps} )| \leq Kt$ for small $t$ as well as \eqref{eq:fdec}, \eqref{FitermAbsch} and \eqref{LasttermAbsch} yields the lemma.
\hfill $\Box$

\subsection{Proof of Theorem \ref{WeakConvThm}.}
Let $X^{(n)}$ denote a semimartingale of the form \eqref{ItoSem} and consider the decomposition  $X^{(n)}_t = Y^{(n)}_t + Z^{(n)}_t$, where
\begin{align*}
Y^{(n)}_t = X^{(n)}_0 + \int \limits_0^t b^{(n)}_s ds + \int \limits_0^t \sigma^{(n)}_s dW^{(n)}_s
\end{align*} 
and $Z^{(n)}_t$ is a pure jump It\=o semimartingale with characteristics $(0,0,\nu_s^{(n)})$. Furthermore we consider the process
\begin{align*}
\Gb^\circ_n(\theta,z)
&=
\frac{1}{\sqrt{k_n}} \sum \limits_{j=1}^{\lfloor n\theta \rfloor} \{ \mathtt 1_{\lbrace \Delta_j^n Z^{(n)} \in \iz \rbrace} - \mathbb P (\Delta_j^n Z^{(n)} \in \iz) \}
\end{align*}
in $\linaeps$. The proof can be divided into two steps:
 \begin{eqnarray}\label{proof1}
 &&\Gb^\circ_n \weak \Gb \\
 \label{proof2}
 && \| \Gb_n -  \Gb^\circ_n \|_{A_\eps} \pto  0.
 \end{eqnarray}
The assertion of  Theorem \ref{WeakConvThm} then follows from Lemma 1.10.2(i) in \cite{VanWel96}.

\eqref{proof1} can be obtained with similar steps as in the first part of the proof of Theorem 2.3 in \cite{BueHofVetDet15} using Theorem 11.16 in \cite{Kos08} which is a central limit theorem for triangular arrays of row-wise i.i.d.\ data. The main difference regards the use of Lemma \ref{WahrschAbLem} which is needed as we work in general with a time-varying kernel $\nu_s^{(n)}$.  

Concerning \eqref{proof2} we have for $(\gseqi, \taili) \in A_\eps$
\begin{equation} \label{neu1}
|\Gb_n(\gseqi, \taili) - \Gb^\circ_n(\gseqi, \taili)| \leq \sqrt{k_n} \Big|U_n(\gseqi, \taili) - U^\circ_n(\gseqi, \taili)\Big| +
\sqrt{k_n} \Big| \Eb U^\circ_n(\gseqi, \taili) - \int \limits_0^{\gseqi} g(y,z) dy \Big|,
\end{equation}
where $U^\circ_n$ denotes the statistic $U_n$ based on the scheme $\{Z^{(n)}_{i \Delta_n} \mid i = 0, \ldots,n\}$. For the first term in \eqref{neu1} we obtain 
$$
\sqrt{k_n} \Big|U_n(\gseqi, \taili) - U^\circ_n(\gseqi, \taili)\Big| = o_\Prob(1),
$$
uniformly on $A_\eps$, with the same arguments as in the second part of the proof of Theorem 2.3 in \cite{BueHofVetDet15}. Furthermore, along the lines of the proof of Corollary 2.5 in \cite{BueHofVetDet15}, but using Lemma \ref{WahrschAbLem} instead, one can show that the second term in \eqref{neu1} is a uniform $o(1)$ on $A_\eps$. \hfill $\Box$

\subsection{Proof of Theorem \ref{GnVorberKonvThm}.}
Recall the decomposition $X^{(n)}_t = Y^{(n)}_t + Z^{(n)}_t$ in the proof of Theorem \ref{WeakConvThm}. The idea of the proof is to show the claim of Theorem \ref{GnVorberKonvThm} for $\hat \Gb_n^\circ$, the process being defined exactly as $\hat \Gb_n$ in \eqref{hatGbndefeq} but based on the increments $\Delta_j^n Z^{(n)}$. This can be done with Theorem 3 in \cite{Kos03}, because we have i.i.d.\ increments of the processes $Z^{(n)}$. Furthermore, by Lemma A.1 in \cite{Buc11} it is then enough to prove $\| \hat \Gb_n - \hat \Gb_n^\circ \|_{A_\eps} =o_\Prob(1)$ in order to show Theorem \ref{GnVorberKonvThm}. For a detailed proof we refer the reader to the proof of Theorem 3.3 in \cite{BueHofVetDet15}, which follows similar lines.
\hfill $\Box$

\subsection{Proof of Corollary \ref{51}.} As the process $\Hb$ is  tight, $\Hb^{(\eps)}$ is also tight
and  Theorem \ref{SchwKtomotv} together with the continuous mapping theorem yield $\Hb_n^{(\eps)} \weak \Hb^{(\eps)}$ in $\linne$. 
The assertion now follows observing the definition of $\Hb_n$ in \eqref{bbhnDef} and the fact that
$D(\kseqi, \gseqi, \taili) $ vanishes whenever $\theta \leq \gseqi^{(\eps)}_0$.
 \hfill $\Box$

\subsection{Proof of Theorem \ref{SchaeKonvThm}.}
The claim follows if we can prove the existence of a constant $K>0$ such that
\begin{align}
&\Prob \Big(\consest < \gseqi^{(\eps)}_0 \Big) = o(1), \label{zzgendeGl1}  \\
&\Prob \Big(\consest> \gseqi^{(\eps)}_0 + K \beta_n \Big)= o(1)  \label{zzgendeGl2},
\end{align}
 where $\beta_n = (\thrle_n/\sqrt{k_n})^{1/\smooi}$.
In order to verify \eqref{zzgendeGl1} we calculate as follows
\begin{align}  \label{calc}
\Prob \Big(\consest < \gseqi^{(\eps)}_0 \Big) &\leq \Prob \Big( k_n^{1/2} \Db^{(\eps)}_n(\gseqi) > \thrle_n \text{ for some }
                                          \gseqi < \gseqi^{(\eps)}_0 \Big) \\
                          &\leq \Prob \Big( k_n^{1/2} \Dc^{(\eps)}(\gseqi) + \Hb_n^{(\eps)}(\gseqi) > \thrle_n \text{ for some } \gseqi < \gseqi^{(\eps)}_0
													                 \Big) \leq \Prob \Big( \Hb_n^{(\eps)}(1) > \thrle_n \Big) = o(1), \nonumber
\end{align}
where $\Hb_n^{(\eps)}(1)$ is defined in \eqref{Fnepsdef}. The third estimate is a consequence of the fact that $\Dc^{(\eps)}(\gseqi)=0$ whenever $\gseqi < \gseqi_0^{(\eps)}$ and the final convergence follows because a weakly converging sequence in $(\R,\Bb)$ is asymptotically tight. 

For a proof of \eqref{zzgendeGl2} we  note that $k_n^{1/2} \Dc^{(\eps)}(\gseqi) - \Hb_n^{(\eps)}(\gseqi) \leq k_n^{1/2} \Db_n^{(\eps)}(\gseqi)$ and we obtain
\begin{align}
\label{FiundProbAb}
\Prob \Big(\consest > \gseqi^{(\eps)}_0 + K \beta_n \Big) &\leq \Prob \Big( k_n^{1/2} \Db_n^{(\eps)}(\gseqi) \leq \thrle_n \text{ for some } \gseqi > \gseqi_0^{(\eps)} + K \beta_n \Big) \nonumber \\
&\leq \Prob \Big ( k_n^{1/2} \Dc^{(\eps)}(\gseqi) - \Hb_n^{(\eps)}(\gseqi) \leq \thrle_n \text{ for some } \gseqi > \gseqi_0^{(\eps)} + K \beta_n \Big).
\end{align}
 Now  it follows  from \eqref{addass} that for sufficiently large $n \in \N$
\begin{align}
\label{smoountabeqn}
\inf \limits_{\gseqi \in [\gseqi_0^{(\eps)} + K \beta_n,1]} \Dc^{(\eps)}(\gseqi) = \Dc^{(\eps)}(\gseqi_0^{(\eps)} + K \beta_n) \geq \frac 12 c^{(\eps)} (K \beta_n)^\smooi.
\end{align}
Therefore with \eqref{FiundProbAb} and by the definition of $\beta_n$ we get for large $n \in \N$ and $K>0$ large enough 
\begin{align*}
\Prob \Big(\consest&> \gseqi^{(\eps)}_0 + K \beta_n \Big) \leq \Prob \Big( \frac 12 \sqrt{k_n} c^{(\eps)}
                                 (K \beta_n)^\smooi - \Hb_n^{(\eps)}(1) \leq \thrle_n \Big) \\
			&\leq \Prob \Big( \frac 12 \sqrt{k_n} c^{(\eps)} (K \beta_n)^\smooi - \Hb_n^{(\eps)}(1) \leq \thrle_n, \Hb_n^{(\eps)}(1) \leq \alpha_n \Big) + \Prob \Big( \Hb_n^{(\eps)}(1) > \alpha_n \Big) =o(1),
\end{align*}
where $\alpha_n \rightarrow \infty$ is a sequence with $\alpha_n/\thrle_n \rightarrow 0$, using asymptotic tightness again. \hfill $\Box$

\subsection{Proof of Theorem \ref{MSEzerlThm}.} 
For a proof of \eqref{MSE1queeqn} note that
\begin{align*}
\big(\consest - \gseqi^{(\eps)}_0 \big)^2 = \Big \{ \int \limits_{\gseqi_0^{(\eps)}}^1 1_{\left\{ k_n^{1/2} \Db_n^{(\eps)}(\gseqi) \leq \thrle_n \right\}} d \gseqi - \int \limits_0^{\gseqi_0^{(\eps)}} \Big(1- 1_{\left\{ k_n^{1/2} \Db_n^{(\eps)}(\gseqi) \leq \thrle_n \right\}} \Big) d\gseqi \Big\}^2
\end{align*}
and furthermore we have for any $\gseqi \in [0,1]$
\begin{align}
\label{triinekonse}
k_n^{1/2} \Dc^{(\eps)}(\gseqi) - \Hb_n^{(\eps)}(1) \leq k_n^{1/2} \Db_n^{(\eps)}(\gseqi) \leq k_n^{1/2} \Dc^{(\eps)}(\gseqi) + \Hb_n^{(\eps)}(1).
\end{align}
Thus, if $\Hb_n^{(\eps)}(1) \leq \alpha_n$, we get for sufficiently large $n \in \N$
\begin{align*}
0 &\leq \int \limits_0^{\gseqi_0^{(\eps)}} \Big(1- 1_{\left\{ k_n^{1/2} \Db_n^{(\eps)}(\gseqi) \leq \thrle_n \right\}} \Big) d \gseqi \leq \int \limits_0^{\gseqi_0^{(\eps)}} \Big(1- 1_{\left\{ k_n^{1/2} \Dc^{(\eps)}(\gseqi) + \alpha_n \leq \thrle_n \right\}} \Big) d \gseqi \nonumber \\
&= \int \limits_0^{\gseqi_0^{(\eps)}} \Big(1- 1_{\left\{ \alpha_n \leq \thrle_n \right\}} \Big) d \gseqi =0,
\end{align*}
because $\Dc^{(\eps)}(\gseqi) =0$ for $\gseqi \leq \gseqi_0^{(\eps)}$. Hence, for $n$ sufficiently large, 
\begin{align}
\label{mse1forlan}
\text{MSE}_1^{(\eps)}(\thrle_n, \alpha_n) = \Eb \bigg [ \Big\{ \int \limits_{\gseqi_0^{(\eps)}}^1 1_{\left\{ k_n^{1/2} \Db_n^{(\eps)}(\gseqi) \leq \thrle_n \right\}} d \gseqi \Big\}^2 1_{\left\{ \Hb_n^{(\eps)}(1) \leq \alpha_n \right\}} \bigg].
\end{align}
In the following let $0< \varphi < 1$ be arbitrary, let $K_1, K_2$ be as in \eqref{KeinKzweDef} and define
\begin{align}
\label{Kastdef}
K_1^\ast \defeq \Big( \frac{1-\varphi/2}{c^{(\eps)}} \Big)^{1/\smooi} \quad \text{ and } \quad K_2^\ast \defeq \Big( \frac{1+\varphi}{c^{(\eps)}} \Big)^{1/\smooi}.
\end{align}
As in \eqref{smoountabeqn} we obtain from \eqref{addass}
\begin{align}
\label{maxAbsch}
\max \limits_{\gseqi \in [\gseqi_0^{(\eps)}, \gseqi_0^{(\eps)} + K_1^\ast \beta_n]} \Dc^{(\eps)}(\gseqi) = \Dc^{(\eps)}(\gseqi_0^{(\eps)} + K_1^\ast \beta_n) \leq \frac{1}{1-\varphi /3} c^{(\eps)} (K_1^{\ast} \beta_n)^\smooi 
\end{align}
and
\begin{align}
\label{minAbsch}
\inf \limits_{\gseqi \in [\gseqi_0^{(\eps)} + K_2^\ast \beta_n,1]} \Dc^{(\eps)}(\gseqi) = \Dc^{(\eps)}(\gseqi_0^{(\eps)} + K_2^\ast \beta_n) \geq \frac 1{1+\varphi/2} c^{(\eps)} (K_2^\ast \beta_n)^\smooi
\end{align}
for $n \in \N$ large enough. Now \eqref{triinekonse} and \eqref{mse1forlan} yield
\begin{align}
\label{upboumse1}
\text{MSE}_1^{(\eps)}(\thrle_n, \alpha_n) &\leq \bigg[ \int \limits_{\gseqi_0^{(\eps)}}^1 1_{\left\{ \sqrt{k_n} \Dc^{(\eps)}(\gseqi) \leq
                               \thrle_n + \alpha_n \right\}} d\gseqi \bigg]^2 \nonumber \\
										&= \bigg[ \int \limits_{\gseqi_0^{(\eps)}}^{\gseqi_0^{(\eps)} + K_2^\ast \beta_n} 1_{\left\{ \sqrt{k_n}
										    \Dc^{(\eps)}(\gseqi) \leq \thrle_n + \alpha_n \right\}} d\gseqi \bigg]^2 \leq (K_2^\ast)^2 \beta_n^2 = K_2 \beta_n^2
\end{align}
for a sufficiently large $n \in \N$ which is the desired bound. Note that the first equation in the second line of \eqref{upboumse1} follows from \eqref{minAbsch}, because for $\gseqi \in [\gseqi_0^{(\eps)} + K_2^\ast \beta_n,1]$ we have
$$\sqrt{k_n} \Dc^{(\eps)}(\gseqi) \leq \thrle_n + \alpha_n \Longrightarrow \frac 1{1+\varphi /2} c^{(\eps)} (K_2^\ast)^\smooi \thrle_n \leq \thrle_n + \alpha_n$$ which cannot hold for large $n \in \N$ due to \eqref{Kastdef}. 

In order to get a lower bound recall \eqref{mse1forlan} and use \eqref{triinekonse} to obtain for $n \in \N$ sufficiently large
\begin{align}
\label{mse1lowbou}
\text{MSE}_1^{(\eps)}(\thrle_n, \alpha_n)  &\geq \Prob\big( \Hb_n^{(\eps)}(1) \leq \alpha_n \big) \bigg( \int
                    \limits_{\gseqi_0^{(\eps)}}^{\gseqi_0^{(\eps)} + K_1^\ast \beta_n} 1_{\left\{ \sqrt{k_n}
										    \Dc^{(\eps)}(\gseqi) \leq \thrle_n - \alpha_n \right\}} d\gseqi \bigg)^2 \nonumber \\
						 &= \Prob\big( \Hb_n^{(\eps)}(1) \leq \alpha_n \big) (K_1^\ast)^2 \beta_n^2,
\end{align}
where the equality follows from the implication (see \eqref{maxAbsch})
$$\frac{1}{1-\varphi /3} c^{(\eps)} (K_1^\ast)^\smooi \thrle_n \leq \thrle_n - \alpha_n \Longrightarrow \sqrt{k_n} \Dc^{(\eps)}(\gseqi) \leq \thrle_n - \alpha_n \quad \text{ for all } \gseqi \in [\gseqi_0^{(\eps)}, \gseqi_0^{(\eps)} + K_1^\ast \beta_n].$$
The left-hand side in the previous display always holds for large $n \in \N$ by the choice of $K_1^\ast$ in \eqref{Kastdef}. Using asymptotical tightness we also have 
$$\Prob\big( \Hb_n^{(\eps)}(1) \leq \alpha_n \big) \geq \Big( \frac{1-\varphi}{1- \varphi /2} \Big)^{2/\smooi} = K_1/(K_1^\ast)^2$$
for a large $n$, which together with \eqref{mse1lowbou} yields
$
\text{MSE}_1^{(\eps)}(\thrle_n, \alpha_n) \geq K_1 \beta_n^2.
$
\hfill $\Box$

\subsection{Proof of Theorem \ref{OptChoofThThm}.}
Similarly to (\ref{calc}) we get
\begin{align}
\label{Prochapoabeq}
\Prob\Big(\consest(\hathrleei) < \chapo \Big) \leq \Prob\Big( \Hb_n^{(\eps)}(\chapo) \geq \hathrleei \Big).
\end{align}
Recall $\Hb^{(\eps)}(\gseqi)$ from \eqref{limhbepsDef}. 
It holds that
$\Hb^{(\eps)} (\chapo) \geq |\Hb(\chapo/2,\chapo, \bar \taili)|,$
with $\bar \taili$ from \eqref{nichtdegenAss}, and by \eqref{HbGrProCov} we have
\begin{align}
\label{Vargrnull}
\operatorname{Var}(\Hb(\chapo/2,\chapo, \bar \taili)) = \frac 14 \int\limits_0^{\chapo} g(y, \bar \taili) dy >0.
\end{align}
Thus $\Hb^{(\eps)} (\chapo)$ is a supremum of a non-vanishing Gaussian process with mean zero. Due to Corollary 1.3 and Remark 4.1 in \cite{GaeMolRos07} $\Hb^{(\eps)} (\chapo)$ then has a continuous distribution function.
As a consequence \eqref{underestabEq} follows from \eqref{Prochapoabeq} and Proposition F.1 in the supplement to \cite{BueKoj14} as soon as we can show
\begin{align}
\label{zzjointwecon}
\Big(\Hb_n^{(\eps)}(\chapo), \bootsuhb_{\scriptscriptstyle n, \xi^{(1)}}(\preest), \ldots, \bootsuhb_{\scriptscriptstyle n, \xi^{(B)}}(\preest)\Big) \weak \Big(\Hb^{(\eps)} (\chapo), \Hb^{(\eps)}_{(1)} (\chapo), \ldots, \Hb^{(\eps)}_{(B)} (\chapo) \Big)
\end{align}
for any fixed $B \in \N$, where $\Hb^{(\eps)}_{(1)} (\chapo), \ldots, \Hb^{(\eps)}_{(B)} (\chapo)$ are independent copies of $\Hb^{(\eps)} (\chapo)$.

In order to establish \eqref{zzjointwecon} we first show that the sample paths of $\Hb^{(\eps)}$ are uniformly continuous on $[0,1]$ with respect to the Euclidean distance. By Theorem \ref{WeakConvThm} and Assumption \ref{mcGDef} the sample paths of the process $\Gb$ in $\linaeps$ satisfy $\Gb(0,\taili)=0$ for all $\taili \in M_\eps$ and they are uniformly continuous with respect to the Euclidean distance on $A_\eps$. Thus the uniform continuity of the sample paths of $\Hb^{(\eps)}$ holds if we can show that for a bounded and uniformly continuous function $f \colon A_\eps \rightarrow \R$ with $f(0,\taili)=0$ for all $\taili \in M_\eps$ the function $H \colon [0,1] \rightarrow \R$ defined via
$$H(\gseqi) := \sup \limits_{\gseqi ' \in [0, \gseqi]} \sup \limits_{\kseqi \in [0,\gseqi ']}\sup \limits_{\taili \in M_\eps} |f(\kseqi, \taili) - \frac \kseqi{\gseqi '} f(\gseqi ',\taili)|$$
is uniformly continuous on $[0,1]$. But since a continuous function on a compact metric space is uniformly continuous it suffices to show continuity of the function 
$$F(\gseqi) := \sup \limits_{\kseqi \in [0,\gseqi]}\sup \limits_{\taili \in M_\eps} |f(\kseqi, \taili) - \frac \kseqi{\gseqi} f(\gseqi,\taili)|$$ 
in every $\gseqi_0 \in [0,1]$. The continuity of $F$ in $\gseqi_0=0$ is obvious, because $f$ is uniformly continuous and satisfies $f(0,\taili)=0$ for all $\taili \in M_\eps$. Therefore only the case $0< \gseqi_0 \leq 1$ remains. Let $U$ be a neighbourhood of $\gseqi_0$ in $[0,1]$ which is bounded away from $0$. Then it is immediate to see that the function $h \colon B_\eps \rightarrow \R$ defined by
$$h(\kseqi,\gseqi,\taili) := f(\kseqi, \taili) - \frac \kseqi\gseqi f(\gseqi,\taili)$$
is uniformly continuous on $B_\eps \cap([0,1] \times U \times M_\eps)$. 

Let $\eta >0$ be arbitrary and choose $\delta>0$ such that 
$
|h(\kseqi_1, \gseqi_1,\taili_1) - h(\kseqi_2,\gseqi_2, \taili_2)| < \eta/2 
$
for all $(\kseqi_1, \gseqi_1,\taili_1), (\kseqi_2,\gseqi_2, \taili_2) \in B_\eps \cap([0,1] \times U \times M_\eps)$ with maximum distance $\|(\kseqi_1, \gseqi_1,\taili_1)^T - (\kseqi_2,\gseqi_2, \taili_2)^T \|_\infty \leq \delta$ and additionally such that $|\gseqi - \gseqi_0| < \delta$ implies $\gseqi \in U$. 
Then, if $|\gseqi - \gseqi_0| < \delta$, there exists $(\kseqi_1, \gseqi_0,\taili_1) \in B_\eps$ with
$
F(\gseqi_0) - \eta < |h(\kseqi_1, \gseqi_0,\taili_1)| - \eta/2 
$
and we can choose a $\kseqi_2 \leq \gseqi$ such that $\|(\kseqi_1, \gseqi_0,\taili_1)^T - (\kseqi_2,\gseqi, \taili_1)^T \|_\infty \leq \delta$ which gives
$
F(\gseqi_0) - \eta < |h(\kseqi_2,\gseqi, \taili_1)| \leq F(\gseqi).
$
In an analogous manner we see that also $F(\gseqi) - \eta < F(\gseqi_0)$ for each $\gseqi \in [0,1]$ with $|\gseqi - \gseqi_0| < \delta$, and therefore $F$ is continuous in $\gseqi_0$. 

Next we show
\begin{multline}
\label{ohneestdisteq}
\Prob \Big( \Big\| \Big(\Hb_n^{(\eps)}(\chapo), \bootsuhb_{\scriptscriptstyle n, \xi^{(1)}}(\preest), \ldots, \bootsuhb_{\scriptscriptstyle n, \xi^{(B)}}(\preest)\Big)^T - \\ 
- \Big(\Hb_n^{(\eps)}(\chapo), \bootsuhb_{\scriptscriptstyle n, \xi^{(1)}}(\chapo), \ldots, \bootsuhb_{\scriptscriptstyle n, \xi^{(B)}}(\chapo)\Big)^T \Big\|_{\infty} > \eta \Big) \rightarrow 0,
\end{multline}
for arbitrary $\eta >0$. By Proposition 10.7 in \cite{Kos08} and Theorem \ref{BootstrTeststThm} we have $\bootsuhb_{\scriptscriptstyle n, \xi^{(i)}} \weakP$ $ \Hb^{(\eps)}$ in $\linne$ for all $i=1, \ldots, B$, which yields $\bootsuhb_{\scriptscriptstyle n, \xi^{(i)}} \weak \Hb^{(\eps)}$ for all $i=1, \ldots, B$ with the same reasoning as in the proof of Theorem 2.9.6 in \cite{VanWel96}. Theorem 1.5.7 and its addendum therein show that $\bootsuhb_{\scriptscriptstyle n, \xi^{(i)}}$ is asymptotically uniformly $\rho$-equicontinuous in probability for each $i$, where $\rho$ denotes the Euclidean metric on $[0,1]$ because the sample paths of $\Hb^{(\eps)}$ are uniformly continuous with respect to $\rho$ and $([0,1], \rho)$ is totally bounded. 

Therefore, for any $\gamma >0$ we can choose a $\delta >0$ such that
\begin{align*}
\max \limits_{i=1, \ldots, B} \limsup \limits_{n \rightarrow \infty} \Prob\Big( \sup \limits_{\rho(\gseqi_1, \gseqi_2) < \delta} \Big| \bootsuhb_{\scriptscriptstyle n, \xi^{(i)}}(\gseqi_1) - \bootsuhb_{\scriptscriptstyle n, \xi^{(i)}}(\gseqi_2) \Big| > \eta \Big) < \gamma/(2B),
\end{align*}
which yields
\begin{align*}
\Prob &\Big( \Big\| \Big(\Hb_n^{(\eps)}(\chapo), \bootsuhb_{\scriptscriptstyle n, \xi^{(1)}}(\preest), \ldots, \bootsuhb_{\scriptscriptstyle n, \xi^{(B)}}(\preest)\Big)^T - \\
&\hspace{6cm} - \Big(\Hb_n^{(\eps)}(\chapo), \bootsuhb_{\scriptscriptstyle n, \xi^{(1)}}(\chapo), \ldots, \bootsuhb_{\scriptscriptstyle n, \xi^{(B)}}(\chapo)\Big)^T \Big\|_{\infty} > \eta \Big) \\
&\hspace{5mm} \leq \Prob \Big( \Big| \bootsuhb_{\scriptscriptstyle n, \xi^{(i)}}(\preest) - \bootsuhb_{\scriptscriptstyle n, \xi^{(i)}}(\chapo) \Big| > \eta \text{ for at least one } i= 1, \ldots, B \text{ and } |\preest - \chapo| < \delta \Big) + \\
&\hspace{6cm}+ \Prob\Big(|\preest - \chapo| \geq \delta \Big) \\
&\hspace{5mm} \leq \Prob\Big(|\preest - \chapo| \geq \delta \Big) + \sum \limits_{i=1}^B \Prob\Big( \sup \limits_{\rho(\gseqi_1, \gseqi_2) < \delta} \Big| \bootsuhb_{\scriptscriptstyle n, \xi^{(i)}}(\gseqi_1) - \bootsuhb_{\scriptscriptstyle n, \xi^{(i)}}(\gseqi_2) \Big| > \eta \Big) < \gamma
\end{align*}
for $n \in \N$ large enough, using consistency of the preliminary estimator. 

Thus, now that we have established \eqref{ohneestdisteq}, by Lemma 1.10.2(i) in \cite{VanWel96} we obtain \eqref{zzjointwecon} if we can show 
\begin{align*}
\Big(\Hb_n^{(\eps)}(\chapo), \bootsuhb_{\scriptscriptstyle n, \xi^{(1)}}(\chapo), \ldots, \bootsuhb_{\scriptscriptstyle n, \xi^{(B)}}(\chapo)\Big) \weak \Big(\Hb^{(\eps)} (\chapo), \Hb^{(\eps)}_{(1)} (\chapo), \ldots, \Hb^{(\eps)}_{(B)} (\chapo) \Big).
\end{align*}
But this is an immediate consequence of the continuous mapping theorem and
\begin{align}
\label{JointConvEq}
(\Gb_n, \hat \Gb_{n, \xi^{(1)}}, \ldots, \hat \Gb_{n, \xi^{(B)}}) &\weak  (\mathbb G,  \mathbb G^{(1)}, \ldots, \mathbb G^{(B)})
\end{align}
in $(\linaeps)^{B+1}$ for all $B \in \N$, where $\mathbb G^{(1)}, \ldots, \mathbb G^{(B)}$ are independent copies of $\mathbb G$, since $\Hb_n^{(\eps)}(\chapo)$, $\bootsuhb_{\scriptscriptstyle n, \xi^{(i)}}(\chapo)$ and $\Hb^{(\eps)} (\chapo)$ are the images of the same continuous functional applied to $\Gb_n$, $\hat \Gb_{\scriptscriptstyle n, \xi^{(i)}}$ and $\Gb$, respectively. \eqref{JointConvEq} follows as in Proposition 6.2 in \cite{BueHofVetDet15}. \hfill $\Box$

\medskip

\subsection{Proof of Theorem \ref{KorThrChoiThm}.}

We start with a proof of
$\beta_{n} \probto 0$
which is equivalent to $\hacothrle/\sqrt{k_n} \probto 0$. Therefore we have to show
\begin{align}
\label{empgrlevaleq}
\Prob(\hacothrle/\sqrt{k_n} \leq x) = \Prob \Big( \frac{1}{B_n} \sum \limits_{i=1}^{B_n} \mathtt 1_{\{\bootsuhb_{\scriptscriptstyle n, \xi^{(i)}}(\preest)  \leq (\sqrt{k_n} x)^{1/r} \}} \geq 1- \alpha_n \Big) \rightarrow 1,
\end{align}
for arbitrary $x>0$, by the definition of $\hacothrle$ in \eqref{thrledefeq}. Since the
\begin{align*}
\mathtt 1_{\{\bootsuhb_{\scriptscriptstyle n, \xi^{(i)}}(\preest)  \leq (\sqrt{k_n} x)^{1/r} \}} - \Prob_\xi\Big( \bootsuhb_{n}(\preest)  \leq (\sqrt{k_n} x)^{1/r}   \Big), \quad i=1, \ldots, B_n,
\end{align*}
are pairwise uncorrelated with mean zero and bounded by $1$, we have
\begin{align}
\label{DistMuExpEq}
\Prob &\Big( \Big| \frac 1{B_n} \sum \limits_{i=1}^{B_n} \mathtt 1_{\{\bootsuhb_{\scriptscriptstyle n, \xi^{(i)}}(\preest)  \leq (\sqrt{k_n} x)^{1/r} \}} - \Prob_\xi\Big( \bootsuhb_{n}(\preest)  \leq (\sqrt{k_n} x)^{1/r} \Big) \Big| > \alpha_n/2 \Big) 
\leq 4 \alpha_n^{-2} B_n^{-1} \rightarrow 0.
\end{align}
Therefore, in order to prove \eqref{empgrlevaleq}, it suffices to verify
\begin{align}
\label{secProbabeq}
&\Prob \Big( \Prob_\xi\Big( \bootsuhb_{n}(\preest)  \leq (\sqrt{k_n} x)^{1/r} \Big)  < 1- \alpha_n/2 \Big) \nonumber 
\leq \frac 2{\alpha_n} \Prob \Big( \bootsuhb_{n}(\preest)  > (\sqrt{k_n} x)^{1/r} \Big) \\ \leq& \frac 2{\alpha_n} \Prob \Big( 2 \sup \limits_{\gseqi \in [0,1]} \sup \limits_{|\taili| \geq \eps} |\hat \Gb_{n}(\gseqi, \taili)|  > (\sqrt{k_n} x)^{1/r} \Big)
\rightarrow 0,
\end{align}
where the first inequality in the above display follows with the Markov inequality and the last inequality in \eqref{secProbabeq} is a consequence of the fact that
$\bootsuhb_{n}(\preest) \le \bootsuhb_{n}(1) \leq 2 \sup \limits_{\gseqi \in [0,1]} \sup \limits_{|\taili| \geq \eps} |\hat \Gb_{n}(\gseqi, \taili)|.$
Furthermore, by the definition of $\hat \Gb_{n}$ in \eqref{hatGbndefeq} we have
\begin{align}
\label{EbsupAbscheq}
\Eb \Big\{\sup \limits_{\gseqi \in [0,1]} \sup \limits_{|\taili| \geq \eps} |\hat \Gb_{n}(\gseqi, \taili)| \Big\}\leq \frac 1{n \sqrt{k_n}} \sum \limits_{j=1}^n \sum \limits_{i=1}^n (\Prob(| \Delta_j^n X^{(n)}| \geq \eps) + \Prob(| \Delta_i^n X^{(n)}| \geq \eps)),
\end{align}
because of $\Eb| \xi_j| \leq 1$ for every $j=1, \ldots, n$. Recall the decomposition $X^{(n)} = Y^{(n)} + Z^{(n)}$ in the proof of Theorem \ref{WeakConvThm} and let $v_n = \Delta_n^{\tau/2} \rightarrow 0$ with $\tau$ from Assumption \ref{Assump1}. Then we have for $i = 1, \ldots, n$ and $n \in \N$ large enough
\begin{align}
\label{FiProbgeAbeq}
\Prob (|\Delta_i^n X^{(n)}| \geq \eps) \leq \Prob (|\Delta_i^n Y^{(n)}| \geq v_n) + \Prob (|\Delta_i^n Z^{(n)}| \geq \eps/2) \leq \Prob (|\Delta_i^n Y^{(n)}| \geq v_n) + K \Delta_n,
\end{align}
where the last inequality follows using Lemma \ref{WahrschAbLem}. By H\"older inequality, the Burkholder-Davis-Gundy inequalities (see for instance page 39 in \citealp{JacPro12}) and the Fubini theorem we have with $p>2$, $1< \alpha <3$ from Assumption \ref{Assump1}, for each $1 \leq j \leq n$,
\begin{align}
\label{driftwahrabeq}
\Eb \Big| \int_{(j-1) \Delta_n}^{j \Delta_n} b^{(n)}_s ds \Big|^{\alpha} &\leq \Delta_n^{\alpha} \Eb \Big( \frac{1}{\Delta_n} \int_{(j-1) \Delta_n}^{j \Delta_n} |b^{(n)}_s|^{\alpha} ds \Big) \leq K \Delta_n^{\alpha}
\end{align}
and
\begin{align}
\label{volatwahrabeq}
\Eb \Big| \int_{(j-1) \Delta_n}^{j \Delta_n} \sigma^{(n)}_s dW^{(n)}_s \Big|^p &\leq K \Delta_n^{p/2} \Eb \Big( \frac{1}{\Delta_n} \int_{(j-1) \Delta_n}^{j \Delta_n} | \sigma^{(n)}_s|^2 ds \Big)^{p/2} \nonumber
\\ &\leq K \Delta_n^{p/2} \Eb \Big( \frac{1}{\Delta_n} \int_{(j-1) \Delta_n}^{j \Delta_n} | \sigma^{(n)}_s|^{p} ds \Big) \leq K \Delta_n^{p/2}.
\end{align}
Together with \eqref{FiProbgeAbeq} and the Markov inequality these estimates yield
\begin{align}
\label{FinalProAbeq}
\Prob (|\Delta_i^n X^{(n)}| \geq \eps) \leq K \Delta_n^{p/2 - p\tau/2} + K \Delta_n^{\alpha - \alpha \tau/2} + K \Delta_n = K \Delta_n^{\frac{2p+p}{2p +2}} + K \Delta_n \leq K \Delta_n.
\end{align}
Therefore due to \eqref{secProbabeq}, \eqref{EbsupAbscheq}, \eqref{FinalProAbeq} and the Markov inequality we obtain
\begin{align*}
\Prob \Big( \Prob_\xi\Big( \bootsuhb_{n}(\preest)  \leq (\sqrt{k_n} x)^{1/r} \Big)  < 1- \alpha_n/2 \Big)  \leq K \frac{n^2 \Delta_n}{\alpha_n n \sqrt{k_n} (\sqrt{k_n})^{1/r}} = K \Big( (n \Delta_n)^{\frac{1-r}{2r}} \alpha_n \Big)^{-1} \rightarrow 0,
\end{align*}
by the assumptions on the involved sequences. Thus we conclude $\beta_n \probto 0$. 

Next we show $\hacothrleor \probto \infty$, which is equivalent to
\begin{align*}
\Prob(\hacothrle \leq x) = \Prob \Big( \frac{1}{B_n} \sum \limits_{i=1}^{B_n} \mathtt 1_{\{\bootsuhb_{\scriptscriptstyle n, \xi^{(i)}}(\preest)  \leq x^{1/r} \}} \geq 1- \alpha_n \Big) \rightarrow 0,
\end{align*}
for each $x >0$. With the same considerations as for \eqref{DistMuExpEq} it is sufficient to show
$$\Prob \Big( \Prob_\xi\Big( \bootsuhb_{n}(\preest) > x^{1/r} \Big)  \leq 2 \alpha_n \Big) \rightarrow 0.$$
By continuity of the function $\kseqi \mapsto \int_0^\kseqi g(y, \bar \taili) dy$ for $\bar \taili$ from \eqref{nichtdegenAss} we can find $\bar \kseqi < \bar \gseqi < \chapo$ with
\begin{align}
\label{ChoibakseEq}
\int \limits_0^{\bar \kseqi} g(y, \bar \taili) dy >0
\end{align}
and because of 
$$\hat \Hb_{n}(\bar \kseqi, \bar \gseqi, \bar \taili) \leq \bootsuhb_{n}(\preest) \Longrightarrow \Prob_\xi\Big( \hat\Hb_{n}(\bar \kseqi, \bar \gseqi, \bar \taili) > x^{1/r} \Big) \leq \Prob_\xi\Big( \bootsuhb_{n}(\preest) > x^{1/r} \Big)$$
on the set $\{\bar \gseqi < \preest\}$ and the consistency of the preliminary estimate it further suffices to prove
\begin{align}
\label{GenzzUngl}
\Prob \Big( \Prob_\xi\Big( \bootsuhb_{n}(\preest) > x^{1/r} \Big)  \leq 2 \alpha_n \text{ and } \bar \gseqi < \preest \Big) \leq \Prob \Big( \Prob_\xi\Big( \hat \Hb_{n}(\bar \kseqi, \bar \gseqi, \bar \taili) > x^{1/r} \Big)  \leq 2 \alpha_n \Big) \rightarrow 0.
\end{align}
In order to show \eqref{GenzzUngl} we want to use a Berry-Esseen type result. Recall 
\begin{align*}
\hat \Hb_{n} (\bar \kseqi, \bar \gseqi, \bar \taili) = \frac 1{\sqrt{n \Delta_n}} \sum \limits_{j=1}^n B_j \xi_j
\end{align*}
from \eqref{hatHbndefeq} with
$
B_j = \Big(\mathtt 1_{\{j \leq \ip{n \bar\kseqi} \}} - \frac{\bar\kseqi}{\bar\gseqi} \mathtt 1_{\{j \leq \ip{n \bar\gseqi}\}}\Big) A_j,
$
where
\begin{align*}
A_j = \mathtt 1_{\{\Delta_j^n X^{(n)} \in \izba\}} - \frac 1n \sum \limits_{i=1}^n \mathtt 1_{\{\Delta_i^n X^{(n)} \in \izba\}}.
\end{align*}
By the assumptions on the multiplier sequence it is immediate to see that
\begin{align*}
\bar W_n^2 := \Eb_\xi (\hat \Hb_{n} (\bar \kseqi, \bar \gseqi, \bar \taili))^2 = \frac 1{n\Delta_n} \sum \limits_{j=1}^n B_j^2.
\end{align*}
Thus Theorem 2.1 in \cite{CheSha01} yields
\begin{align}
\label{BerEssBoundEq}
\sup \limits_{x \in \R} \Big| \Prob_\xi\Big( \hat \Hb_{n}(\bar \kseqi, \bar \gseqi, \bar \taili) > x \Big) - (1- \Phi(x/\bar W_n)) \Big| \leq K \Big\{ \sum \limits_{i=1}^n \Eb_\xi U_i^2 \mathtt 1_{\{|U_i| >1\}} + \sum \limits_{i=1}^n \Eb_\xi |U_i|^3 \mathtt 1_{\{|U_i| \leq 1\}} \Big\},
\end{align}
with
$
U_i = \frac{B_i \xi_i}{\sqrt{n \Delta_n} \bar W_n}
$
and where $\Phi$ denotes the standard normal distribution function. Before we proceed further in the proof of \eqref{GenzzUngl}, we first show
\begin{align}
\label{OpeinszzEq}
\frac 1{\bar W_n^2} = \frac{n \Delta_n}{\sum \limits_{j=1}^n B_j^2} = O_\Prob (1),
\end{align}
which is
\begin{align*}
\lim \limits_{M \rightarrow \infty} \limsup \limits_{n \rightarrow \infty} \Prob \Big( n \Delta_n > M \sum \limits_{j=1}^n B_j^2 \Big) =0.
\end{align*}
Let $M>0$. Then a straightforward calculation gives
\begin{align}
\label{ProbKlAbschEq}
\Prob \Big( n \Delta_n &> M \sum \limits_{j=1}^n B_j^2 \Big) \leq 
			 \Prob \Big( n \Delta_n > M' \sum \limits_{j=1}^{\ip{n \bar \kseqi}} A_j^2 \Big) \nonumber \\
			&= \Prob \Big( n \Delta_n > M' \frac 1{n^2} \sum \limits_{j=1}^{\ip{n \bar \kseqi}} \sum \limits_{i=1}^n \sum \limits_{k=1}^n \Big( \mathtt 1_{\{\Delta_j^n X^{(n)} \in \izba \}} + \mathtt 1_{\{\Delta_i^n X^{(n)} \in \izba \}} \mathtt 1_{\{\Delta_k^n X^{(n)} \in \izba \}}              \nonumber \\
			&\hspace{35mm} - \mathtt 1_{\{\Delta_i^n X^{(n)} \in \izba \}} \mathtt 1_{\{\Delta_j^n X^{(n)} \in \izba \}} - \mathtt 1_{\{\Delta_j^n X^{(n)} \in \izba \}} \mathtt 1_{\{\Delta_k^n X^{(n)} \in \izba \}} \Big) \Big),
\end{align}
with $M' = M(1- \bar \kseqi /\bar \gseqi)^2$. Now consider again the decomposition $X^{(n)} = Y^{(n)} + Z^{(n)}$ of the underlying It\=o semimartingale as in the proof of Theorem \ref{WeakConvThm} and the sequence $v_n = \Delta_n^{\tau/2} \rightarrow 0$ with $\tau$ from Assumption \ref{Assump1}.
With \eqref{driftwahrabeq}, \eqref{volatwahrabeq} and Lemma \ref{WahrschAbLem} it is immediate to see that for $1 \leq i,k \leq n$
\begin{align*}
\Eb \Big| \mathtt 1_{\{\Delta_i^n X^{(n)} \in \izba \}} \mathtt 1_{\{\Delta_k^n X^{(n)} \in \izba \}} - \mathtt 1_{\{\Delta_i^n Z^{(n)} \in \izba \}} \mathtt 1_{\{\Delta_k^n Z^{(n)} \in \izba \}} \Big| 
= o(\Delta_n).
\end{align*}
Setting
\begin{align*}
D_n &:= \frac 1{n^3 \Delta_n} \sum \limits_{j=1}^{\ip{n \bar \kseqi}} \sum \limits_{i=1}^n \sum \limits_{k=1}^n \Big( \mathtt 1_{\{\Delta_i^n X^{(n)} \in \izba \}} \mathtt 1_{\{\Delta_k^n X^{(n)} \in \izba \}} - \mathtt 1_{\{\Delta_i^n X^{(n)} \in \izba \}} \mathtt 1_{\{\Delta_j^n X^{(n)} \in \izba \}} -\\
&\hspace{9cm} - \mathtt 1_{\{\Delta_j^n X^{(n)} \in \izba \}} \mathtt 1_{\{\Delta_k^n X^{(n)} \in \izba \}}  \Big)
\end{align*}
it is easy to deduce
$\Eb |D_n|
= o(1),
$
using Lemma \ref{WahrschAbLem} again as well as independence of the increments of $Z^{(n)}$. Combining this result with \eqref{ProbKlAbschEq} we have 
\begin{align*}
\Prob \Big( n \Delta_n &> M \sum \limits_{j=1}^n B_j^2 \Big) \leq \\
&\leq\Prob(|D_n| > 1/M') + \Prob \Big( 1/M' > D_n + \frac 1{n \Delta_n} \sum \limits_{j=1}^{\ip{n \bar \kseqi}} \mathtt 1_{\{\Delta_j^n X^{(n)} \in \izba \}} \text{ and } |D_n| \leq 1/M' \Big) \\
&\leq \Prob(|D_n| > 1/M') + \Prob \Big( 2/M' > \frac 1{n \Delta_n} \sum \limits_{j=1}^{\ip{n \bar \kseqi}} \mathtt 1_{\{\Delta_j^n X^{(n)} \in \izba \}} \Big)
\end{align*}
for $M>0$. Thus with \eqref{ChoibakseEq} we obtain \eqref{OpeinszzEq}, because
by Theorem \ref{WeakConvThm} we have
$$\frac 1{n \Delta_n} \sum \limits_{j=1}^{\ip{n \bar \kseqi}} \mathtt 1_{\{\Delta_j^n X^{(n)} \in \izba \}} = \int \limits_0^{\bar \kseqi} g(y,\bar \taili) dy + o_\Prob(1).$$
Recall that our main objective is to show \eqref{GenzzUngl} and thus we consider the Berry-Esseen bound on the right-hand side of \eqref{BerEssBoundEq}. For the first summand we distinguish two cases according to the assumptions on the multiplier sequence.

Let us discuss the case of bounded multipliers first. For $M>0$ we have
$$|U_i| \leq \frac {\sqrt{M} K}{\sqrt{n\Delta_n}}$$
for all $i= 1, \ldots,n$ on the set $\{1/\bar W_n^2 \leq M\}$, since $|B_i|$ is bounded by $1$. As a consequence
\begin{align}
\label{vanishEq}
\sum \limits_{i=1}^n \Eb_\xi U_i^2 \mathtt 1_{\{|U_i| >1\}} =0
\end{align}
for large $n \in \N$ on the set $\{1/\bar W_n^2 \leq M\}$. 

In the situation of normal multipliers, recall that there exist constants $K_1, K_2 >0$ such that for $\xi \sim \mathcal N(0,1)$ and $y >0$ large enough we have
\begin{align}
\label{GaussmomAbsch}
\Eb_\xi \xi^2 \mathtt 1_{\{|\xi| >y\}} = \frac{2}{\sqrt{2 \pi}} \int \limits_y^\infty z^2 e^{-z^2/2} dz \le K \Prob(\mathcal N(0,2) >y) \leq K_1 \exp(-K_2 y^2).
\end{align}
Thus we can calculate for $n \in \N$ large enough on the set $\{1/\bar W_n^2 \leq M\}$
\begin{align*}
\sum \limits_{i=1}^n \Eb_\xi U_i^2 \mathtt 1_{\{|U_i| >1\}} &= \sum \limits_{i=1}^n \Big( \sum \limits_{j=1}^n B_j^2 \Big)^{-1} B_i^2 \Eb_\xi \xi_i^2 \mathtt 1_{\{|\xi_i| > (\sum \limits_{j=1}^n B_j^2)^{1/2}/|B_i| \}} \\
&\leq \sum \limits_{i=1}^n \Big( \sum \limits_{j=1}^n B_j^2 \Big)^{-1} \Eb_\xi \xi_i^2 \mathtt 1_{\{|\xi_i| > (\sum \limits_{j=1}^n B_j^2)^{1/2}\}} \\
&\leq \frac{M}{n \Delta_n} \sum \limits_{i=1}^n \Eb_\xi \xi_i^2 \mathtt 1_{\{|\xi_i| > (n\Delta_n /M)^{1/2}\}} \leq \frac{K_1}{\Delta_n} \exp(-K_2 n \Delta_n),
\end{align*}
where $K_1$ and $K_2$ depend on $M$. The first inequality in the above display uses $|B_i| \leq 1$ again and the last one follows with \eqref{GaussmomAbsch}. Now let $\rho >0$ with $n \Delta_n^{1+\rho} \rightarrow \infty$ and define $\bar p := 1/\rho$. Then, for $n \geq N(M) \in \N$ on the set $\{1/\bar W_n^2 \leq M\}$, using $\exp(-K_2 n \Delta_n) \leq (n \Delta_n)^{-\bar p}$, we conclude
\begin{align}
\label{FirstTerklnuf}
\sum \limits_{i=1}^n \Eb_\xi U_i^2 \mathtt 1_{\{|U_i| >1\}} \leq K_1 \Delta_n^{-1} (n \Delta_n)^{-\bar p} = K_1 (n \Delta_n^{1+ \rho})^{-\bar p}.
\end{align}

We now consider the second term on the right-hand side of \eqref{BerEssBoundEq}, for which 
\begin{align*}
\sum \limits_{i=1}^n \Eb_\xi |U_i|^3 \mathtt 1_{\{|U_i| \leq 1\}} &\leq \sum \limits_{i=1}^n \Big(\sum \limits_{j=1}^n B_j^2\Big)^{-3/2} |B_i|^3 \Eb_\xi |\xi_i|^3 \leq \frac K{(n \Delta_n)^{3/2}} \sum \limits_{i=1}^n |B_i|
\end{align*}
holds on $\{1/\bar W_n^2 \leq M\}$, using $|B_i| \leq 1$ again. With \eqref{FinalProAbeq} we see that
\begin{align*}
\Eb\Big( \sum \limits_{i=1}^n |B_i|\Big) \leq \Eb \Big(\sum \limits_{i=1}^n |A_i| \Big)\leq 2n \max \limits_{i=1, \ldots,n} \Prob (|\Delta_i^n X^{(n)}| \geq \eps) \leq K n \Delta_n.
\end{align*}
Consequently,
\begin{align}
\label{SecProbNufEq}
\Prob\Big( 1/\bar W_n^2 \leq M  \text{ and } K \sum \limits_{i=1}^n \Eb_\xi |U_i|^3 \mathtt 1_{\{|U_i| \leq 1\}} > (n\Delta_n)^{-1/4} \Big) &\leq \Prob\Big( \frac K{(n \Delta_n)^{3/2}} \sum \limits_{i=1}^n |B_i|> (n\Delta_n)^{-1/4} \Big) \nonumber \\
&\leq K (n\Delta_n)^{-1/4}
\end{align}
follows. Thus from \eqref{vanishEq}, \eqref{FirstTerklnuf} and \eqref{SecProbNufEq} we see that with $K>0$ from \eqref{BerEssBoundEq} for each $M>0$ there exists a $K_3>0$ such that
\begin{multline}
\label{BoundnufEq}
\Prob\Big( 1/\bar W_n^2 \leq M \text{ and } K \Big\{ \sum \limits_{i=1}^n \Eb_\xi U_i^2 \mathtt 1_{\{|U_i| >1\}} + \sum \limits_{i=1}^n \Eb_\xi |U_i|^3 \mathtt 1_{\{|U_i| \leq 1\}} \Big\} \\
> K_3((n\Delta_n)^{-1/4} + (n \Delta_n^{1+ \rho})^{-\bar p}) \Big) \rightarrow 0.
\end{multline}

Now we can show \eqref{GenzzUngl}. Let $\eta >0$ and according to \eqref{OpeinszzEq} choose an $M >0$ with
$\Prob(1/\bar W_n^2 >M) < \eta/2$
for all $n \in \N$. For this $M>0$ choose a $K_3 >0$ such that the probability in \eqref{BoundnufEq} is smaller than $\eta/2$ for large $n$. Then for $n \in \N$ large enough we have
\begin{multline*}
\Prob \Big( \Prob_\xi\Big( \hat \Hb_{n}(\bar \kseqi, \bar \gseqi, \bar \taili) > x^{1/r} \Big)  \leq 2 \alpha_n \Big) < \\
\Prob\Big( (1- \Phi(x^{1/r}/\bar W_n)) \leq 2 \alpha_n + K_3((n\Delta_n)^{-1/4} + (n \Delta_n^{1+ \rho})^{-\bar p}) \text{ and } 1/\bar W_n^2 \leq M \Big) + \eta = \eta,
\end{multline*}
using \eqref{BerEssBoundEq} and the fact, that if $1/\bar W_n^2 \leq M$ there exists a $c>0$ with $ (1- \Phi(x^{1/r}/\bar W_n)) >c$. 

Thus we have shown $\hacothrleor \probto \infty$ and we are only left with proving \eqref{overestabEq}. Let 
$$K =\Big(  (1+\varphi)/c^{(\eps)} \Big)^{1/\smooi} > \Big(  1/c^{(\eps)} \Big)^{1/\smooi}$$
 for some $\varphi >0$. Then 
\begin{align*}  
\Prob\Big(\consest(\hacothrleor) &> \chapo + K \beta_n \Big) \leq \Prob \Big( \sqrt{n \Delta_n} \Db_n^{(\eps)}(\gseqi) \leq \hacothrleor \text{ for some } \gseqi > \chapo + K \beta_n \Big) \\
&\leq \Prob\Big( \sqrt{n \Delta_n} \Dc^{(\eps)}(\gseqi) - \Hb^{(\eps)}_n(1) \leq \hacothrleor \text{ for some } \gseqi > \chapo + K \beta_n \Big).
\end{align*}
By \eqref{addass} there exists a $y_0 >0$ with
$$\inf \limits_{\gseqi \in [\chapo + K y_1,1]} \Dc^{(\eps)}(\gseqi) = \Dc^{(\eps)}(\chapo + K y_1) \geq (c^{(\eps)}/(1+\varphi/2)) (K y_1)^\smooi$$
for all $0 \leq y_1 \leq y_0$. Distinguishing the cases $\{\beta_n > y_0\}$ and $\{\beta_n \le y_0\}$ we get due to $\beta_n \probto 0$
\begin{align*}
&\Prob\Big(\consest(\hacothrleor) > \chapo + K \beta_n \Big) \\
   \leq& \Prob\Big( \sqrt{n\Delta_n}  (c^{(\eps)}/(1+\varphi/2)) (K \beta_n)^\smooi - \Hb^{(\eps)}_n(1) \leq \hacothrleor \Big) + o(1) \leq P^{(1)}_n + P^{(2)}_n + o(1)
\end{align*}
with 
\begin{align*}
P^{(1)}_n &= \Prob\Big( \sqrt{n\Delta_n}  (c^{(\eps)}/(1+\varphi/2)) (K \beta_n)^\smooi - \Hb^{(\eps)}_n(1) \leq \hacothrleor \text{ and } \Hb^{(\eps)}_n(1) \leq b_n \Big), \\
P^{(2)}_n &= \Prob\Big(\Hb^{(\eps)}_n(1) > b_n \Big),
\end{align*}
where $b_n := \sqrt{\hacothrleor}$. Due to the choice $K = \Big(  (1+\varphi)/c^{(\eps)} \Big)^{1/\smooi}$ and the definition of $\beta_n$ it is clear that $P^{(1)}_n = o(1)$, because $\hacothrleor \probto \infty$. 

Concerning $P^{(2)}_n$ let $H_n$ be the distribution function of $\Hb^{(\eps)}_n(1)$ and let $H$ be the distribution function of $\Hb^{(\eps)}(1)$. Then as we have seen in \eqref{Vargrnull} in the proof of Theorem \ref{OptChoofThThm} the function $H$ is continuous and by Theorem \ref{SchwKtomotv} and the continuous mapping theorem $H_n$ converges pointwise to $H$. Thus for $\eta >0$ choose an $x>0$ with $1-H(x) < \eta/2$ and conclude
\begin{align*}
P^{(2)}_n \leq \Prob(b_n \leq x) + 1- H_n(x) \leq \Prob(b_n \leq x) + 1- H(x) + | H_n(x) - H(x) | < \eta,
\end{align*}
for $n \in \N$ large enough, because of $\hacothrleor \probto \infty$.
\hfill $\Box$

\medskip

\subsection{Proof of Proposition \ref{CorrGCP}.}
Under the null hypothesis ${\bf H}_0(\eps)$ we have 
$k_n^{1/2} \Db_n^{(\eps)}(1) = \Hb_n^{(\eps)}(1)$.
Furthermore,
\begin{align*}
\text{Var}(\Hb(\bar \kseqi, 1, \bar \taili)) 
  &= \int \limits_0^{\bar \kseqi} g(y, \bar \taili) dy -2 \bar \kseqi \int \limits_0^{\bar \kseqi} g(y,\bar \taili) dy + \bar{\kseqi}^2 \int \limits_0^1 g(y, \bar \taili) dy \\
	&= (1- \bar \kseqi)^2 \int \limits_0^{\bar \kseqi} g(y, \bar \taili) dy + \bar{\kseqi}^2 \int \limits_{\bar \kseqi}^1 g(y, \bar \taili) dy >0.
\end{align*}
Therefore, as in the proof of Theorem \ref{OptChoofThThm}, $\Hb^{(\eps)}(1)$ has a continuous cdf and 
\begin{align*}
\Big(\Hb_n^{(\eps)}(1), \bootsuhb_{\scriptscriptstyle n, \xi^{(1)}}(1), \ldots, \bootsuhb_{\scriptscriptstyle n, \xi^{(B)}}(1)\Big) \weak \Big(\Hb^{(\eps)} (1), \Hb^{(\eps)}_{(1)} (1), \ldots, \Hb^{(\eps)}_{(B)} (1) \Big).
\end{align*}
holds in $(\R^{B+1}, \Bb^{B+1})$ for every $B \in \N$, where $\Hb^{(\eps)}_{(1)} (1), \ldots, \Hb^{(\eps)}_{(B)} (1)$ are independent copies of $\Hb^{(\eps)} (1)$. As a consequence the assertion follows with Proposition F.1 in the supplement to \cite{BueKoj14}. The result for the test \eqref{testlokal} follows in the same way.
\hfill $\Box$

\medskip

\subsection{Proof of Proposition \ref{CorConsi}.} 
If ${\bf H}_1(\eps)$  holds, then \eqref{neuA}  is a simple consequence of 
$\lim \limits_{n \rightarrow \infty} \mathbb P(k_n^{1/2} \Db_n^{(\eps)}(1) \geq K) =1$ for all $K>0$ and $\lim \limits_{K \rightarrow \infty} \limsup \limits_{n \rightarrow \infty} \mathbb P \Big( \hat{\Hb}_{n, \xi^{(b)}}^{(\eps)}(1) > K \Big) = 0$ which follow from Theorem \ref{SchwKtomotv} and Theorem \ref{BootstrTeststThm} by similar arguments as in the previous proofs. The second claim can be shown in the same way.
%
\qed

\bigskip

\subsection{Proof of the results in  Example \ref{ex2} and Example \ref{ex3}(\ref{ex3nr2}).} \label{proofex2}
\begin{enumerate}[(1)]
\item First we show that a kernel as in \eqref{locbestagDef} belongs to the set $\Gc$. Using the uniqueness theorem for measures we see that $g(y,dz)$ is the measure with Lebesgue density $h_y(z) = A(y)\beta(y)/|z|^{1+\beta(y)}$ for each $y \in [0,1]$ and $z \neq 0$. This function is continuously differentiable with derivative $h_y'(z) = -\text{sgn}(z)A(y)\beta(y)(1+\beta(y))/|z|^{2+\beta(y)}$, and we obtain
$$\sup \limits_{y \in [0,1]} \sup \limits_{|z| \ge \eps} \Big( h_y(z) + |h'_y(z)| \Big) < \infty$$
for any $\eps >0$ so that Assumption \ref{mcGDef}\eqref{mcGDef3} is satisfied. Assumption \ref{mcGDef}\eqref{mcGDef20} is obvious, 
and by definition it is also clear that $g(y,dz)$ does not charge $\{0\}$ for any $y \in [0,1]$, Finally, a simple calculation using symmetry of the integrand yields
\begin{align*}
\sup \limits_{y \in [0,1]} \Big (\int (1 \wedge z^2) g(y,dz) \Big) &= \sup \limits_{y \in [0,1]} \Big ( 2A(y)\beta(y) \Big\{  \int \limits_0^1  z^{1-\beta(y)} dz + \int \limits_1^\infty z^{-1-\beta(y)} dz \Big\} \Big) \\
&= \sup \limits_{y \in [0,1]} \Big ( 2A(y)\beta(y) \Big\{ \frac 1{2-\beta(y)} + \frac 1{\beta(y)}  \Big\} \Big) < \infty,
\end{align*}
by the assumptions on $A$ and $\beta$. Thus also Assumption \ref{mcGDef}\eqref{mcGDef2} is valid.
\item Now we show that if additionally \eqref{erstkonstGl} and \eqref{analyticass} are satisfied, both $k_{0,\eps}< \infty$ holds for every $\eps >0$ and $g_k(\taili)$ is a bounded function on $M_\eps$ as stated in Example \ref{ex3}\eqref{ex3nr2}. By shrinking the interval $U$ if necessary we may assume without loss of generality that the functions $\bar A, \bar\beta$ are bounded away from $0$ on $U$. But then it is well known from complex analysis that there exist a domain $U \subset U^\ast \subset \Cb$ and holomorphic functions $A^\ast, \beta^\ast \colon U^\ast \rightarrow \Cb^+ := \{ u \in \Cb \mid \operatorname{Re}(u) >0 \}$ such that $\bar A, \bar \beta$ are the restrictions of $A^\ast$ and $\beta^\ast$ to $U$. Therefore for any $\taili \in \Ron$ the function $g^\ast(y,z) = A^\ast(y) \exp\{-\beta^\ast(y) \log(|z|)\}$ is holomorphic in $y \in U^\ast$ as a concatenation of holomorphic functions and thus its restriction $\bar g(y,\taili)$ to $y \in U$ is real analytic. Consequently, by shrinking $U$ again if necessary, we have the power series expansion
\begin{align}
\label{bargexpan}
\bar g(y,\taili) = \sum \limits_{k=0}^\infty \frac{g_k(z)}{k!} (y- \gseqi_0)^k,
\end{align}
for every $y \in U$ and $\taili \in \Ron$. If $k_{0,\eps} = \infty$ for some $\eps >0$, then for any $k \in \N$ and $\taili \in M_\eps$ we have $g_k(z)=0$. Thus we obtain
\begin{align*}
\bar g(y,\taili) = g_0(\taili) \Longleftrightarrow \log(\bar A(y)) = \log(g_0(\taili)) + \bar \beta(y) \log(|\taili|)
\end{align*}
for every $y \in U$ and $\taili \in M_\eps$. Taking the derivative with respect to $y$ for a fixed $\taili$ yields
$
(\log(\bar A(y)))' = \log(|\taili|) \bar \beta'(y)
$
for each $y \in U$ and $\taili \in M_\eps$. But since it is assumed that at least one of the functions $\bar A$ and $\bar \beta$ is non-constant, there is a $y_0 \in U$ such that one derivative (and therefore both) are different from zero. Varying $\taili$ for this $y_0 \in U$ yields a contradiction. 

In order to show that for each $k \in \N$ the function $g_k(\taili)$ is bounded in $\taili \in M_\eps$, we use
\begin{align*}
\frac{\partial^\ell}{\partial y^\ell} \Big(\frac{1}{|z|^{\bar \beta(y)}}\Big)(\gseqi_0) = \frac{\partial^\ell}{\partial y^\ell} (\exp(-\bar\beta(y) \log(|z|)))(\gseqi_0) = (-1)^\ell \frac{\log^\ell(|z|)}{|z|^{\beta_0}} \beta_1^\ell,
\end{align*}
for $\ell \in \N$, where $\bar \beta(y) = \beta_0 + \beta_1(y - \gseqi_0)$. Furthermore let the series expansion of the real analytic function $\bar A$ be given by 
\begin{equation}
\label{barAExpan}
\bar A(y) = \sum_{\ell=0}^{\infty} A_\ell (y-\gseqi_0)^\ell, \quad y \in U.
\end{equation}
Then, using the generalization of the product formula for higher derivatives 
\begin{align}
\label{gkAbschEq2}
\sup \limits_{z \in M_\eps} |g_k(z)| &= \sup \limits_{z \in M_\eps} \Big|\sum \limits_{\ell=0}^k \binom{k}{\ell} \bar A^{(\ell)}(\gseqi_0) \frac{\partial^{(k-\ell)}}{\partial y^{(k-\ell)}} \Big(\frac{1}{|z|^{\bar \beta(y)}}\Big)(\gseqi_0) \Big| \nonumber \\
&\leq  \sum \limits_{\ell=0}^k \binom{k}{\ell} |A_\ell| \ell ! (K(k-\ell))^{k-\ell} 
\end{align}
follows. The inequality in the display above holds because for $\ell\in\N$ the continuously differentiable function $f_\ell(z) = \log^\ell (z) /z^{\beta_0}$ on $(\eps,\infty)$ satisfies $\lim_{z \rightarrow \infty} f_\ell(z) =0$ and the only possible roots of its derivative are $z=1$ and $z= \exp\{ \ell/\beta_0\}$. Therefore we obtain
\begin{align*}
\sup \limits_{z \in M_\eps} \frac{|\log^{\ell}(|z|)|}{|z|^{\beta_0}} = \max \Big\{\frac{|\log^\ell(\eps)|}{\eps^{\beta_0}}, \Big(\frac{\ell}{\beta_0}\Big)^\ell e^{-\ell} \Big\} \leq (K \ell)^{\ell},
\end{align*}
for some suitable $K>0$ which does not depend on $\ell$. $\bar A$ is real analytic, thus the power series in \eqref{barAExpan} has a positive radius of convergence and by the Cauchy-Hadamard formula this fact is equivalent to the existence of a $K>0$ such that $|A_\ell| \leq K^{\ell}$ for each $\ell \in \N$. As a consequence \eqref{gkAbschEq2} yields
\begin{align}
\label{gkAbschEq}
\sup \limits_{z \in M_\eps} |g_k(z)| &\leq  (Kk)^k + \sum \limits_{\ell=1}^k \binom{k}{\ell} (K\ell)^\ell (K(k-\ell))^{k-\ell}  \nonumber \\
&\leq  (Kk)^k  \sum \limits_{\ell=0}^k \binom{k}{\ell} \leq (Kk)^k.
\end{align}

\item Finally, we show the expansion \eqref{smoexpanGl}, and 
we prove first that it suffices to verify it for $\tilde \Dc^{(\eps)}$ from \eqref{tildDcepsDef}. If $\tilde\Dc^{(\eps)}$ satisfies \eqref{smoexpanGl} with a function $\aleph(\gseqi)$, we have for $\gseqi \in [\gseqi_0, \gseqi_0 + \lambda]$
\begin{align*}
c^{(\eps)}&(\gseqi - \gseqi_0)^{k_{0,\eps}+1} - K(\gseqi - \gseqi_0)^{k_{0,\eps}+2} \leq \sup \limits_{\gseqi_0 \leq \gseqi ' \leq \gseqi} c^{(\eps)}(\gseqi ' - \gseqi_0)^{k_{0,\eps}+1} - \sup \limits_{\gseqi_0 \leq \gseqi ' \leq \gseqi} |\aleph(\gseqi')| \\
&\leq \sup \limits_{\gseqi_0 \leq \gseqi ' \leq \gseqi} |c^{(\eps)}(\gseqi ' - \gseqi_0)^{k_{0,\eps}+1} + \aleph(\gseqi')| = \sup \limits_{\gseqi_0 \leq \gseqi ' \leq \gseqi} \tilde\Dc^{(\eps)}(\gseqi') = \Dc^{(\eps)}(\gseqi) \\
&\leq \sup \limits_{\gseqi_0 \leq \gseqi ' \leq \gseqi} c^{(\eps)}(\gseqi ' - \gseqi_0)^{k_{0,\eps}+1} + \sup \limits_{\gseqi_0 \leq \gseqi ' \leq \gseqi} |\aleph(\gseqi')| \leq c^{(\eps)}(\gseqi - \gseqi_0)^{k_{0,\eps}+1} + K(\gseqi - \gseqi_0)^{k_{0,\eps}+2},
\end{align*}
where 
$$c^{(\eps)} = \frac{1}{(k_{0,\eps} +1)!} \sup \limits_{|z| \geq \eps}|g_{k_{0,\eps}}(z)|.$$

			A power series can be integrated term by term within its radius of convergence. Therefore, \eqref{bargexpan} gives for $\gseqi_0 \leq \gseqi \in U$
			\begin{align*}
			\int \limits_{\gseqi_0}^\gseqi g(y,\taili) dy = \sum \limits_{k=0}^\infty \frac{g_k(z)}{(k+1)!} (\gseqi - \gseqi_0)^{k+1},
			\end{align*}
			which yields 
			\begin{align}
			\label{Dvereinf}
			D(\kseqi,\gseqi,\taili) &= 
				 \begin{cases}
				   -\frac \kseqi\gseqi \sum \limits_{k=1}^\infty \frac{g_k(z)}{(k+1)!}(\gseqi - \gseqi_0)^{k+1}, \quad & \text{ if }\kseqi \leq \gseqi_0 \\
					\frac 1\gseqi \sum \limits_{k=1}^\infty \frac{g_k(z)}{(k+1)!} [\gseqi(\kseqi - \gseqi_0)^{k+1} - \kseqi(\gseqi - \gseqi_0)^{k+1}], \quad & \text{ if }\gseqi_0 < \kseqi \leq \gseqi.
					\end{cases}
			\end{align}
For any set $T$ and bounded functions $g,h \colon T \rightarrow \R$ we have
$
| \sup_{t \in T} |g(t)| - \sup_{t \in T} |h(t)|| \leq \sup_{t \in T} |g(t) - h(t)|.
$			
Together with \eqref{Dvereinf} this yields for $\gseqi_0 < \gseqi \in U$
\begin{align}
\label{FiFueTerm}
\sup \limits_{\kseqi \in [0, \gseqi_0]} \sup \limits_{z \in M_\eps} |D(\kseqi,\gseqi,\taili)| = \sup \limits_{\kseqi \in [0, \gseqi_0]} \sup \limits_{z \in M_\eps} \Big|\frac \kseqi\gseqi \frac{g_{k_{0,\eps}}(z)}{(k_{0,\eps}+1)!}(\gseqi - \gseqi_0)^{k_{0,\eps}+1} \Big| + O((\gseqi - \gseqi_0)^{k_{0, \eps} +2}),
\end{align}
and
\begin{multline}
\label{SecFueTerm}
\sup \limits_{\kseqi \in (\gseqi_0, \gseqi]} \sup \limits_{z \in M_\eps} |D(\kseqi,\gseqi,\taili)| =	\sup \limits_{\kseqi \in (\gseqi_0, \gseqi]} \sup \limits_{z \in M_\eps}	\Big| \frac 1\gseqi \frac{g_{k_{0,\eps}}(z)}{(k_{0,\eps}+1)!} [\gseqi(\kseqi - \gseqi_0)^{k_{0,\eps}+1} - \kseqi(\gseqi - \gseqi_0)^{k_{0,\eps}+1}] \Big| +\\+ O((\gseqi - \gseqi_0)^{k_{0, \eps} +2}),
\end{multline}
for $\gseqi \downarrow \gseqi_0$, as soon as we can show
\begin{align}
\label{FiRestAbsch}
\sup \limits_{\kseqi \in (\gseqi_0, \gseqi]} \sup \limits_{z \in M_\eps} \Big| \frac 1\gseqi \sum \limits_{k=k_{0,\eps} +1}^\infty \frac{g_k(z)}{(k+1)!} [\gseqi(\kseqi - \gseqi_0)^{k+1} - \kseqi(\gseqi - \gseqi_0)^{k+1}] \Big| =O((\gseqi - \gseqi_0)^{k_{0,\eps}+2})
\end{align}
and
\begin{align}
\label{SecRestAbsch}
\sup \limits_{\kseqi \in [0,\gseqi_0]} \sup \limits_{z \in M_\eps} \Big| \frac \kseqi\gseqi \sum \limits_{k=k_{0,\eps}+1}^\infty \frac{g_k(z)}{(k+1)!}(\gseqi - \gseqi_0)^{k+1} \Big| =O((\gseqi - \gseqi_0)^{k_{0,\eps}+2}).
\end{align}
To prove \eqref{FiRestAbsch} note that for $k \in \N$ and $\gseqi_0 < \kseqi \leq \gseqi$
\begin{align*}
&|\gseqi(\kseqi - \gseqi_0)^{k+1} - \kseqi (\gseqi - \gseqi_0)^{k+1}| \leq \gseqi |[(\kseqi-\gseqi_0)^{k+1} - (\gseqi - \gseqi_0)^{k+1}]| + (\gseqi - \kseqi)(\gseqi - \gseqi_0)^{k+1} \\
  =& \gseqi \Big| \sum \limits_{j=0}^k (\gseqi - \gseqi_0)^{k-j}(\kseqi - \gseqi)(\kseqi - \gseqi_0)^j \Big| + (\gseqi - \kseqi)(\gseqi - \gseqi_0)^{k+1} \leq 2\gseqi (k+1)(\gseqi-\gseqi_0)^{k+1}
\end{align*}
holds, and by \eqref{gkAbschEq} we have 
$
\limsup \limits_{k \rightarrow \infty} ({\bar g_k}/{k!})^{1/k} < \infty$  for $\bar g_k = \sup \limits_{z \in M_\eps} |g_k(z)|.
$
Consequently, 
\begin{align*}
&\sup \limits_{\kseqi \in (\gseqi_0, \gseqi]} \sup \limits_{z \in M_\eps} \Big| \frac 1\gseqi \sum \limits_{k=k_{0,\eps} +1}^\infty \frac{g_k(z)}{(k+1)!} [\gseqi(\kseqi - \gseqi_0)^{k+1} - \kseqi(\gseqi - \gseqi_0)^{k+1}] \Big|
\leq \sum \limits_{k=k_{0,\eps} +1}^\infty \frac{2\bar g_k}{k!} (\gseqi - \gseqi_0)^{k+1} \\&= (\gseqi - \gseqi_0)^{k_{0,\eps}+2} \sum \limits_{k=0}^\infty \frac{2\bar g_{k+k_{0,\eps}+1}}{(k+k_{0,\eps}+1)!} (\gseqi - \gseqi_0)^{k} 
= O((\gseqi - \gseqi_0)^{k_{0,\eps}+2})
\end{align*}
for $\gseqi \downarrow \gseqi_0$, because the latter power series has a positive radius of convergence around $\gseqi_0$. 

For the same reason, in order to prove \eqref{SecRestAbsch} we use  
\begin{align*}
\sup \limits_{\kseqi \in [0,\gseqi_0]} \sup \limits_{z \in M_\eps} \Big| \frac \kseqi\gseqi \sum \limits_{k=k_{0,\eps}+1}^\infty \frac{g_k(z)}{(k+1)!}(\gseqi - \gseqi_0)^{k+1} \Big| &\leq \frac{\gseqi_0}\gseqi (\gseqi - \gseqi_0)^{k_{0,\eps}+2} \sum \limits_{k=0}^\infty \frac{\bar g_{k+k_{0,\eps}+1}}{(k+k_{0,\eps}+2)!} (\gseqi - \gseqi_0)^{k} \nonumber \\
  &=O((\gseqi - \gseqi_0)^{k_{0,\eps}+2})
\end{align*}
as $\gseqi \downarrow \gseqi_0$.
Now, because of 
$$\tilde\Dc^{(\eps)}(\gseqi) = \max \Big\{\sup \limits_{\kseqi \in [0, \gseqi_0]} \sup \limits_{z \in M_\eps} |D(\kseqi,\gseqi,\taili)|, \sup \limits_{\kseqi \in (\gseqi_0, \gseqi]} \sup \limits_{z \in M_\eps} |D(\kseqi,\gseqi,\taili)| \Big\}$$
\eqref{FiFueTerm} and \eqref{SecFueTerm} yield the desired expansion \eqref{smoexpanGl} for $\tilde \Dc^{(\eps)}$, if we can show for $\gseqi \downarrow \gseqi_0$
\begin{align*}
& \sup \limits_{\kseqi \in [0, \gseqi_0]} \sup \limits_{z \in M_\eps} \Big|\frac \kseqi\gseqi \frac{g_{k_{0,\eps}}(z)}{(k_{0,\eps}+1)!}(\gseqi - \gseqi_0)^{k_{0,\eps}+1} \Big| = \frac{\bar g_{k_{0, \eps}}}{(k_{0,\eps} +1)!} (\gseqi - \gseqi_0)^{k_{0,\eps}+1} + O((\gseqi - \gseqi_0)^{k_{0,\eps}+2}), \\
&\sup \limits_{\kseqi \in (\gseqi_0, \gseqi]} \sup \limits_{z \in M_\eps}	\Big| \frac 1\gseqi \frac{g_{k_{0,\eps}}(z)}{(k_{0,\eps}+1)!} [\gseqi(\kseqi - \gseqi_0)^{k_{0,\eps}+1} - \kseqi(\gseqi - \gseqi_0)^{k_{0,\eps}+1}] \Big| \\
&\hspace{65mm}= \frac{\bar g_{k_{0, \eps}}}{(k_{0,\eps} +1)!} (\gseqi - \gseqi_0)^{k_{0,\eps}+1} + O((\gseqi - \gseqi_0)^{k_{0,\eps}+2}).
\end{align*}
The first assertion is obvious, since 
$| \gseqi_0/\gseqi -1 | \leq K(\gseqi - \gseqi_0)$
by $0< \gseqi_0 <\gseqi$. 
In order to prove the latter claim, consider for $0< \gseqi_0 < \gseqi < 1$ and $k \in \N$ the function $h_k:  [\gseqi_0, \gseqi] \rightarrow \R$ with 
$h_k(\kseqi) = \gseqi(\kseqi - \gseqi_0)^{k+1} - \kseqi (\gseqi - \gseqi_0)^{k+1}.
$
Its derivative is given by 
$h'_k (\kseqi) = \gseqi(k+1) (\kseqi-\gseqi_0)^k - (\gseqi - \gseqi_0)^{k+1}$
and it has a unique root at $\kseqi_0$ with
$$\gseqi_0 < \kseqi_0 = \gseqi_0 + \frac{(\gseqi - \gseqi_0)^{1+1/k}}{(\gseqi(k+1))^{1/k}} < \gseqi.$$
Thus because of $h_k(\gseqi_0) <0$, $h_k(\gseqi)=0$, $h'_k(\kseqi) <0$ for $\kseqi < \kseqi_0$ and $h'_k(\kseqi) >0$ for $\kseqi > \kseqi_0$ we obtain the result, since for $\gseqi \downarrow \gseqi_0$
\begin{align*}
\frac{1}\gseqi &\sup \limits_{\kseqi \in (\gseqi_0, \gseqi]} |h_{k_{0,\eps}}(\kseqi)| = \frac{1}\gseqi |h_{k_{0,\eps}}(\kseqi_0)| \nonumber \\
&= \Big| \frac{\gseqi_0}{\gseqi} + \frac{(\gseqi - \gseqi_0)^{1+1/k_{0,\eps}}}{\gseqi^{1+1/k_{0,\eps}}(k_{0,\eps}+1)^{1/k_{0,\eps}}} - \frac{(\gseqi - \gseqi_0)^{1+1/k_{0,\eps}}}{\gseqi^{1+1/k_{0,\eps}}(k_{0,\eps}+1)^{1+1/k_{0,\eps}}} \Big| (\gseqi - \gseqi_0)^{k_{0,\eps}+1} \nonumber \\
&= \frac{\gseqi_0}\gseqi  (\gseqi - \gseqi_0)^{k_{0,\eps}+1} + O( (\gseqi - \gseqi_0)^{k_{0,\eps}+2+1/k_{0,\eps}}) = (\gseqi - \gseqi_0)^{k_{0,\eps}+1} + O( (\gseqi - \gseqi_0)^{k_{0,\eps}+2}).
\end{align*}
\hfill $\Box$
\end{enumerate}

\end{document}